\documentclass[11pt,reqno]{amsart}%

\usepackage{amsthm,amsmath,amsfonts,amssymb,amsxtra,appendix,bookmark,euscript,delarray,dsfont,mathrsfs}
\usepackage[utf8]{inputenc}
\usepackage{bbm}
\usepackage{hyperref}
\usepackage{soul}

\allowdisplaybreaks

\newtheorem{theorem}{Theorem}[section]
\newtheorem{proposition}[theorem]{Proposition}
\newtheorem{lemma}[theorem]{Lemma}

\theoremstyle{definition}

\newtheorem{definition}[theorem]{Definition}

\theoremstyle{remark}

\newtheorem{remark}[theorem]{Remark}

\numberwithin{equation}{section}

\newcommand{\1}{\mathbbm{1}}

\newcommand{\CC}{\mathcal{C}}

\renewcommand{\epsilon}{\varepsilon}

\newcommand{\N}{\mathbb{N}}

\renewcommand{\phi}{\varphi}
\newcommand{\R}{\mathbb{R}}

\newcommand{\Sph}{\mathbb{S}}

\newcommand{\Z}{\mathbb{Z}}
\newcommand{\bP}{\mathbb{P}}

\newcommand{\cE}{{\mathcal{E}}}	
\newcommand{\cK}{{\mathcal{K}}}	
\newcommand{\lin}{{\rm lin}}
	
\newcommand{\andf}{\quad\hbox{and}\quad}
\newcommand{\with}{\quad\hbox{with}\quad}
\newcommand{\al}{\alpha}
\newcommand{\wt}{\widetilde}
\newcommand{\sE}{\mathscr{E}}

\newcommand{\sY}{\mathscr{Y}}

\DeclareMathOperator{\dist}{dist}
\DeclareMathOperator{\supp}{supp}

\title[Quantitative bounds]{Quantitative bounds for bounded solutions \\ of 3-D Navier-Stokes equations \\ in the endpoint critical Besov spaces}
\author[R. Hu]{Ruilin Hu}
\address[R. Hu]{Department of Mathematical Sciences, University of Bath, Bath, BA2 7AY, UK.}
\email{rh2488@bath.ac.uk}

\author[P.-T. Nguyen]{Phuoc-Tai Nguyen}
\address[P.-T. Nguyen]{Department of Mathematics and Statistics, Masaryk University, Brno, Czech Republic}
\email{ptnguyen@math.muni.cz}

\author[Q-.H. Nguyen]{Quoc-Hung Nguyen}
\address[Q-.H. Nguyen]{State Key Laboratory of Mathematical Sciences, Academy of Mathematics and Systems Science, Chinese Academy of Sciences, Beijing, 100190, China}
\email{qhnguyen@amss.ac.cn}

\author[P. Zhang]{Ping Zhang}
\address[P. Zhang]{State Key Laboratory of Mathematical Sciences, Academy of Mathematics $\&$ Systems Science, The Chinese Academy of Sciences, Beijing 100190, China;\\
School of Mathematical Sciences, University of Chinese Academy of Sciences,
Beijing 100049, China}
\email{zp@amss.ac.cn}

\begin{document}

	\begin{abstract}
Building on Tao's quantitative regularity theory and triple-logarithmic blow-up estimate in $L^3$ in \cite{Tao_20}, we consider classical solutions $(u,P)$ of the three-dimensional incompressible Navier--Stokes equations on $[0,T]\times\mathbb{R}^3$. For $3<p<\infty$, under simultaneous uniform control of the two scaling-critical quantities $\|u\|_{L_T^\infty(\dot B_{p,\infty}^{-1+\frac{3}{p}})}$ and $\||D|^{-1+\frac{3}{p}}u\|_{L_T^\infty(L^p)}$, we obtain explicit quantitative estimates for all spatial derivatives of $u$. As a consequence, we derive a mixed blow-up criterion coupling a double exponential of the critical Besov norm with the $L^p$ norm of $|D|^{-1+\frac{3}{p}}u$, which forces quantified growth of at least one of these two critical quantities near any finite blow-up time. The low regularity and lack of dyadic summability in the endpoint Besov space are handled through a finite iterative decomposition that successively improves spatial integrability and produces an energy-class remainder, together with refined nonlinear energy estimates. The nonlocal signed quantity $|D|^{-1+\frac{3}{p}}u$ is treated by localized mean-zero vector tests and almost orthogonality across geometrically separated concentration scales.
    \medskip
		
		\noindent\textit{Keywords: quantitative regularity, blowup criteria, critical Besov spaces.}
		
		\medskip
		
		\noindent\textit{Mathematics Subject Classification: 35Q30, 37N10, 76D05.}
		
	\end{abstract}
	
	\maketitle
	\tableofcontents
	
	\section{Introduction}
	\subsection{Overview on the relevant literature}
In this paper, we consider the following three-dimensional incompressible Navier–Stokes equations  (with the viscosity normalized to be $1$):
\begin{equation} \label{prob:NS}
\left\{
\begin{aligned}
\partial_t u - \Delta u + u \cdot \nabla u &= - \nabla P, \quad (t,x)\in\R^+\times\R^3,\\
\nabla \cdot u &= 0,
\end{aligned}
\right.
\end{equation}
where $u = (u^1, u^2, u^3) : [0,T] \times \mathbb{R}^3 \to \mathbb{R}^3$ and $P : [0,T] \times \mathbb{R}^3 \to \mathbb{R}$ denote the velocity vector and the scalar pressure function of the viscous fluid, respectively. Here $u \cdot \nabla = \sum_{j=1}^3 u^j \partial_j$ and $\nabla \cdot u = \sum_{j=1}^3 \partial_j u^j$.

The Navier–Stokes equations model the motion of viscous incompressible fluids and are a central subject in the field of mathematical fluid mechanics. The system consists of momentum equations, describing the time evolution of the velocity field $u$ under the influence of diffusion, convection, and pressure gradient, along with an incompressibility condition $\nabla \cdot u = 0$ that enforces volume preservation. This setting in three-dimensional space is known for its mathematical complexity, particularly regarding the existence and regularity of solutions over time.

We shall restrict our attention to \textit{classical} solutions, namely solutions that are smooth and such that all derivatives of $u, P$ are in $L_T^\infty (L^2)$. We say that $u$ blows up at a finite time $T_*$ if $u$ admits no smooth extension to $[0,T_*]$. $T_*$ is referred to as the \textit{maximal time of existence} of $u$.

A conspicuous feature of \eqref{prob:NS} lies in its scaling symmetry
\begin{equation} \label{scaling}
u(t,x) \mapsto u_\lambda(t,x) := \lambda u(\lambda^2 t, \lambda x), \qquad
P(t,x) \mapsto P_\lambda(t,x) := \lambda^2 P(\lambda^2 t, \lambda x), \quad \forall \lambda > 0,
\end{equation}
which implies in particular that the following spaces
\begin{equation} \label{critical spaces}
L^3(\mathbb{R}^3) \hookrightarrow L^{3,q}(\mathbb{R}^3) \hookrightarrow \dot{B}_{p,q}^{-1+\frac{3}{p}}(\mathbb{R}^3) \hookrightarrow \dot{B}_{p,\infty}^{-1+\frac{3}{p}}(\mathbb{R}^3),
\end{equation}
with $3 < p, q < \infty$, are critical to \eqref{prob:NS} in the sense that the norm $\|u\|_{L^\infty(\mathbb{R}^+; X)}$ is invariant under the scaling transformation \eqref{scaling}, where $X$ is any of the  functional spaces listed in \eqref{critical spaces}.

The research on blowup criteria for problem \eqref{prob:NS} in terms of critical norms is a long-standing challenging topic in nonlinear PDEs from both qualitative and quantitative perspectives, initiated by the landmark achievements of Leray \cite{Ler_34} and Prodi–Serrin–Ladyzhenskaya \cite{Pro_59, Ser_62, Lad_67} for non-endpoint critical Lebesgue spaces (namely for $L_t^q(L_x^p)$ with $3 < p \leq \infty$ and $2 \leq q < \infty$ such that $\frac{2}{q} + \frac{3}{p} = 1$). The endpoint case $L_t^\infty (L_x^3)$ remained open for a long time and was settled in a celebrated work by Escauriaza, Seregin and Šverák \cite{EscSerSve_03}. By using a contradiction argument based on a blowup procedure and unique continuation results, they asserted that if a solution $u$ of \eqref{prob:NS} blows up at a finite time $T_*$, then
\begin{equation} \label{limsup-L3}
\limsup_{t \nearrow T_*} \|u(t)\|_{L^3} = +\infty.
\end{equation}

The research topic has recently taken a noticeable leap forward, evidenced by numerous impactful contributions in various contexts and different directions; see, e.g., \cite{Alb-2018, GalIfPla_03, GalKocPla_13, GalKocPla_16, KenKoc_11, Phuc_15, Ser_12}. Among many others, we refer to Seregin \cite{Ser_12}, where the superior limit in \eqref{limsup-L3} was replaced by the limit; to Phuc \cite{Phuc_15} for the generalization of \eqref{limsup-L3} to non-endpoint borderline Lorentz spaces $L^{3,r}$ for any $r < \infty$; and to Gallagher, Koch and Planchon \cite{GalKocPla_16} and Albritton \cite{Alb-2018} for respectively extending \eqref{limsup-L3} to the superior limit and limit in non-endpoint critical Besov spaces $\dot{B}_{p,q}^{-1+\frac{3}{p}}$ with $3 < p, q < \infty$. The endpoint case has appeared to be more challenging and was treated by Barker in \cite{Bar_17}, where a qualitative blowup criterion regarding $B_{p,\infty}^{-1+\frac{3}{p}}$ ($3 < p < \infty$) was derived under an additional condition imposing the vanishing of $u$ at the blowup time $T_*$ in the limit of a rescaling procedure. We also mention the works \cite{CGZ2, CZ4, CZZ1} for related works concerning one-component regularity criteria and \cite{CDP26,CH26,CP26,Pal_26} for recent ill-posedness results for system \eqref{prob:NS}. Nevertheless, all the cited works provide only qualitative blowup criteria, and the proofs therein, which rely on some form of compactness argument and qualitative ingredients, cannot be employed to extract quantitative information near the blowup time.

Quantitative bounds for critical norms of solutions were first established in a breakthrough paper by Tao \cite{Tao_20}, in which compactness arguments and qualitative results were replaced by quantitative ingredients such as Carleman inequalities and quantitative regularity estimates. Specifically, Tao obtained bounds in terms of a triple exponential function of the critical $L^3$-norm. This result implies that the $L^3$-norm of $u$ must blow up at least at a triple logarithmic rate as the time approaches a finite blowup time, namely
\begin{equation} \label{est:quantitative-L3}
\limsup_{t \nearrow T_*} \frac{\|u(t)\|_{L^3}}{\bigl( \ln \ln |\ln (T_* - t)| \bigr)^c} = +\infty,
\end{equation}
for some absolute constant $c > 0$.

Following Tao's work, various significant improvements and refinements of \eqref{est:quantitative-L3} and related results have been developed. These include localized quantitative blowup rates in $L^3$ in \cite{BarPra_21, Bar_23}, an improvement to a double logarithmic blowup rate for axisymmetric solutions in \cite{Pal_21}, a high-dimensional analog of \eqref{est:quantitative-L3} in \cite{Pal_22}, and extensions of \eqref{est:quantitative-L3} to Lorentz spaces in \cite{FenHeWan, OzaPal_23}.

Despite the recent advancements, the study of quantitative blowup rates involving critical Besov spaces still remains widely under-explored. In this paper, we aim to establish analogous quantitative regularity and blowup criteria for solutions to \eqref{prob:NS} in the setting of critical Besov spaces.

\subsection{Main results and strategy of the proofs}
In this paper, we consider the Navier–Stokes equations in the framework of critical Besov spaces $\dot{B}_{p,\infty}^{-1+\frac3p}$ (see the definitions of these spaces in Subsection \ref{subsec:space}).
Our first main result states as follows:

\begin{theorem}[Quantitative regularity] \label{thm:main1}
Let $3 < p < \infty$ and $(u,P)$ be a classical solution of \eqref{prob:NS} in $[0,T] \times \mathbb{R}^3$ such that
\begin{align} \label{assumption}
&\|u\|_{L_T^\infty (\dot{B}_{p,\infty}^{-1+\frac{3}{p}})} \leq M,\\
\label{assumption-Du}
&\sup_{t \in (0,T)} \bigl\| |D|^{-1+\frac{3}{p}} u(t) \bigr\|_{L^p}^p \leq A,
\end{align}
for some $M \geq 2$ and $A \geq 2$. Then for any multi-index $\alpha \in \mathbb{N}^3$ and any $t \in (0,T]$, there holds
\begin{equation} \label{est:main-quadruple}
t^{\frac{1+|\alpha|}{2}} \| \nabla^\alpha u(t) \|_{L^\infty} \lesssim_{|\alpha|}
\exp\left( 2^{|\alpha|} A \exp\bigl( \exp(M^{c_p}) \bigr) \right),
\end{equation}
where $c_p$ is a positive constant depending only on $p$.
\end{theorem}

As a consequence of Theorem \ref{thm:main1}, we derive the following result, which can be interpreted as a quantitative blowup criterion coupling a double exponential of $\|u(t)\|_{\dot{B}_{p,\infty}^{-1+\frac{3}{p}}}$ with $\bigl\| |D|^{-1+\frac{3}{p}} u(t) \bigr\|_{L^p}$, which implies that at least one of these norms must exhibit a specific logarithmic growth rate as time approaches the blowup time.

\begin{theorem} \label{thm:main2}
Let $(u,P)$ be a classical solution of \eqref{prob:NS} on $[0,T_*)$ such that $u$ blows up at a finite time $T_* \in (0,\infty)$. Then for $p \in (3,\infty)$, we have
\begin{equation} \label{est:quantitative-Besovinfty}
\limsup_{t \nearrow T_*} \left[ \frac{ \exp\bigl( \exp( \|u(t)\|_{\dot{B}_{p,\infty}^{-1+\frac{3}{p}}}^{c_p} ) \bigr) }{ \bigl(  |\ln (T_* - t)| \bigr)^b } \bigl\| |D|^{-1+\frac{3}{p}} u(t) \bigr\|_{L^p} \right] = +\infty,
\end{equation}
where $c_p$ is a positive constant depending only on $p$ and $b$ is a positive absolute constant.
\end{theorem}

In order to prove our main results, we substantially improve the strategy in \cite{Tao_20} (see \cite{BarPra_21, Pal_22, OzaPal_23} for comprehensive explanations of Tao's strategy). Let us recall that the idea in \cite{Tao_20} is to show that if a solution $u$ concentrates at a particular time then there will be moderate $L^3$ concentrations at widely separated spatial scales. Precisely, if $\|u\|_{L_t^\infty (L^3)} \leq A$ for some constant $A \geq 2$ and $u$ concentrates at $(t_0,x_0)$ and at frequency $N_0$ for which
\begin{equation} \label{est:PN}
N_0^{-1} |P_{N_0} u(t_0,x_0)| \geq A^{-C_0}
\end{equation}
where $C_0$ is a large absolute constant and $P_{N_0}$ is a Littlewood–Paley projection to frequencies equivalent to $N_0$, then $t_0 N_0^2 \lesssim \exp \exp \exp (A^c)$. In other words, concentrations cannot occur for high frequencies. This key estimate was obtained in \cite{Tao_20} by using Carleman inequalities in combination with several auxiliary results such as bounded total speed, epochs of regularity, back propagation, iterated back propagation, and annuli of regularity. Subsequently, the triple exponential bounds in \cite[Theorem 1.2]{Tao_20} were derived through an adaptation of the local well-posedness and regularity theory for Navier–Stokes as shown in \cite[Section 6]{Tao_20}.

A difficulty that we encounter in our analysis while using the aforementioned strategy stems from the low regularity of the Besov spaces $\dot{B}_{p,\infty}^{-1+\frac{3}{p}}$ with $3 < p < \infty$. Therefore, unlike \cite{Tao_20} where it is sufficient to use the simple decomposition $u = e^{t\Delta} u_0 + w$ to obtain regularity results of $u$, in the present paper, to tackle this difficulty, we perform a more subtle decomposition via iteration (as in \cite{GalKocPla_16}) as
\begin{equation} \label{decomposition-v}
u = u_{1L} + u_{2L} + \cdots + u_{mL} + v,
\end{equation}
with $m$ being a large enough integer (precisely, $m = [p] + 3$; see Lemma \ref{expdl}), where $u_{1L} := e^{t\Delta} u_0$ and $u_{kL}$ ($1 \leq k \leq m$) is determined iteratively by the terms $u_{iL}$, $1 \leq i \leq k-1$ (see Subsection \ref{subsec:decompose} for more detail). Furthermore, in order to exploit assumption \eqref{assumption} and stick with critical quantities, at many points in our analysis, we involve appropriate Kato spaces in which the norm of solutions is invariant under the scaling transformation \eqref{scaling}. We note that similar bounds have been obtained in \cite{Pal_21, Pal_22}.

On one hand, the advantage of \eqref{decomposition-v} is that it allows us to improve the regularity in critical Kato spaces step by step in the sense that the spatial integrability exponent of $u_{kL}$ in the $k$-th step equals $\frac{p}{k}$. Eventually, at the $m$-th step, we deduce that $u_{mL}$ belongs to a critical Kato space with spatial exponent $1$. In this procedure, the higher the spatial regularity of $u_{kL}$ is, the more singular with respect to time it is at the initial time (see Lemma \ref{expdl} for more detail). On the other hand, due to the decomposition \eqref{decomposition-v}, in the equation satisfied by $v$, in addition to the nonlinear term as in \cite{Tao_20}, there also appear a linear term and a force term both involving all $u_{kL}$, $1 \leq k \leq m$ (see the system \eqref{prob:NS-toy}). This complicates significantly the analysis and requires new ingredients to go through the argument. In particular, we establish several delicate estimates on $u_{kL}$ (see Lemma \ref{expdl}) and $v$ (see Lemma \ref{lem:energy} and Lemma \ref{refksv}) which enable us to handle the intricacies induced by the lack of regularity in the proof of the back propagation in Proposition \ref{backprop} and annuli of regularity in Proposition \ref{anreg}.

\medskip
\noindent \textbf{Main novelties and proof strategy.}
Most previous quantitative implementations of the concentration–propagation strategy are formulated under critical Lebesgue or Lorentz-space control. Although qualitative regularity and blowup criteria in critical Besov spaces are available, their compactness-based arguments do not provide the explicit bounds sought here. The endpoint space $\dot{B}_{p,\infty}^{-1+\frac{3}{p}}$ presents a genuinely different difficulty: its norm controls every dyadic block at the critical size, but it provides no summability over the frequency index. Consequently, physical-space estimates that follow directly from an $L^3$ bound cannot be used without further structure.

The first ingredient is the finite iterative decomposition \eqref{decomposition-v}. At each stage one extracts a component $u_{kL}$ satisfying sharp heat-decay estimates in critical Kato spaces, while its spatial integrability improves from $p$ to approximately $p/k$. After finitely many steps the last component reaches the integrability level needed for energy and localized kernel estimates. The remainder $v$ has enough $L^2$-based regularity to recover bounded total speed, epochs of regularity, back propagation, and annuli of regularity. Thus the decomposition converts the non-summable Besov information into a finite collection of quantitatively controlled pieces plus an energy-class remainder.

The second new ingredient concerns the final concentration argument. A positive-kernel approach naturally would require an assumption on $|D|^{-1+\frac{3}{p}} |u|$, which discards the sign and direction of the velocity field. Our hypothesis \eqref{assumption-Du} is instead imposed directly on the signed vector field $|D|^{-1+\frac{3}{p}} u$; in particular, no corresponding bound with $|u|$ is assumed. Near each concentration scale we test this signed quantity against a mean-zero vector field obtained from a localized curl. The mean-zero property produces an additional decay factor. By choosing geometrically separated scales, the associated test functions are almost orthogonal in $L^{p'}$, and duality with $|D|^{-1+\frac{3}{p}} u \in L^p$ bounds the number of possible concentration scales. This yields the key frequency estimate \eqref{inftybound} and hence the quantitative regularity bound \eqref{est:main-quadruple}.

Let us explain the use of the hypothesis \eqref{assumption-Du} in more detail. Put $\sigma = 1 - \frac{3}{p} \in (0,1)$. If one would assume that
\begin{equation}\label{est:strong-positive-assumption}
|D|^{-\sigma} |u| \in L^p,
\end{equation}
then the Riesz-potential representation
\[
|D|^{-\sigma} |u|(x) = c_\sigma \int_{\mathbb{R}^3} \frac{|u(y)|}{|x-y|^{3-\sigma}} \, dy
\]
has a positive kernel and a nonnegative integrand. Consequently, every concentration of $u$ at a given scale produces a positive contribution to this potential, and contributions coming from suitably separated scales cannot cancel one another. One can therefore sum the resulting lower bounds and continue with the scale-counting argument in the spirit of \cite{Tao_20}. Thus, under \eqref{est:strong-positive-assumption}, the last step is a rather direct adaptation of Tao's argument.

This positivity is unavailable under our weaker and more natural signed hypothesis
\begin{equation}\label{est:signed-assumption-intro}
|D|^{-\sigma} u \in L^p.
\end{equation}
Indeed, neither a large value of $u$ nor a concentration of its vorticity gives a pointwise lower bound for $|D|^{-\sigma} u$, because different spatial regions and different vector components may cancel. We overcome this obstruction by converting each vorticity concentration into a nontrivial signed pairing $\langle u, \Psi_i \rangle$, where $\Psi_i$ is a localized curl and hence has zero mean. By self-adjointness,
\[
\langle u, \Psi_i \rangle = \langle |D|^{-\sigma} u, |D|^{\sigma} \Psi_i \rangle.
\]
After the critical normalization $\phi_i = r_i^{-2} |D|^{\sigma} \Psi_i$, the functions $\phi_i$ have uniformly bounded $L^{p'}$ norms. Their mean-zero origin yields improved far-field decay, while the geometric separation of the scales makes the family almost orthogonal in $L^{p'}$. Choosing the sign of each $\phi_i$ according to its pairing and then using Hölder duality gives
\[
K \varepsilon^3 \lesssim
\left\langle |D|^{-\sigma} u,
\sum_{i=1}^K \operatorname{sgn}\bigl( \langle |D|^{-\sigma} u, \phi_i \rangle \bigr) \phi_i \right\rangle
\lesssim \bigl\| |D|^{-\sigma} u \bigr\|_{L^p} K^{1/p'},
\]
which bounds $K$, the number of concentration scales (see \eqref{est:K-epsi} for more detail). This dual almost-orthogonality argument is precisely what replaces the positive summation available for $|D|^{-\sigma} |u|$.

In this way, the proof uses both components of the assumptions in complementary roles: the critical Besov bound drives the parabolic decomposition and all propagation estimates, while the signed $L^p$ bound rules out an arbitrarily long cascade of separated concentration scales. Finally, an argument by contradiction combining Theorem \ref{thm:main1} with the classical Prodi–Serrin–Ladyzhenskaya criterion yields Theorem \ref{thm:main2}.

\vspace{0.1cm}

\textbf{Organization of the paper.} The rest of the paper is organized as follows. In Section \ref{sec:pre}, we recall some basic facts on Littlewood–Paley theory and present an elementary decomposition of solutions via iteration. In Section \ref{sec:basic}, we present some preliminary estimates of the linear and nonlinear components of solutions. Section \ref{Sect4} is devoted to establishing regularity and propagation results. In Section \ref{sec:mainthm}, we employ the result in Section \ref{Sect4} to prove key Lemma \ref{mt}, Theorem \ref{thm:main1} and Theorem \ref{thm:main2}. In Appendixes \ref{Appa} and \ref{Appb}, we present respectively the proof of estimates of linear components \eqref{S2eq12} and \eqref{pwsm1}.

	\section{Preliminaries} \label{sec:pre}
\subsection{Notations}
We list below some notations that will be used frequently throughout the paper.

$\bullet$ The notation $A \gtrsim B$ (resp. $A \lesssim B$) means $A \geq c B$ (resp. $A \leq c B$), where the implicit $c$ is a positive constant depending on some initial parameters. We write $A \lesssim_\kappa B$ if the constant $c = c(\kappa)$ depends on $\kappa$. If both $A \gtrsim B$ and $A \lesssim B$, we write $A \asymp B$.

$\bullet$ For a Schwartz function $f$, the Fourier transform of $f$, denoted by $\hat f$ or $\mathcal{F}(f)$, is defined as
\[
\mathcal{F}(f)(\xi) = (2\pi)^{-\frac{3}{2}} \int_{\mathbb{R}^3} e^{-i x \cdot \xi} f(x) \, dx.
\]
The inverse Fourier transform of $f$ is denoted by $\mathcal{F}^{-1}(f)$.

$\bullet$ We designate by $c_\beta$ a generic positive constant depending only on $\beta$, which may vary on each occurrence. In particular, $c_0$ denotes a universal positive constant. For $c_p \geq 1$ and  $M > 0$ defined in \eqref{assumption}, we denote by $M_1, M_2, \ldots$ large constants such that for any $i \in \mathbb{N}$,
\begin{equation} \label{Mi}
M^{c_p} \leq M_1 \leq M^{d_p c_p}, \qquad M_i^{c_p} \leq M_{i+1} \leq M_i^{d_p c_p},
\end{equation}
where $d_p \gg 1$ is a large constant depending only on $p$.

$\bullet$ For $x \in \mathbb{R}^3$ and $R > 0$, we denote $B(x,R) := \{ y \in \mathbb{R}^3 : |x-y| < R \}$, and for $0 < R_1 < R_2$, $\mathcal{C}(x, R_1, R_2) := \{ y \in \mathbb{R}^3 : R_1 < |x-y| < R_2 \}$. When $x = 0$, we simply write $\mathcal{C}(R_1, R_2)$ for $\mathcal{C}(0, R_1, R_2)$.

$\bullet$ For $T > 0$, we denote $Q_T := (0,T) \times \mathbb{R}^3$.

\subsection{Function spaces} \label{subsec:space}
For $1 \leq p, q \leq \infty$, a time interval $I \subset \mathbb{R}$, and a domain $\Omega \subset \mathbb{R}^3$, the space-time Lebesgue space $L^q(I; L^p(\Omega))$ consists of all functions $f$ defined on $I \times \Omega$ such that
\[
\|f\|_{L^q(I; L^p(\Omega))} := \left( \int_I \|f(t)\|_{L^p(\Omega)}^q \, dt \right)^{\frac{1}{q}} < +\infty,
\]
with
\[
\|f(t)\|_{L^p(\Omega)} := \left( \int_{\Omega} |f(t,x)|^p \, dx \right)^{\frac{1}{p}},
\]
with the usual modifications when $p = \infty$ or $q = \infty$. When $I = [0,T]$ for $T > 0$ and $\Omega = \mathbb{R}^3$, we shall abbreviate $L^q([0,T]; L^p(\mathbb{R}^3))$ by $L_T^q(L^p)$ and $\|f\|_{L^p(\mathbb{R}^3)}$ by $\|f\|_{L^p}$.

Next, for the reader's convenience, we recall some basic facts on Littlewood–Paley theory from \cite{BahCheDan}. Let $\phi$ be a smooth radial function in $\mathbb{R}^3$ satisfying
\begin{equation*}
\begin{aligned}
0 \leq \phi \leq 1 \quad \text{in } \mathbb{R}^3, \qquad \operatorname{supp} \phi \subset \mathcal{C}\left( \frac{3}{4}, \frac{8}{3} \right), \qquad \sum_{j \in \mathbb{Z}} \phi(2^{-j} \xi) = 1, \quad \forall \xi \in \mathbb{R}^3 \setminus \{0\}.
\end{aligned}
\end{equation*}
For any $j \in \mathbb{Z}$, we define the homogeneous dyadic block $\dot{\Delta}_j$ and the homogeneous partial sum operator $\dot{S}_j$ via
\begin{equation*}
\dot{\Delta}_j u := \mathcal{F}^{-1}\bigl( \phi(2^{-j} \xi) \hat{u}(\xi) \bigr), \qquad \dot{S}_j u := \sum_{\ell \leq j-1} \dot{\Delta}_\ell u.
\end{equation*}
We remark that both $\dot{\Delta}_j$ and $\dot{S}_j$ are bounded operators on $L^r(\mathbb{R}^3)$ ($1 \leq r \leq \infty$) with norm independent of $j$.

The dyadic block $\dot{\Delta}_j$ has the quasi-orthogonality property, namely $\dot{\Delta}_j \dot{\Delta}_k = 0$ if $|j-k| > 1$. Hence $\dot{\Delta}_j = \sum_{|j-k| \leq 1} \dot{\Delta}_j \dot{\Delta}_k$.

Let
\[
\mathcal{S}_h' := \left\{ u \in \mathcal{S}'(\mathbb{R}^3) : \lim_{j \to -\infty} \| \dot{S}_j u \|_{L^\infty} = 0 \right\}
\]
be the homogeneous distribution space. It is obvious that $u = \sum_{j \in \mathbb{Z}} \dot{\Delta}_j u$ for any $u \in \mathcal{S}_h'$.

\begin{definition}[see [see {\cite[Subsection 2.3]{BahCheDan}}] ] \label{def2.1}
Let $(p,q) \in [1, +\infty]^2$, $s \in \mathbb{R}$, and $u \in \mathcal{S}_h'$. We set
\[
\|u\|_{\dot{B}_{p,q}^s} := \bigl\| 2^{js} \| \dot{\Delta}_j u \|_{L^p} \bigr\|_{\ell^q(\mathbb{Z})}.
\]

$\bullet$ For $s < \frac{3}{p}$ (or $s = \frac{3}{p}$ if $q = 1$), we define
\[
\dot{B}_{p,q}^s(\mathbb{R}^3) := \left\{ u \in \mathcal{S}_h' : \|u\|_{\dot{B}_{p,q}^s} < \infty \right\}.
\]

$\bullet$ If $k \in \mathbb{N}$ and $\frac{3}{p} + k \leq s < \frac{3}{p} + k + 1$ (or $s = \frac{3}{p} + k + 1$ if $q = 1$), then $\dot{B}_{p,q}^s(\mathbb{R}^3)$ is defined as the set of distributions $u \in \mathcal{S}_h'$ such that $\nabla^\beta u \in \dot{B}_{p,q}^{s-k-1}(\mathbb{R}^3)$ for every multi-index $\beta \in \mathbb{N}^3$ with $|\beta| = k + 1$.
\end{definition}

For $T > 0$, $s \in \mathbb{R}$, and $1 \leq p, q, r \leq \infty$, we denote
\[
L_T^r(\dot{B}_{p,q}^s) = L^r((0,T); \dot{B}_{p,q}^s(\mathbb{R}^3))
:= \left\{ u : Q_T \to \mathbb{R}^3 : \bigl\| \|u(t)\|_{\dot{B}_{p,q}^s} \bigr\|_{L^r(0,T)} < \infty \right\},
\]
with the norm
\[
\|u\|_{L_T^r(\dot{B}_{p,q}^s)} := \bigl\| \|u(t)\|_{\dot{B}_{p,q}^s} \bigr\|_{L^r(0,T)}.
\]

Let us also recall a useful characterization of homogeneous Besov spaces with negative index in terms of Kato-type norms of the heat flow (see e.g. \cite[Theorem 2.34]{BahCheDan}). For $s \in \mathbb{R}$, $1 \leq p, q \leq \infty$, and $0 < T \leq \infty$, we define the Kato space $\mathcal{K}_{p,q}^s(Q_T)$ by
\[
\mathcal{K}_{p,q}^s(Q_T) := \left\{ u \in L_{\mathrm{loc}}^1(Q_T) : \|u\|_{\mathcal{K}_{p,q}^s(Q_T)} < \infty \right\},
\]
where
\begin{equation} \label{Kato-norm}
\begin{aligned}
&\|u\|_{\mathcal{K}_{p,q}^s(Q_T)} := \bigl\| t^{-\frac{s}{2}} \|u(t,\cdot)\|_{L^p} \bigr\|_{L^q((0,T); dt/t)} \quad \text{if } q < \infty, \\
&\|u\|_{\mathcal{K}_{p,\infty}^s(Q_T)} := \sup_{t \in (0,T)} \left( t^{-\frac{s}{2}} \|u(t,\cdot)\|_{L^p} \right).
\end{aligned}
\end{equation}
For any $s < 0$ and $1 \leq p, q \leq \infty$, there exists a positive constant $C$ depending only on $s$ such that
\begin{equation} \label{Besov-characterization}
C^{-1} \|u\|_{\dot{B}_{p,q}^s} \leq \left\| e^{t\Delta} u \right\|_{\mathcal{K}_{p,q}^s(Q_\infty)} \leq C \|u\|_{\dot{B}_{p,q}^s}, \quad \forall u \in \mathcal{S}_h'.
\end{equation}

Next, we recall Bony's para-product decomposition. Formally, the product of $f, g \in \mathcal{S}_h'(\mathbb{R}^3)$ may be decomposed into
\begin{equation} \label{bony}
fg = T_f g + T_g f + R(f,g),
\end{equation}
where
\begin{equation*}
\begin{aligned}
&T_f g := \sum_{j \in \mathbb{Z}} \dot{S}_{j-1} f \, \dot{\Delta}_j g, \\
&R(f,g) := \sum_{j \in \mathbb{Z}} \dot{\Delta}_j f \, \tilde{\Delta}_j g, \qquad \text{with} \quad \tilde{\Delta}_j g := \sum_{j' = j-1}^{j+1} \dot{\Delta}_{j'} g.
\end{aligned}
\end{equation*}
The operator $T$ is called the homogeneous para-product of $g$ by $f$, which represents low–high frequency interactions, while the operator $R$ is called the homogeneous remainder of $f$ and $g$, which corresponds to the high–high frequency interactions.

As a consequence of Bony's decomposition, we have the following useful formula:
\begin{equation} \label{form:decomposition}
\dot{\Delta}_j (fg) = \dot{\Delta}_j \Bigl( \sum_{|j'-j| \leq 4} \dot{S}_{j'-1} f \, \dot{\Delta}_{j'} g + \sum_{|j'-j| \leq 4} \dot{S}_{j'-1} g \, \dot{\Delta}_{j'} f + \sum_{j' \geq j-4} \dot{\Delta}_{j'} f \, \tilde{\Delta}_{j'} g \Bigr).
\end{equation}

We recall the well-known Bernstein lemma, which will be used frequently throughout the paper.
\begin{lemma}[{\cite[Lemma 2.1]{BahCheDan}}] \label{lem:Bernstein}
Let $\mathcal{B}$ be a ball and $\mathcal{C}$ be an annulus of $\mathbb{R}^3$ with center $0$. There exists a positive constant $C$ such that for any $k \in \mathbb{N}$, any $1 \leq p \leq q \leq \infty$, and any function $u \in L^p$, the following estimates hold.

i) If $\operatorname{supp} \hat{u} \subset \lambda \mathcal{B}$ for some $\lambda > 0$, then
\begin{equation} \label{est:Bernstein-1}
\|D^k u\|_{L^q} := \sum_{|\alpha| = k} \| \nabla^\alpha u \|_{L^q} \leq C^{k+1} \lambda^{k + 3\left( \frac{1}{p} - \frac{1}{q} \right)} \|u\|_{L^p}.
\end{equation}

ii) If $\operatorname{supp} \hat{u} \subset \lambda \mathcal{C}$ for some $\lambda > 0$, then
\begin{equation} \label{est:Bernstein-2}
C^{-k-1} \lambda^{k} \|u\|_{L^p} \leq \|D^k u\|_{L^p} \leq C^{k+1} \lambda^{k} \|u\|_{L^p}.
\end{equation}
\end{lemma}

It is easy to observe from Definition \ref{def2.1} and Lemma \ref{lem:Bernstein} that for any $s \in \mathbb{R}$, $1 \leq p_1 \leq p_2 \leq \infty$, and $1 \leq q_1 \leq q_2 \leq \infty$, one has
\[
\dot{B}_{p_1,q_1}^s(\mathbb{R}^3) \hookrightarrow \dot{B}_{p_2,q_2}^{s - 3\left( \frac{1}{p_1} - \frac{1}{p_2} \right)}(\mathbb{R}^3).
\]

Next we recall the following local multiplier theorem from \cite{Tao_20}.

\begin{lemma}[{\cite[Lemma 2.1]{Tao_20}}] \label{localY}
Let $\Omega$ be an open subset of $\mathbb{R}^3$. For $N > 0$ and $\nu \geq 1$, we denote
\[
\Omega_{\nu/N} := \bigl\{ x \in \mathbb{R}^3 : \operatorname{dist}(x, \Omega) \leq \nu/N \bigr\}.
\]
Let $\mathcal{T}_{\mathfrak{m}}(f) := \mathcal{F}^{-1}(\mathfrak{m}(\xi) \hat{f}(\xi))$ be the Fourier multiplier with symbol $\mathfrak{m}$, which satisfies $\operatorname{supp} \mathfrak{m} \subset B(0,N)$ and
\[
|\nabla^j \mathfrak{m}(\xi)| \lesssim \mu N^{-j}, \qquad 0 \leq j \leq 100.
\]
Then for any $1 \leq p_i \leq q_i \leq \infty$, $i = 1, 2$, with $q_2 \geq q_1$, we have
\[
\|\mathcal{T}_{\mathfrak{m}} f\|_{L^{q_1}(\Omega)}
\lesssim \mu N^{\frac{3}{p_1} - \frac{3}{q_1}} \|f\|_{L^{p_1}(\Omega_{\nu/N})}
+ \nu^{-1000} \mu |\Omega|^{\frac{1}{q_1} - \frac{1}{q_2}} N^{\frac{3}{p_2} - \frac{3}{q_2}} \|f\|_{L^{p_2}(\mathbb{R}^3)}.
\]
\end{lemma}
Note that $\nu^{-1000}$ above can be replaced by any negative power of $\nu$.

\subsection{Decomposition via iteration} \label{subsec:decompose}
In this subsection, we will decompose the solution to \eqref{prob:NS} into linear and nonlinear components.

Let $\mathbb{P}$ be the Leray projection operator defined by $\mathbb{P} := \mathrm{Id} + \nabla (-\Delta)^{-1} \nabla \cdot$. This is a Fourier multiplier with matrix-valued symbol which is homogeneous of degree zero and smooth away from the origin. By applying $\mathbb{P}$ to the momentum equations of \eqref{prob:NS}, we obtain
\begin{equation} \label{eq:projected-NS}
\partial_t u - \Delta u = - \mathbb{P} \nabla \cdot (u \otimes u).
\end{equation}
Assume that $u$ satisfies the initial condition $u|_{t=0} = u_0$. By Duhamel's formula, we write
\begin{equation} \label{Duhamel}
u(t) = e^{t\Delta} u_0 - \int_0^t e^{(t-\tau)\Delta} \mathbb{P} \nabla \cdot (u \otimes u)(\tau) \, d\tau, \qquad t \in (0,T).
\end{equation}
Let $m \in \mathbb{N}$ be a large integer. We decompose $u$ as
\begin{equation}\label{decompv}
u = u_{1L} + u_{2L} + \cdots + u_{mL} + v,
\end{equation}
where $u_{1L}, \ldots, u_{mL}$ solve respectively the following equations:
\begin{equation} \label{prob:u1L}
\left\{
\begin{aligned}
&\partial_t u_{1L} - \Delta u_{1L} + \nabla P_{1L} = 0 \quad \text{in } Q_T, \\
&\nabla \cdot u_{1L} = 0 \quad \text{in } Q_T, \qquad u_{1L}|_{t=0} = u_0,
\end{aligned}
\right.
\end{equation}
and
\begin{equation*}
\left\{
\begin{aligned}
&\partial_t u_{2L} - \Delta u_{2L} + \nabla P_{2L} = - u_{1L} \cdot \nabla u_{1L} \quad \text{in } Q_T, \\
&\nabla \cdot u_{2L} = 0 \quad \text{in } Q_T, \qquad u_{2L}|_{t=0} = 0,
\end{aligned}
\right.
\end{equation*}
and
\begin{equation*}
\left\{
\begin{aligned}
&\partial_t u_{3L} - \Delta u_{3L} + \nabla P_{3L} = - u_{2L} \cdot \nabla u_{2L} - u_{1L} \cdot \nabla u_{2L} - u_{2L} \cdot \nabla u_{1L} \quad \text{in } Q_T, \\
&\nabla \cdot u_{3L} = 0 \quad \text{in } Q_T, \qquad u_{3L}|_{t=0} = 0,
\end{aligned}
\right.
\end{equation*}
and for any $3 \leq k \leq m-1$,
\begin{equation}\label{prob:lin}
\left\{
\begin{aligned}
&\partial_t u_{(k+1)L} - \Delta u_{(k+1)L} + \nabla P_{(k+1)L} = - u_{kL} \cdot \nabla u_{kL} - \sum_{i=1}^{k-1} \bigl( u_{iL} \cdot \nabla u_{kL} + u_{kL} \cdot \nabla u_{iL} \bigr), \\
&\nabla \cdot u_{(k+1)L} = 0 \quad \text{in } Q_T, \qquad u_{(k+1)L}|_{t=0} = 0.
\end{aligned}
\right.
\end{equation}
Then, in view of \eqref{prob:NS}, \eqref{decompv}, and \eqref{prob:u1L}–\eqref{prob:lin}, we deduce that $v$ solves
\begin{equation} \label{prob:NS-toy}
\left\{
\begin{aligned}
&\partial_t v - \Delta v + v \cdot \nabla v + V_1 \cdot \nabla v + v \cdot \nabla V_1 + \nabla \bar{P} \\
&\hspace{2cm} = - V_1 \cdot \nabla V_3 - V_3 \cdot \nabla V_2 \quad \text{in } Q_T, \\
&\nabla \cdot v = 0 \quad \text{in } Q_T, \qquad v|_{t=0} = 0,
\end{aligned}
\right.
\end{equation}
where
\begin{equation}\label{def:Vi}
V_1 := \sum_{i=1}^{m} u_{iL}, \qquad V_2 := \sum_{i=1}^{m-1} u_{iL}, \qquad V_3 := u_{mL}.
\end{equation}

Let $N_0 \in \mathbb{N}$ be a large integer. As in \cite{CHN1}, we use Duhamel's formula for $u_{iL}$, the definition of Kato norms \eqref{Kato-norm}, and the characterization \eqref{Besov-characterization} to obtain the following estimates:
\[
\sup_{\substack{0 \leq |\alpha| \leq N_0 \\ q \in [p,\infty]}}
\| \nabla^\alpha u_{1L} \|_{\mathcal{K}_{q,\infty}^{-\left(|\alpha| + 1 - \frac{3}{q}\right)}(Q_T)}
\lesssim \|u_0\|_{\dot{B}_{p,\infty}^{-1+\frac{3}{p}}},
\]
and for any $1 \leq k \leq m$,
\begin{equation}\label{S2eq12}
\sup_{\substack{0 \leq |\alpha| \leq N_0 \\ q \in [\max\{ \frac{p}{k}, 1 \}, \infty]}}
\| \nabla^\alpha u_{kL} \|_{\mathcal{K}_{q,\infty}^{-\left(|\alpha| + 1 - \frac{3}{q}\right)}(Q_T)}
\lesssim \left( 1 + \|u_0\|_{\dot{B}_{p,\infty}^{-1+\frac{3}{p}}} \right)^{2^{k-1}}.
\end{equation}
We shall postpone the proof of \eqref{S2eq12} to Appendix \ref{Appa}.

Therefore, if $\|u\|_{L_T^\infty(\dot{B}_{p,\infty}^{-1+\frac{3}{p}})} \leq M$ for some $M \geq 1$, and for $m$ large enough, we deduce from \eqref{S2eq12} that
\begin{equation} \label{condVi}
\sum_{i=1}^2 \sup_{\substack{0 \leq |\alpha| \leq N_0 \\ q \in [p,\infty]}}
\| \nabla^\alpha V_i \|_{\mathcal{K}_{q,\infty}^{-\left(|\alpha| + 1 - \frac{3}{q}\right)}(Q_T)}
+ \sup_{\substack{0 \leq |\alpha| \leq N_0 \\ q \in [1,\infty]}}
\| \nabla^\alpha V_3 \|_{\mathcal{K}_{q,\infty}^{-\left(|\alpha| + 1 - \frac{3}{q}\right)}(Q_T)}
\lesssim M^{2^{m-1}}.
\end{equation}		
		
		\section{Some preliminary  estimates} \label{sec:basic}
        In this section, we establish various estimates of $u_{kL}$ and $v$ in Besov spaces with time weights. We start with an interpolation inequality.

		\begin{lemma} \label{lem:interpolation} Assume $p>3$. Then for any  $2< r \leq 3$, we have
			\begin{equation} \label{est:interpolation-Lq} \|w\|_{L^r}\lesssim\|w\|_{L^2}^{\alpha_1}\|\nabla w\|_{L^2}^{\alpha_2}\|w\|_{\dot B_{p,\infty}^{-1+\frac{3}{p}}}^{\alpha_3}, \quad \forall w \in H^1(\R^3) \cap \dot{B}_{p,\infty}^{-1+\frac{3}{p}}(\R^3),
			\end{equation}
			where $\alpha_1,\alpha_2,\alpha_3\in ( 0,1)$ are given by
			\begin{equation*}
				\alpha_1=\frac{2(4p+r-pr-6)}{pr+2p-4r}, \quad  \alpha_2=\frac{6(p-r)(r-2)}{r(pr+2p-4r)}, \quad  \alpha_3=\frac{3p(r-2)^2}{r(pr+2p-4r)},
			\end{equation*}
			and $\alpha_1+\alpha_2+\alpha_3=1.$
		\end{lemma}
		\begin{remark}
			We note that the condition $2<r\leq 3<p$ ensures $\alpha_i \in (0,1)$.	
		\end{remark}

		\begin{proof}
			Let $j_0 \in \Z$ to be determined later on. We have
			\begin{equation} \label{eq:decompose-j0}
				\|w\|_{L^r}\lesssim\sum_{j\leq j_0}\|\dot{\Delta}_j w\|_{L^r}+\sum_{j\geq j_0}\|\dot{\Delta}_j w\|_{L^r}.
			\end{equation}
			Since $r>2$, we deduce from  \eqref{est:Bernstein-1} that
			\begin{equation} \label{est:u-j<j0}
				\sum_{j\leq j_0}\|\dot{\Delta}_j w\|_{L^r}\lesssim \sum_{j\leq j_0}2^{\frac{3(r-2)}{2r}j}\|\dot{\Delta}_j w\|_{L^2}\lesssim 2^{\frac{3(r-2)}{2r}j_0}\|w\|_{L^2}.
			\end{equation}
To handle  the second term on the right-hand side of \eqref{eq:decompose-j0}, we get, by using H\"older's inequality, that
			\begin{equation} \label{est:u-j>j_0}
				\sum_{j\geq j_0}\|\dot{\Delta}_j w\|_{L^r}
				\lesssim \sum_{j\geq j_0}\|\dot{\Delta}_j  w\|_{L^2}^{\frac{2(p-r)}{(p-2)r}}\|\dot{\Delta}_j w\|_{L^p}^{\frac{p(r-2)}{(p-2)r}}.
			\end{equation}
From  \eqref{est:Bernstein-2}, we deduce
			$$ \|\dot{\Delta}_j  w\|_{L^2} \lesssim 2^{-j}\|\dot{\Delta}_j  \nabla w\|_{L^2} \lesssim 2^{-j}\|\nabla w \|_{L^2}.
			$$ 	
Keeping in mind that $pr-4p-r+6<0$ due to the assumption $r \leq 3<p$, we thus obtain
			\begin{equation} \label{est:u-j>j0-1}	
				\begin{aligned}
					\sum_{j\geq j_0}\|\dot{\Delta}_j w\|_{L^r} &\lesssim  \sum_{j\geq j_0}\bigl(2^{-j} \| \nabla w\|_{L^{2}}\bigr)^{\frac{2(p-r)}{(p-2)r}} \bigl(2^{(1-\frac{3}{p})j}\|w\|_{\dot B_{p,\infty}^{-1+\frac{3}{p}}}\bigr)^{\frac{p(r-2)}{(p-2)r}} \\
					&\lesssim \sum_{j\geq j_0} 2^{\frac{pr-4p-r+6}{(p-2)r}j}\| \nabla w\|_{L^{2}}^{\frac{2(p-r)}{(p-2)r}} \|w\|_{\dot B_{p,\infty}^{-1+\frac{3}{p}}}^{\frac{p(r-2)}{(p-2)r}} \\
					&\lesssim 2^{\frac{pr-4p-r+6}{(p-2)r}j_0} \| \nabla w\|_{L^{2}}^{\frac{2(p-r)}{(p-2)r}} \|w\|_{\dot B_{p,\infty}^{-1+\frac{3}{p}}}^{\frac{p(r-2)}{(p-2)r}}.
				\end{aligned}
			\end{equation}
		By	inserting \eqref{est:u-j<j0} and \eqref{est:u-j>j0-1} into \eqref{eq:decompose-j0}, we achieve
			$$
			\|w\|_{L^r} \lesssim 2^{\frac{3(r-2)}{2r}j_0}\|w\|_{L^2} + 2^{\frac{pr-4p-r+6}{(p-2)r}j_0} \| \nabla w\|_{L^{2}}^{\frac{2(p-r)}{(p-2)r}} \|w\|_{\dot B_{p,\infty}^{-1+\frac{3}{p}}}^{\frac{p(r-2)}{(p-2)r}}, \quad \forall j_0 \in \Z.
			$$
			By minimizing the right-hand side of the above estimate over $j_0 \in \Z$, we obtain \eqref{est:interpolation-Lq}.
		\end{proof}
		
	Next we are going to use the standard energy method to establish energy estimates for the function $v$ in \eqref{decompv}.
		
		\begin{lemma}\label{lem:energy}
			Assume that $v$ solves \eqref{prob:NS-toy} and satisfies
			\begin{equation} \label{est:v-Besov}
				\| v \|_{L_T^\infty (\dot{B}_{p,\infty}^{-1+\frac{3}{p}})} \leq M^{c_p}.	
			\end{equation}	
			Then one has
			\begin{equation} \label{est:LinftyL2v}
				\| v \|_{\cK_{2,\infty}^{\frac{1}{2}}(Q_T)} +\sup_{t\in (0,T)}(t^{-\frac{1}{4}}\|\nabla v\|_{L^2((0,t) \times\mathbb{R}^3)}) \leq M^{c_p}.			\end{equation}
		\end{lemma}
		\begin{proof}
			By taking $L^2$ inner product of \eqref{prob:NS-toy} with $v$ and using the integration by parts, we obtain, for any $\tau \in (0,T)$, that
			\begin{align} \nonumber
				&\frac{1}{2} \frac{d}{d\tau} \| v(\tau) \|_{L^2}^2  + \|\nabla v(\tau) \|_{L^2}^2 \\ \nonumber
				&\lesssim  \int_{\R^3} |v(\tau,x)|^2 |\nabla V_1(\tau,x)| dx  + \int_{\R^3}|v(\tau,x)\|V_1(\tau,x)\|\nabla V_3(\tau,x)|dx \\ \nonumber
				&\quad+ \int_{\R^3} |v(\tau,x)\|V_3(\tau,x)\|\nabla V_2(\tau,x)| dx \\ \label{est:energy-1}
				&\lesssim   \| v(\tau) \|_{L^{\frac{2p}{p-1}}}^2 \| \nabla V_1(\tau) \|_{L^p} + \| v(\tau) \|_{L^2} \| V_1(\tau)\|_{L^p} \| \nabla V_3(\tau) \|_{L^{\frac{2p}{p-2}}} \\ \nonumber
				&\quad + \| v(\tau)\|_{L^2} \| \| V_3(\tau)\|_{L^{\frac{2p}{p-2}}} \| \nabla V_2\|_{L^p}.
			\end{align}
Next we shall handle term by term above.
			
		\begin{equation*}\label{energy2}
		\begin{aligned}
		\| v(\tau) \|_{L^{\frac{2p}{p-1}}}^2 \| \nabla V_1(\tau) \|_{L^p} &\leq CM^{\frac{2p^2-5p+6}{p(2p-5)}} \tau^{-\frac{2p-3}{2p}} \| v(\tau)\|_{L^2}^{\frac{4(p-3)(p-1)}{p(2p-5)}} \| \nabla v(\tau)\|_{L^2}^{\frac{6(p-3)}{p(2p-5)}}\\
		&\leq CM^{\frac{2p^2-5p+6}{2p^2-8p+9}}\tau^{-\frac{(2p-3)(2p-5)}{2(2p^2-8p+9)}}\|v(\tau)\|_{L^2}^{\frac{4(p-3)(p-1)}{2p^2-8p+9}} + \frac{1}{2}\|\nabla v(\tau)\|_{L^2}^2\\
		&\leq CM^{\frac{2p^2-5p+6}{2p^2-8p+9}}\tau^{-\frac{1}{2}}\|v\|_{\cK_{2,\infty}^{\frac{1}{2}}(Q_\tau)}^{\frac{4(p-3)(p-1)}{p(2p-5)}} + \frac{1}{2}\|\nabla v(\tau)\|_{L^2}^2.
		\end{aligned}
		\end{equation*}
		Note that $\frac{4(p-3)(p-1)}{p(2p-5)}=:\alpha_p<2.$ We deduce from  \eqref{condVi} that
				\begin{equation*}
				\| V_1(\tau)\|_{L^p} \| \nabla V_3(\tau) \|_{L^{\frac{2p}{p-2}}} +  \| V_3(\tau) \|_{L^{\frac{2p}{p-2}}} \| \nabla V_2(\tau)\|_{L^p}\lesssim M^{2^{m-1}}\tau^{-\frac{3}{4}},
				\end{equation*}
which leads to
			\begin{align*}\label{energy3}
				\| v(\tau) \|_{L^2} \bigl(\| V_1(\tau)\|_{L^p} \| \nabla V_3(\tau) \|_{L^{\frac{2p}{p-2}}} +  \| V_3(\tau) \|_{L^{\frac{2p}{p-2}}} \| \nabla V_2(\tau)\|_{L^p}\bigr) \leq CM^{c_p} \tau^{-\frac{1}{2}}\| v \|_{\cK_{2,\infty}^{\frac{1}{2}}(Q_\tau)}.
			\end{align*}
By substituting the above estimates into  \eqref{est:energy-1},  we obtain
	\begin{equation*}
	\begin{aligned}
	\frac{1}{2}\frac{d}{d\tau}\|v(\tau)\|^2_{L^2}+\|\nabla v(\tau)\|_{L^2}^2\leq CM^{c_p}\tau^{-\frac{1}{2}}\bigl(\|v\|_{\mathcal{K}_{2,\infty}^{\frac{1}{2}}(Q_\tau)}^{\alpha_p}
+\|v\|_{\mathcal{K}_{2,\infty}^{\frac{1}{2}}(Q_\tau)}\bigr)+\frac{1}{2}\|\nabla v(\tau)\|_{L^2}^2
	\end{aligned}
	\end{equation*}
with $\alpha_p<2$. Note that $\|v(0)\|_{L^2}=0$, hence we get, by  integrating the above inequality over $(0,t)$, that
	\begin{equation*}
	\|v(t)\|_{L^2}^2+\int_0^t\|\nabla v(\tau)\|_{L^2}^2d\tau\lesssim M^{c_p}t^{\frac{1}{2}}\bigl(\|v\|_{\mathcal{K}_{2,\infty}^{\frac{1}{2}}(Q_t)}^{\alpha_p}+\|v\|_{\mathcal{K}_{2,\infty}^{\frac{1}{2}}(Q_t)}\bigr),
	\end{equation*}
	from which, we infer
	\begin{equation*}
	\|v\|_{\mathcal{K}_{2,\infty}^{\frac{1}{2}}(Q_T)}^2+\sup_{t\in(0,T)}\Bigl(t^{-\frac{1}{2}}\|\nabla v\|_{L^2(Q_t)}^2\Bigr)\lesssim M^{c_p}\bigl(\|v\|_{\mathcal{K}_{2,\infty}^{\frac{1}{2}}(Q_T)}^{\alpha_p}+\|v\|_{\mathcal{K}_{2,\infty}^{\frac{1}{2}}(Q_T)}\bigr).
	\end{equation*}
Since $\alpha_p<2$, by using Young's inequality, we conclude the proof of  \eqref{est:LinftyL2v}.
		\end{proof}
		
Next we deal with the functions $u_{kL}$ in \eqref{decompv}.
	
		\begin{lemma}\label{expdl}
			Assume that $\| u_0 \|_{\dot B_{p,\infty}^{-1+\frac{3}{p}}} \leq M$ with $p\in (3,\infty)$. Let $m=[p]+3$ and $u_{kL}$ be as in \eqref{prob:u1L}--\eqref{prob:lin} for $1\leq k\leq m-1$. Then for any $1 \leq k \leq m$, $ r_* \geq \max\{3,\frac{p}{k}\}$ and $j \in \Z$, we have
			\begin{equation}\label{est:ulindecay}
				\|\dot\Delta_ju_{kL}(t)\|_{L^{r_*}}\lesssim 2^{\left(1-\frac{3}{r_*}\right)j}e^{-c_{k-1}t2^{2j}}M^{\widetilde{c}_{k-1}},
			\end{equation}
			where $c_k$, $\widetilde{c}_k$ are positive constants depending only on $k$. Furthermore, there holds
			\begin{equation}\label{est:forcev}
				\sup_{t \in (0,T)}\Bigl(t^{-1}\|u_{mL}(t)\|_{\dot B_{1,\infty}^0}+\|u_{mL}(t)\|_{\dot B_{1,\infty}^2}\Bigr)\lesssim M^{c_p}.
			\end{equation}
		\end{lemma}
		\begin{proof} For any $t\in[0,T]$, we divide the proof into several steps.

\noindent{\sc Step 1.} The proof of \eqref{est:ulindecay} for the case with $1\leq k\leq [\frac{p}{3}]+1.$

We shall prove \eqref{est:ulindecay} by induction. We first observe from  Lemma \ref{lem:Bernstein}
that it suffices to prove \eqref{est:ulindecay} for
 the case $r_*=\max\{3,\frac{p}{k}\}$.
 In view of \eqref{prob:u1L}, we have $u_{1L}(t)=e^{t\Delta}u_0.$ Then it follows from \cite[Lemma 2.4]{BahCheDan}
 that \begin{equation}\label{S3eq1}
			\| \dot{\Delta}_j u_{1L}(t) \|_{L^p} \lesssim e^{-c_0 t 2^{2j}}\| \dot{\Delta}_j u_0\|_{L^p} \lesssim 2^{\left(1-\frac{3}{p}\right)j} e^{-c_0 t 2^{2j}}M, \quad \forall t\geq 0,
\end{equation}
	from which and  Lemma \ref{lem:Bernstein}, we deduce that for any $r\geq  p$,
			\begin{equation*}
				\|\dot{\Delta}_ju_{1L}\|_{L_t^\infty (L^r)}+2^{2j}\|\dot{\Delta}_ju_{1L}\|_{L_t^1(L^r)}\lesssim 2^{(1-\frac{3}{r})j}M, \quad  \forall j \in \Z.
			\end{equation*}
			
			Next we shall prove \eqref{est:ulindecay} for $k=2$. By Duhamel's formula and Bony's decomposition \eqref{form:decomposition}, for any $j \in \Z$, we have
			\begin{equation} \label{est:induc-n=2} \begin{aligned}
					\| \dot{\Delta}_j u_{2L}(t) \|_{L^{\frac{p}{2}}}\lesssim & \sum_{|j'-j|\leq 4} \int_0^t \| e^{(t-s)\Delta} \bP \nabla \cdot \dot{\Delta}_j(\dot{S}_{j'-1} u_{1L}(s) \otimes  \dot{\Delta}_{j'} u_{1L}(s)) \|_{L^{\frac{p}{2}}}ds \\
					&+\sum_{|j'-j|\leq 4} \int_0^t \| e^{(t-s)\Delta}\bP \nabla \cdot  \dot{\Delta}_j( \dot{\Delta}_{j'} u_{1L}(s)\otimes \dot{S}_{j'-1} u_{1L}(s)) \|_{L^{\frac{p}{2}}}ds \\
					&+ 	\sum_{\ell\geq j-4} \int_0^t \| e^{(t-s)\Delta}\bP \nabla \cdot \dot{\Delta}_j(\dot{\Delta}_\ell u_{1L}(s) \otimes  \tilde{\Delta}_{\ell}u_{1L}(s)) \|_{L^{\frac{p}{2}}}ds.
				\end{aligned}
			\end{equation}
It is easy to observe from \eqref{S3eq1} that
$\|\dot{S}_{\ell} u_{1L}\|_{L_t^\infty (L^p)}\lesssim 2^{\left(1-\frac{3}{p}\right){\ell}} M$ for any $\ell \in \Z.$ Then we deduce from
\cite[Lemma 2.4]{BahCheDan} that
			\begin{align*}
				&\sum_{|j'-j|\leq 4} \int_0^t \| e^{(t-s)\Delta}\bP \nabla \cdot \dot{\Delta}_j(\dot{S}_{j'-1} u_{1L}(s) \otimes  \dot{\Delta}_{j'}u_{1L}(s)) \|_{L^{\frac{p}{2}}}ds \\
				&\lesssim \left( \int_0^t  2^j e^{-c_0(t-s)2^{2j}} 2^{(1-\frac{3}{p})j} 2^{(1-\frac{3}{p})j}e^{-c_0s2^{2j}}ds \right) M^2\lesssim e^{-2^{-8}c_0t2^{2j}}2^{(1-\frac{6}{p})j}M^2.
			\end{align*}
The second term on the right-hand side of \eqref{est:induc-n=2} shares the same estimate.\\
By \eqref{S3eq1}, we have
			\begin{equation*}
				\sum_{\ell\geq j-4}\|\dot\Delta_\ell u_{1L}(s)\|_{L^p}\|\tilde\Delta_{\ell}u_{1L}(s)\|_{L^p}\lesssim e^{-2^{-8}c_0s2^{2j}}\sum_{\ell\geq j-4}e^{-c_0s2^{2\ell}}2^{\left(2-\frac{6}{p}\right)\ell} M^2,
			\end{equation*}
from which, we infer
			\begin{align*}
				&\sum_{\ell\geq j-4} \int_0^t \| e^{(t-s)\Delta}\bP\cdot\dot{\Delta}_j(\dot{\Delta}_\ell u_{1L}(s) \otimes  \tilde{\Delta}_{\ell} u_{1L}(s)) \|_{L^{\frac{p}{2}}}ds \\
				&\lesssim M^2\int_0^t2^je^{-c_0(t-s)2^{2j}}e^{-2^{-8}c_0s2^{2j}}\Bigl(\sum_{\ell\geq j-4}e^{-c_0s2^{2\ell}}2^{\left(2-\frac{6}{p}\right)\ell}\Bigr)ds  \\
				&\lesssim M^2 2^j e^{-2^{-8}c_0t2^{2j}}\Bigl(\int_0^t \sum_{\ell\geq j-4}e^{-c_0s2^{2\ell}}2^{\left(2-\frac{6}{p}\right)\ell}\,ds \Bigr) \lesssim  e^{-c_2t2^{2j}} 2^{\left(1-\frac{6}{p}\right)j}M^2.
			\end{align*}
By substituting the above estimates into \eqref{est:induc-n=2}, we achieve  \eqref{est:ulindecay} for $k=2$.\\
Now  we suppose that estimate  \eqref{est:ulindecay} holds for any $2 \leq k\leq n-1$. We shall show that \eqref{est:ulindecay}  also holds for $k=n$. In view of \eqref{prob:lin}, we write
			\begin{equation} \label{est:Deltaj-ukL}
				\begin{aligned}
					\dot{\Delta}_j u_{nL}(t) &= - \int_0^t  e^{(t-s)\Delta}\bP \nabla \cdot \dot\Delta_j(u_{(n-1)L}\otimes u_{(n-1)L})(s)ds	\\
					&-\sum_{i=1}^{n-2} \int_0^t  e^{(t-s)\Delta} \bP \nabla \cdot \dot\Delta_j \bigl(u_{iL} \otimes u_{(n-1)L}+u_{(n-1)L} \otimes u_{iL}\bigr)(s)\,ds, \quad t \in [0,T].
				\end{aligned}
			\end{equation}
We shall only deal  with the typical term involving  $u_{1L} \otimes  u_{(n-1)L}.$  The other terms can be treated along the same line.
Once again by using Bony's decomposition \eqref{form:decomposition}, we write
\begin{equation}\label{S3eq3}
			\begin{split}
				&\int_0^t  e^{(t-s)\Delta}\bP \nabla \cdot \dot{\Delta}_j(u_{(n-1)L}\otimes u_{1L})(s)\,ds=\text{I}_j^n(t)+\text{II}_j^n(t)+\text{III}_j^n(t) \with\\
&\text{I}_j^n(t):=\text{I}_j(u_{(n-1)L},u_{1L})(t),\quad \text{II}_j^n(t):=\text{II}_j(u_{(n-1)L},u_{1L})(t),\quad \text{III}_j^n(t):=\text{III}_j(u_{(n-1)L},u_{1L})(t),\\
			&\text{I}_j(f,g)(t):=\sum_{|j'-j|\leq 4}\int_0^t e^{(t-s)\Delta}\bP \nabla \cdot \dot{\Delta}_j(\dot{S}_{j'-1} f\otimes \dot{\Delta}_{j'} g)(s)\,ds, \\
				&\text{II}_j(f,g)(t):=\sum_{|j'-j|\leq 4}\int_0^t  e^{(t-s)\Delta}\bP \nabla \cdot \dot{\Delta}_j (\dot{\Delta}_{j'}f\otimes \dot{S}_{j'-1} g)(s)ds, \\
				&\text{III}_j(f,g)(t):=\sum_{\ell\geq j-4}\int_0^t e^{(t-s)\Delta}\bP \nabla \cdot \dot{\Delta}_j (\dot{\Delta}_{\ell} f\otimes \tilde {\Delta}_{\ell} g)(s)\,ds.
			\end{split}\end{equation}
			Observing from the inductive assumption that
			$\|\dot{S}_{\ell} u_{(n-1)L}(s)\|_{L^{\frac{p}{n-1}}} \lesssim 2^{\left(1-\frac{3(n-1)}{p}\right)\ell}M^{\wt{c}_{n-2}},$
	from which and \eqref{S3eq1}, we infer
			\begin{align*}
				\|\text{I}_j^n(t)\|_{L^{\frac{p}{n}}} &\lesssim
				\sum_{|j'-j|\leq 4} \int_0^t2^je^{-c_0(t-s)2^{2j}}e^{-c_0 s2^{2j}}2^{\left(2-\frac{3n}{p}\right)j'}M^{\wt{c}_{n-2}+1}
\lesssim e^{-2^{-8}c_0 t2^{2j}}2^{\left(1-\frac{3n}{p}\right)j}M^{\wt{c}_{n-2}+1}.
			\end{align*}
			Similarly, we deduce from the inductive assumption that
			\begin{align*}
				\|\text{II}_j^n(t)\|_{L^{\frac{p}{n}}} &\lesssim
				\sum_{|j'-j|\leq 4} \int_0^t2^je^{-c_0(t-s)2^{2j}}e^{-c_{n-2}s2^{2j'}}2^{\left(2-\frac{3n}{p}\right)j'}M^{\wt{c}_{n-2}+1}\,ds \\
				&\lesssim e^{-2^{-8}c_{n-2}t2^{2j}}2^{\left(1-\frac{3n}{p}\right)j}M^{\wt{c}_{n-2}+1}.
			\end{align*}
	Finally observing that for any $0<s<t$,
			\begin{equation*}
				\sum_{\ell\geq j-4}\|\dot\Delta_\ell u_{1L}(s)\|_{L^p}\|\tilde \Delta_{\ell}u_{(n-1)L}(s)\|_{L^{\frac{p}{n-1}}}\lesssim e^{-2^{-8}c_0s2^{2j}}\sum_{\ell\geq j-4}e^{-c_{n-2}s2^{2\ell}}2^{\left(2-\frac{3n}{p}\right)\ell}M^{\wt{c}_{n-2}+1},
			\end{equation*}
	we infer
			\begin{align*}
				\|\text{III}_j^n(t)\|_{L^{\frac{p}{n}}} &\lesssim  \int_0^t2^j e^{-c_0(t-s)2^{2j}}e^{-2^{-8}c_0s2^{2j}}\sum_{\ell\geq j-4}e^{-c_{n-2}s2^{2\ell}}2^{\left(2-\frac{3n}{p}\right)\ell}ds M^{\wt{c}_{n-2}+1} \\
				&\lesssim 2^j e^{-2^{-8}c_0 t 2^{2j}} \int_0^t \sum_{\ell\geq j-4}e^{-c_{n-2}s2^{2\ell}}2^{\left(2-\frac{3n}{p}\right)\ell}ds M^{\wt{c}_{n-2}+1} \\
				&\lesssim e^{-2^{-8}c_{n-2} t 2^{2j}} 2^{\left(1-\frac{3n}{p}\right)j}M^{\wt{c}_{n-2}+1}.
			\end{align*}
		By summarizing the above estimates, we arrive at
			\begin{align*}
				&\left\|\int_0^t  e^{(t-s)\Delta}\bP \nabla \cdot \dot{\Delta}_j(u_{(n-1)L}\otimes u_{1L})(s)ds\right \|_{L^{\frac{p}{n}}} \lesssim e^{-c_{n-1}t2^{2j}}2^{\left(1-\frac{3n}{p}\right)j}M^{\wt{c}_{n-2}+1}.
			\end{align*}
Similar estimate holds for the other terms on the right-hand side of \eqref{est:Deltaj-ukL}. This leads to \eqref{est:ulindecay} for $k=n.$
  By the induction, \eqref{est:ulindecay} holds for for any $1 \leq k \leq [\frac{p}{3}]+1$.

  \noindent{\sc Step 2.} The proof of \eqref{est:ulindecay} for the case for $ [\frac{p}{3}]+1<k\leq m.$\\
  We shall first prove by induction that for any $l \in \N$ and
  $\max\{1,\frac{3p}{p+3l}\}<r<3$, there holds
			\begin{equation}\label{est:inductionv0}
				\|\dot\Delta_ju_{(k_*+l)L}(t)\|_{L^r}\lesssim M^{\wt{c}_{k_\ast+l-1}}\min\bigl\{t^{\frac{1}{2}},2^{-j}\bigr\}^{\left(\frac{3}{r}-1\right)}e^{-c_{k_*+l-1}t2^{2j}},
			\end{equation}
			where $k_*:=[\frac{p}{3}]+1$ and $c_{k_*+l-1}$ is a positive constant depending only on $p$ and $l$.
			
To handle  the case $l=1$ and  $\frac{3p}{p+3}<r<3$ in \eqref{est:inductionv0}, we only consider the typical term involving $u_{k_*L}\otimes u_{1L}$ in \eqref{est:Deltaj-ukL} for $n=k_*+1.$ Notice that $r>\frac{3p}{p+3},$  for $\text{I}_j^{k_*+1}$ given by \eqref{S3eq3} for $n=k_*+1,$ we get, by applying the result in {\bf Step 1}, that
\begin{align*}
				\|\text{I}_j^{k_*+1}(t)\|_{L^r}&\lesssim\sum_{|j'-j|\leq 4}\int_0^t2^je^{-c_0(t-s)2^{2j}}\|\dot S_{j'-1}u_{k_*L}(s)\|_{L^{\frac{pr}{p-r}}}\|\dot\Delta_{j'}u_{1L}(s)\|_{L^p}ds\\
				&\lesssim \sum_{|j'-j|\leq 4}\int_0^t2^je^{-c_0(t-s)2^{2j}}2^{\left(1-\frac{3}{r}+\frac{3}{p}\right)j}2^{\left(1-\frac{3}{p}\right)j}e^{-c_0s2^{2j}}\,ds M^{\wt{c}_{k_\ast-1}+1}\\
				&\lesssim \sum_{|j'-j|\leq 4}e^{-2^{-8}c_0t2^{2j}}\int_0^t2^{\left(3-\frac{3}{r}\right)j}e^{-(1-2^{-8})c_0(t-s)2^{2j}}\,ds M^{\wt{c}_{k_\ast-1}+1}\\
				&\lesssim \min\bigl\{t^{\frac{1}{2}},2^{-j}\big\}^{\frac{3}{r}-1}e^{-2^{-8}c_0t2^{2j}}M^{\wt{c}_{k_\ast-1}+1},
			\end{align*}
		where in the last step, we have used the facts that
			\begin{equation}\label{intexpest}
				\begin{aligned}
					&\int_0^t2^{\alpha j}e^{-cs2^{2j}}ds\lesssim 2^{(\alpha-2)j},\quad \forall \alpha\in\mathbb{R},\\
					&\int_0^t2^{\alpha j}e^{-cs2^{2j}}ds\lesssim \int_0^ts^{-\frac{\alpha}{2}}ds\lesssim t^{1-\frac{\alpha}{2}},\quad 0\leq \alpha<2.
				\end{aligned}
			\end{equation}
Similarly, in view of \eqref{S3eq3} for $n=k_*+1$, we infer
			\begin{align*}
				\|\text{II}_j^{k_*+1}(t)\|_{L^r}&\lesssim\sum_{|j'-j|\leq 4}\int_0^t2^je^{-c_0(t-s)2^{2j}}\|\dot S_{j'-1}u_{1L}(s)\|_{L^p}\|\dot\Delta_{j'}u_{k_*L}(s)\|_{L^{\frac{pr}{p-r}}}ds\\
				&\lesssim \sum_{|j'-j|\leq 4}\int_0^t2^je^{-c_0(t-s)2^{2j}}2^{\left(1-\frac{3}{p}\right)j}2^{\left(1-\frac{3}{r}+\frac{3}{p}\right)j}e^{-c_{k_*-1}s2^{2j}}\,ds M^{\wt{c}_{k_\ast-1}+1}\\
				&\lesssim \sum_{|j'-j|\leq 4}2^{-8}e^{-c_{k_*-1}t2^{2j}}\int_0^t2^{\left(3-\frac{3}{r}\right)j}e^{-(1-2^{-8})c_{k_*-1}(t-s)2^{2j}}\,ds M^{\wt{c}_{k_\ast-1}+1}\\
				&\lesssim \min\bigl\{t^{\frac{1}{2}},2^{-j}\big\}^{\frac{3}{r}-1}e^{-2^{-8}c_{k_*-1}t2^{2j}}M^{\wt{c}_{k_\ast-1}+1}.
			\end{align*}
	Finally we deduce from the result in {\bf Step 1} that
			\begin{align*}
				\sum_{\ell\geq j-4}\|\dot{\Delta}_\ell u_{1L}(s)\otimes \tilde \Delta_{\ell}u_{k_*L}(s)\|_{L^{r}}&\lesssim 2^{\left(1-\frac{3}{r}+\frac{3}{p}\right)j}e^{-2^{-8}c_0s2^{2j}}\sum_{\ell\geq j-4}e^{-c_{k_*-1}s2^{2\ell}}2^{\left(1-\frac{3}{p}\right)\ell}M^{\wt{c}_{k_\ast-1}+1}\\
				&\lesssim 2^{\left(1-\frac{3}{r}+\frac{3}{p}\right)j}e^{-2^{-8}c_0s2^{2j}}s^{-\frac{1}{2}}\sum_{\ell\geq j-4}2^{-\frac{3}{p}\ell}M^{\wt{c}_{k_\ast-1}+1}\\
				&\lesssim s^{-\frac{1}{2}}2^{\left(1-\frac{3}{r}\right)j}e^{-2^{-8}c_0s2^{2j}}M^{\wt{c}_{k_\ast-1}+1},
			\end{align*}
from which and \eqref{S3eq3} for $n=k_*+1,$ we infer
			\begin{align*}
				\|\text{III}_j^{k_*+1}(t)\|_{L^r}&\lesssim e^{-2^{-8}c_0t2^{2j}}\int_0^t2^{\left(2-\frac{3}{r}\right)j}e^{-\left(1-2^{-8}\right)c_0(t-s)2^{2j}}s^{-\frac{1}{2}}\,ds M^{\wt{c}_{k_\ast-1}+1}\\
				&\lesssim \min\bigl\{t^{\frac{1}{2}},2^{-j}\big\}^{\frac{3}{r}-1}e^{-2^{-8}c_0t2^{2j}}M^{\wt{c}_{k_\ast-1}+1}.
			\end{align*}
By summarizing the above estimates and thanks to \eqref{S3eq3}, we obtain
\begin{equation*}
\Bigl\|\int_0^te^{(t-s)\Delta}\bP \nabla \cdot \dot{\Delta}_j(u_{k_*L}\otimes u_{1L})(s)ds\Bigr\|_{L^r}
\lesssim \min\bigl\{t^{\frac{1}{2}},2^{-j}\big\}^{\frac{3}{r}-1}e^{-2^{-8}c_0t2^{2j}}M^{\wt{c}_{k_\ast-1}+1},
\end{equation*}  for any $\frac{3p}{p+3}<r<3.$  The other terms in  \eqref{est:Deltaj-ukL} for $n=k_*+1$ share the same estimate as above.
As a consequence, we obtain, for any $\frac{3p}{p+3}<r<3$,
			\begin{equation}\label{S3eq4} \|\dot\Delta_ju_{(k_*+1)L}(t)\|_{L^r}\lesssim\min\bigl\{t^{\frac{1}{2}},2^{-j}\big\}^{\frac{3}{r}-1}e^{-c_{k_*}t2^{2j}}M^{\wt{c}_{k_\ast}},
			\end{equation}
			which is \eqref{est:inductionv0} for $l=1$.
			
			Next we suppose that \eqref{est:inductionv0} holds for any $1 \leq l\leq n-1$. Let $\frac{3p}{p+3n}<r<3$, we shall show that \eqref{est:inductionv0} holds for $l=n$. Notice that if
			$r>\frac{3p}{p+3n},$  $\frac{pr}{p-r}>\frac{3p}{p+3(n-1)}$. Then in view of \eqref{S3eq3} corresponding to $k_*+n,$
 we get, by using the inductive assumption, that
			\begin{align}
				\|\text{I}_j^{k_*+n}(t)\|_{L^r}&\lesssim\sum_{|j'-j|\leq 4}\int_0^t2^je^{-c_0(t-s)2^{2j}}\|\dot S_{j'-1}u_{(k_*+n-1)L}(s)\|_{L^{\frac{pr}{p-r}}}\|\dot\Delta_{j'}u_{1L}(s)\|_{L^p}ds\label{eq3.18}\\
				&\lesssim \sum_{|j'-j|\leq 4}\int_0^t2^je^{-c_0(t-s)2^{2j}}\Bigl(\sum_{j'\leq j+2}2^{(\frac{3}{\gamma}-\frac{3}{r}+\frac{3}{p})j'}s^{\frac{1}{2}(\frac{3}{\gamma}-1)}\Bigr) 2^{\left(1-\frac{3}{p}\right)j'}e^{-c_0s2^{2j}}\,ds M^{\wt{c}_{k_\ast+n-2}+1}\nonumber\\
				&\lesssim \sum_{|j'-j|\leq 4}e^{-2^{-9}c_0t2^{2j}}\int_0^t2^{\left(3-\frac{3}{r}\right)j}e^{-\left(1-2^{-9}\right)c_0(t-s)2^{2j}}\,ds M^{\wt{c}_{k_\ast+n-2}+1}\nonumber\\
				&\lesssim \min\bigl\{t^{\frac{1}{2}},2^{-j}\big\}^{\frac{3}{r}-1}e^{-2^{-9}c_0t2^{2j}}M^{\wt{c}_{k_\ast+n-2}+1},\nonumber
			\end{align}
			for some $\gamma\in \bigl(\frac{3p}{p+3(n-1)},\frac{pr}{p-r}\bigr)$.\\
Similarly,  we deduce from the inductive assumption that
			\begin{align}
				\|\text{II}_j^{k_*+n}(t)\|_{L^r}&\lesssim\sum_{|j'-j|\leq 4}\int_0^t2^je^{-c_0(t-s)2^{2j}}\|\dot S_{j'-1}u_{1L}(s)\|_{L^p}\|\dot\Delta_{j'}u_{(k_*+n-1)L}(s)\|_{L^{\frac{pr}{p-r}}}ds\label{eq3.19}\\
				&\lesssim \sum_{|j'-j|\leq 4}\int_0^t2^je^{-c_0(t-s)2^{2j}}2^{\left(1-\frac{3}{p}\right)j}2^{\left(1-\frac{3}{r}+\frac{3}{p}\right)j}e^{-c_{k_*+n-2}s2^{2j}}\,ds M^{\wt{c}_{k_\ast+n-2}+1}\nonumber\\
				&\lesssim \sum_{|j'-j|\leq 4}e^{-2^{-8}c_{k_*+n-2}t2^{2j}}\int_0^t2^{\left(3-\frac{3}{r}\right)j}e^{-(1-2^{-8})c_0(t-s)2^{2j}}\,ds M^{\wt{c}_{k_\ast+n-2}+1}\nonumber\\
				&\lesssim \min\bigl\{t^{\frac{1}{2}},2^{-j}\big\}^{\frac{3}{r}-1}e^{-2^{-8}c_{k_*+n-2}t2^{2j}} M^{\wt{c}_{k_\ast+n-2}+1},\nonumber
			\end{align}
where we have used once again  \eqref{intexpest} in the last step.\\
Finally due to $r<3$, we obtain for any $0<s<t$,
			\begin{align*}
				&\sum_{\ell\geq j-4}\|\Delta_\ell u_{1L}(s)\otimes\tilde \Delta_{\ell}u_{(k_*+n-1)L}(s)\|_{L^{r}} \\
&\leq \sum_{\ell\geq j-4} \|\Delta_\ell u_{1L}(s)\|_{L^p}\| \tilde \Delta_{\ell}u_{k_*+n-1)L}(s)\|_{L^{\frac{pr}{p-r}}}\\
				&\lesssim e^{-2^{-9}c_0s2^{2j}}\sum_{\ell\geq j-4}e^{-2^{-1}c_0s2^{2\ell}}2^{\left(1-\frac{3}{p}\right)\ell}2^{\left(1-\frac{3}{r}+\frac{3}{p}\right)\ell}M^{\wt{c}_{k_\ast+n-2}+1}\\
				&\lesssim e^{-2^{-9}c_0s2^{2j}}s^{-\frac{1}{2}}\sum_{\ell\geq j-4}2^{\left(1-\frac{3}{r}\right)\ell}M^{\wt{c}_{k_\ast+n-2}+1}\\
				&\lesssim s^{-\frac{1}{2}}2^{\left(1-\frac{3}{r}\right)j}e^{-2^{-9}c_0s2^{2j}}M^{\wt{c}_{k_\ast+n-2}+1}.
			\end{align*}
		Then	it follows from \eqref{S3eq3} corresponding to $k_*+n$  that
			\begin{align}
				\|\text{III}_j^{k_*+n}(t)\|_{L^r}&\lesssim e^{-2^{-9}c_0t2^{2j}}\int_0^t2^{\left(2-\frac{3}{r}\right)j}e^{-\left(1-2^{-9}\right)c_0(t-s)2^{2j}}s^{-\frac{1}{2}}\,ds M^{\wt{c}_{k_\ast+n-2}+1}\label{eq3.20}\\
				&\lesssim\min\bigl\{t^{\frac{1}{2}},2^{-j}\big\}^{\frac{3}{r}-1}e^{-2^{-9}c_0t2^{2j}}M^{\wt{c}_{k_\ast+n-2}+1}.\nonumber
			\end{align}
By substituting the above estimates into \eqref{S3eq3} corresponding to $k_*+n,$  we obtain
\begin{equation*}
\Bigl\|\int_0^te^{(t-s)\Delta}\bP \nabla \cdot \dot{\Delta}_j(u_{(k_*+n-1)L}\otimes u_{1L})(s)ds\Bigr\|_{L^r}
\lesssim \min\bigl\{t^{\frac{1}{2}},2^{-j}\big\}^{\frac{3}{r}-1}e^{-2^{-9}c_0t2^{2j}}M^{\wt{c}_{k_\ast+n-2}+1},
\end{equation*}  for any $\frac{3p}{p+3n}<r<3.$  The other terms in
\eqref{est:Deltaj-ukL} with output index $k_*+n$ share the same estimate.
As a consequence, we obtain, for any $\frac{3p}{p+3n}<r<3$,
			\begin{equation*}
				 \|\dot\Delta_ju_{(k_*+n)L}(t)\|_{L^r}\lesssim\min\bigl\{t^{\frac{1}{2}},2^{-j}\big\}^{\frac{3}{r}-1}e^{-c_{k_*+n-1}t2^{2j}}M^{\wt{c}_{k_\ast+n-2}+1}.
			\end{equation*}
			By induction, we conclude that \eqref{est:inductionv0} holds for any $l \in \N$.\\
We infer from \eqref{est:inductionv0} and Lemma \ref{lem:Bernstein} that \eqref{est:ulindecay} holds for $ [\frac{p}{3}]+1<k\leq m$. Finally, note that \eqref{eq3.18},\eqref{eq3.19} and \eqref{eq3.20} all holds for $r=1$ if $\frac{3p}{p+3n}<1$, and the estimate \eqref{est:forcev} follows. We thus complete the proof of Lemma \ref{expdl}.
\end{proof}

\begin{remark}
It is easy to observe from \eqref{est:ulindecay} that
\begin{equation}\label{S3eq5}
\|u_{kL}\|_{L_T^\infty(\dot B_{r,\infty}^{-1+\frac{3}{r}})}\lesssim M^{c_p}\quad \mbox{for} \ 1\leq k\leq m \andf r\geq\max\bigl\{3,\frac{p}k\bigr\}.
\end{equation}
\end{remark}

Thanks to the bounds on $u_{kL}$ in Lemma \ref{expdl} and in \eqref{S3eq5}, we can derive estimates on $v$ in suitable Besov spaces with temporal weight.
		
		\begin{lemma}\label{refksv}
			Assume that $u$ solves \eqref{prob:NS} on $[0,T]$ and satisfies \eqref{assumption}. Let $m=[p]+3$
 and $v$ solve \eqref{prob:NS-toy}. Then one has
			\begin{equation}\label{est:vbesov}
		\|v\|_{X_T}:=		\sup_{t\in(0,T)}\left(t^{-1}\|v(t)\|_{\dot B_{1,\infty}^0}+\|v(t)\|_{\dot B_{1,\infty}^2}\right)\lesssim M^{c_p}.
			\end{equation}
		\end{lemma}
		\begin{proof}
		In view of \eqref{prob:NS-toy},	 we  write
			\begin{equation}\label{decopv}
				v=v_N + v_L + v_F \with v_N(t):=v_N(0,t),\ v_L(t):=v_L(0,t), \ v_F(t):=v_F(0,t),
			\end{equation}
			where
			\begin{equation} \label{form:vF}
				\begin{aligned}
					&v_N(t_0,t):=-\int_{t_0}^te^{(t-s)\Delta}\bP \nabla \cdot (v\otimes v)(s)ds, \\
					&v_{L}(t_0,t):=-\int_{t_0}^te^{(t-s)\Delta}\bP \nabla \cdot (v\otimes V_1+V_1\otimes v)(s)ds,\\
					&v_F(t_0,t):=-\int_{t_0}^te^{(t-s)\Delta}\bP \nabla \cdot (V_3\otimes V_1 + V_2\otimes V_3)(s)ds,
				\end{aligned}
			\end{equation}
			and $V_i$, $i=1,2,3$, are defined by \eqref{def:Vi}.\\
			Below we handle term by term in \eqref{decopv}.
			
			\noindent \textbf{The estimate of $v_N$.} For any $j \in \Z$, by virtue of \eqref{S3eq3}, we decompose $\dot{\Delta}_j v_N$ as
			\begin{equation} \label{Delta_jv_N}
				\begin{aligned}
					\dot\Delta_jv_N(t)=\text{I}_j(v,v)(t)+ \text{II}_j(v,v)(t)+\text{III}_j(v,v)(t).
				\end{aligned}
			\end{equation}
where $\text{I}_j(v,v),$ $\text{II}_j(v,v)$ and $\text{III}_j(v,v)$ are given by \eqref{S3eq3}. Let
			\begin{equation} \label{r-alpha}
				r \in \Bigl(1,\frac{p}{p-1}\Bigr)  \text{ be so close to $1$ that }  \kappa:=\frac{p(r-1)}{r(p-1)} \in (0,1).
			\end{equation}
We first observe from  \eqref{assumption} and \eqref{S3eq5} that
\begin{equation} \label{S3eq7} \|v\|_{L_T^\infty(\dot B_{p,\infty}^{-1+\frac{3}{p}})}\lesssim M^{c_p},\end{equation}
	so that we have for any $0<s<t$,
			\begin{equation*} \label{est:I-a}\begin{aligned}
 \sum_{|j'-j|\leq 4}\|\dot S_{j'-1} v(s)\|_{L^p}\|\dot\Delta_{j'}v(s)\|_{L^\frac{p}{p-1}} &\lesssim \sum_{|j'-j|\leq 4} M^{c_p}2^{\left(1-\frac{3}{p}\right)j'}2^{\left(\frac{3}{r}+\frac{3}{p}-3\right)j'}\|\dot\Delta_{j'} v(s)\|_{L^r}\\
					&\lesssim M^{c_p}2^{\left(\frac{3}{r}-2\right)j}\sum_{|j'-j|\leq 4} \|\dot\Delta_{j'} v\|_{L_T^\infty(L^{1})}^{1-\kappa}\|\dot\Delta_{j'} v\|_{L_T^\infty(L^{p})}^\kappa\\
& \lesssim M^{c_p}\|v\|_{X_T}^{1-\kappa}2^{-j},
				\end{aligned}
			\end{equation*}
from which, we infer
			\begin{equation*}\label{est:vnl}
				\|\text{I}_j(v,v)(t)\|_{L^1} \lesssim \|v\|_{X_T}^{1-\kappa}M^{c_p}\int_0^te^{-c_0(t-s)2^{2j}}ds\lesssim \|v\|_{X_T}^{1-\kappa}M^{c_p}\min\bigl\{t,2^{-2j}\bigr\}.
			\end{equation*}
We note that $\text{II}_j(v,v)$ shares the same estimate.

Similarly, with $r \in \bigl(1,\frac{p}{p-1}\bigr)$ and $\kappa=\frac{p(r-1)}{r(p-1)} \in (0,1),$ we get, by using H\"older's inequality
and Lemma \ref{lem:Bernstein}, that for any $0<s<t$,
			\begin{align*}
				\sum_{\ell\geq j-4}\|\dot\Delta_{\ell} v(s) \otimes \tilde\Delta_{\ell}v(s)\|_{L^1} &\leq
				\sum_{\ell\geq j-4}\|\dot\Delta_{\ell} v(s)\|_{L^p}\|\tilde\Delta_{\ell}v(s)\|_{L^{\frac{p}{p-1}}} \\
				&\lesssim \sum_{\ell\geq j-4}\|\dot\Delta_{\ell}v(s)\|_{L^p}\|\tilde\Delta_{\ell}v(s)\|_{L^{r}}2^{\left(\frac{3}{r}+\frac{3}{p}-3\right)\ell} \\
				&\lesssim \sum_{\ell\geq j-4}\|\dot\Delta_{\ell}v(s)\|_{L^p}^{1+\kappa}\|\tilde\Delta_{\ell}v(s)\|_{L^1}^{1-\kappa}2^{\left(\frac{3}{r}+\frac{3}{p}-3\right)\ell} \lesssim 2^{-j}\|v\|_{X_T}^{1-\kappa}M^{c_p},
			\end{align*}
			which implies
			\begin{equation*}
				\|\text{III}_j(v,v)(t)\|_{L^1}\lesssim \|v\|_{X_T}^{1-\kappa}M^{c_p}\int_0^te^{-c_0(t-s)2^{2j}}ds\lesssim \|v\|_{X_T}^{1-\kappa}M^{c_p}\min\{t,2^{-2j}\}.
			\end{equation*}
By summarizing the above estimates, we obtain
			\begin{equation*} \label{est:DeltajvN}
				\|\dot\Delta_jv_N(t)\|_{L^1}\lesssim \|v\|_{X_T}^{1-\kappa}M^{c_p}\min\bigl\{t,2^{-2j}\bigr\},
			\end{equation*}
			which, together with Definition \ref{def2.1}, ensures that
			\begin{equation} \label{est:vN_normX}
				\|v_N\|_{X_{T}}\lesssim \|v\|_{X_{T}}^{1-\kappa}M^{c_p}.
			\end{equation}
			
			\noindent \textbf{The estimate of $v_L$.} Once again we get, by using decomposition \eqref{S3eq3}, that
			\begin{equation*}
				\begin{aligned}
					\dot\Delta_jv_L(t)=\sum_{i=1}^m\bigl(\text{I}_{j}(v,u_{iL})(t)+\text{II}_{j}(v,u_{iL})(t)+\text{III}_{j}(v,u_{iL})(t)\bigr).
				\end{aligned}
			\end{equation*}
Let $r$ and $\kappa$ be given by \eqref{r-alpha}. We get, by using Lemma \ref{lem:Bernstein}, \eqref{est:ulindecay} and \eqref{S3eq5}, that
			\begin{align*}
				\|\dot S_{\ell-1} v(s)\otimes\dot\Delta_\ell u_{iL}(s)\|_{L^1}&\lesssim\|\dot S_{\ell-1} v(s)\|_{L^{r}}\|\dot\Delta_\ell u_{iL}(s)\|_{L^{\frac{r}{r-1}}}\\
				&\lesssim\|\dot S_{\ell-1} v(s)\|_{L^{r}}\|\dot\Delta_\ell u_{iL}(s)\|_{L^p}2^{\left(\frac{3}{r}+\frac{3}{p}-3\right)\ell}\\
				&\lesssim \sum_{\ell'\leq \ell-2}\|\dot\Delta_{\ell'} v\|_{L^\infty_T(L^1)}^{1-\kappa}\|\dot\Delta_{\ell'}v\|_{L^\infty_T(L^p)}^{\kappa}\|\dot\Delta_\ell u_{iL}(s)\|_{L^p}2^{\left(\frac{3}{r}+\frac{3}{p}-3\right)\ell}\\
				&\lesssim\sum_{\ell'\leq \ell-2}s^{1-\kappa}\|v\|_{X_T}^{1-\kappa}  2^{\left(1-\frac{3}{p}\right)\kappa \ell'} 2^{\left(\frac{3}{r}-2\right)\ell}e^{-c_{i-1}s2^{2\ell}}M^{c_p}\\
				&\lesssim 2^{-2(1-\kappa)\ell}2^{\left(1-\frac{3}{p}\right)\kappa \ell}2^{\left(\frac{3}{r}-2\right)\ell}M^{c_p}\|v\|_{X_T}^{1-\kappa} \lesssim M^{c_p}\|v\|_{X_T}^{1-\kappa}2^{-\ell},
			\end{align*}
from which and \eqref{S3eq3}, we deduce that
			\begin{equation*}
				\sum_{i=1}^m \|\text{I}_{j}(v,u_{iL})(t)\|_{L^1}\lesssim M^{c_p}\|v\|_{X_T}^{1-\kappa}\int_0^te^{-c_0(t-s)2^{2j}}ds\lesssim M^{c_p}\|v\|_{X_T}^{1-\kappa}\min\bigl\{t,2^{-2j}\bigr\}.
			\end{equation*}
Since $r$ close to $1,$ $\frac{r}{r-1}>3,$ we get, by applying \eqref{est:ulindecay}, that for any $0<s<t$,
			\begin{align*}
				\|\dot S_{\ell-1}u_{iL}(s)\otimes\dot\Delta_\ell v(s)\|_{L^1}&\lesssim\|\dot S_{\ell-1} u_{iL}(s)\|_{L^\frac{r}{r-1}}\|\dot\Delta_\ell v(s)\|_{L^r}\\
				&\lesssim M^{c_p}2^{(\frac{3}{r}-2)\ell}\|\dot\Delta_\ell v\|_{L^\infty_T(L^1)}^{1-\kappa}\|\dot\Delta_\ell v\|_{L^\infty_T(L^p)}^{\kappa}\\
				&\lesssim M^{c_p}\|v\|_{X_T}^{1-\kappa}2^{-\ell},
			\end{align*}
			which, together with \eqref{S3eq3}, implies
			\begin{equation*}
				\sum_{i=1}^m\|\text{II}_{j}(v,u_{iL})(t)\|_{L^1}\lesssim M^{c_p}\|v\|_{X_T}^{1-\kappa}\int_0^te^{-(t-s)2^{2j}}ds\lesssim M^{c_p}\|v\|_{X_T}^{1-\kappa}\min\bigl\{t,2^{-2j}\bigr\}.
			\end{equation*}
	
Finally notice that for any $0<s<t$
			\begin{align*}
				\sum_{\ell\geq j-4}\|\dot\Delta_\ell u_{iL}(s)\|_{L^{\frac{r}{r-1}}}\|\dot\Delta_{\ell}v(s)\|_{L^{r}}&\lesssim \sum_{\ell\geq j-4}2^{\left(\frac{3}{r}+\frac{3}{p}-3\right)\ell}\|\dot\Delta_\ell u_{iL}(s)\|_{L^p}\|\dot\Delta_{\ell}v\|_{L^\infty_T(L^1)}^{1-\kappa}\|\dot\Delta_{l}v\|_{L^\infty_T(L^p)}^{\kappa}\\
				&\lesssim 2^{-j}\|v\|_{X_T}^{1-\kappa}M^{c_p},
			\end{align*}
from which we infer that
			\begin{equation*}
				\sum_{i=1}^m\|\text{III}_{j}(v,u_{iL})(t)\|_{L^1}\lesssim M^{c_p}\|v\|_{X_T}^{1-\kappa}\int_0^te^{-c_0(t-s)2^{2j}}ds\lesssim M^{c_p}\|v\|_{X_T}^{1-\kappa}\min\bigl\{t,2^{-2j}\bigr\}.
			\end{equation*}
By summarizing the above estimates, we achieve
			\begin{equation*}
				\|\dot\Delta_jv_L(t)\|_{L^1}\lesssim M^{c_p}\|v\|_{X_T}^{1-\kappa}\min\bigl\{t,2^{-2j}\bigr\}.
			\end{equation*}
As a result, it comes out
			\begin{equation} \label{est:vL_normX}
				\|v_{L}\|_{X_T}\lesssim \|v\|_{X_T}^{1-\kappa}M^{c_p}.
			\end{equation}

			\noindent \textbf{The estimate of $v_F$.} We are going to prove that
			\begin{equation} \label{est:vF_normX}
				\|v_F\|_{X_T}\lesssim M^{c_p}.
			\end{equation}
 Indeed thanks to \eqref{def:Vi} and \eqref{form:vF}, it suffices to prove that for any $1\leq i\leq m$,
	\begin{equation}\label{S3eq6}
	\Bigl\|\int_0^te^{(t-s)\Delta}\mathbb{P}\nabla\cdot(u_{mL}\otimes u_{iL})(s)ds \Bigr\|_{X_T}\lesssim M^{c_p}.
	\end{equation}
	By using decomposition \eqref{S3eq3}, we have
	\begin{equation*}
	\dot{\Delta}_j\left(\int_0^te^{(t-s)\Delta}\mathbb{P}\nabla\cdot(u_{mL}\otimes u_{iL})(s)ds\right)=\text{I}_j(u_{mL},u_{iL})+\text{II}_j(u_{mL},u_{iL})+\text{III}_j(u_{mL},u_{iL}).
	\end{equation*}
	For $r$ given by \eqref{r-alpha}, we get, by applying Lemma \ref{lem:Bernstein}, \eqref{est:ulindecay} and \eqref{est:forcev}, that
	\begin{equation*}
	\begin{aligned}
	\|\dot S_{\ell-1}u_{mL}(s)\otimes\dot{\Delta}_\ell u_{iL}(s)\|_{L^1}\lesssim& \|\dot S_{\ell-1}u_{mL}(s)\|_{L^r}\|\dot{\Delta}_\ell u_{iL}(s)\|_{L^{\frac{r}{r-1}}}\\
	\lesssim&\sum_{\ell'\leq \ell-2}\|\dot \Delta_{\ell'}u_{mL}(s)\|_{L^r}\|\dot{\Delta}_\ell u_{iL}(s)\|_{L^{\frac{r}{r-1}}}\\
	\lesssim &M^{c_p}\sum_{\ell' \leq \ell-2} s2^{3(1-\frac{1}{r})\ell'}2^{(\frac{3}{r}-2)\ell} e^{-c_{i}s2^{2\ell}}\\
	\lesssim &M^{c_p}s2^{\ell}e^{-c_is2^{2\ell}} \lesssim M^{c_p}2^{-\ell},
	\end{aligned}
	\end{equation*}
	from which and \eqref{S3eq3}, we deduce that
	\begin{equation*}
	\|\text{I}_j(u_{mL},u_{iL})(t)\|_{L^1}\lesssim M^{c_p}\int_0^te^{-c_0(t-s)2^{2j}}ds\lesssim M^{c_p}\min\bigl\{t,2^{-2j}\bigr\}.
	\end{equation*}

It follows from \eqref{est:forcev} and \eqref{S3eq5} that for any $0<s<t$,
	\begin{align*}
	\|\dot S_{\ell-1}u_{iL}(s)\otimes\dot{\Delta}_\ell u_{mL}(s)\|_{L^1}\lesssim& \|\dot S_{\ell-1}u_{iL}(s)\|_{L^p}\|\dot{\Delta}_\ell u_{mL}(s)\|_{L^{\frac{p}{p-1}}}\\
	\lesssim&\Bigl(\sum_{\ell'\leq \ell-2}\|\dot \Delta_{\ell'}u_{iL}(s)\|_{L^p}\Bigr)\|\dot{\Delta}_\ell u_{mL}(s)\|_{L^{\frac{p}{p-1}}}\\
	\lesssim &M^{c_p}\Bigl(\sum_{\ell' \leq \ell-2} 2^{\left(1-\frac{3}{p}\right)\ell'}\Bigr)2^{(\frac{3}{p}-2)\ell}\lesssim M^{c_p}2^{-\ell},
	\end{align*}
from	which and \eqref{S3eq3}, we infer
	\begin{equation*}
	\|\text{II}_j(u_{mL},u_{iL})(t)\|_{L^1}\lesssim M^{c_p}\int_0^te^{-c_0(t-s)2^{2j}}ds\lesssim M^{c_p}\min\{t,2^{-2j}\}.
	\end{equation*}

Finally it is easy to observe from \eqref{est:forcev} and \eqref{S3eq5} that
	\begin{align*}
	\sum_{\ell\geq j-4}\|\dot{\Delta}_{\ell}u_{mL}(s)\|_{L^r}\|\tilde{\Delta}_\ell u_{iL}\|_{L^{\frac{r}{r-1}}}\lesssim M^{c_p}\sum_{\ell\geq j-4}2^{(1-\frac{3}{r})\ell}2^{(\frac{3}{r}-2)\ell}\lesssim M^{c_p}2^{-j},
	\end{align*}
which implies
	\begin{equation*}
	\|\text{III}_j(u_{mL},u_{iL})(t)\|_{L^1}\lesssim M^{c_p}\int_0^te^{-c_0(t-s)2^{2j}}ds\lesssim M^{c_p}\min\bigl\{t,2^{-2j}\bigr\}.
	\end{equation*}
	By summarizing the above estimates, we achieve \eqref{S3eq6}. We complete the proof of \eqref{est:vF_normX}.

By substituting \eqref{est:vN_normX}, \eqref{est:vL_normX} and \eqref{est:vF_normX} into \eqref{decopv}, we arrive at
			\begin{equation*}
				\|v\|_{X_{T}}\lesssim \|v\|_{X_{T}}^{1-\kappa}M^{c_p}+M^{c_p},
			\end{equation*}
			with $\kappa=\frac{p(r-1)}{r(p-1)} \in (0,1)$ provided $r$ close to $1$. Then   \eqref{est:vbesov} follows from
Young's inequality. This completes the proof of Lemma \ref{refksv}.
		\end{proof}
		
		\section{Regularity and propagation}\label{Sect4}
		In this section, we shall establish necessary ingredients for performing Tao's strategy, including pointwise derivative estimates, bounded total speed, epochs of regularity, back propagation, iterated back propagation and annuli of regularity. We shall designate  $\{c_{j,q}\}_{j \in \Z}$ to be a nonnegative generic sequence in $\ell^q(\Z)$ with $\bigl\|(c_{j,q})_{j\in\Z}\bigr\|_{\ell^q(\Z)}\leq 1,$ which may be different on each occurrence below. In particular, we shall denote
$\{c_{j,1}\}_{j \in \Z}$ by $\{d_{j}\}_{j \in \Z}.$

		\begin{lemma}[Point-wise estimate]\label{lem:ptw} Let
$u:[t_0-T,t_0] \times \R^3 \to \R^3$ be a classical solution of \eqref{prob:NS} and  satisfy \eqref{assumption}. Then for any $t\in[t_0-T,t_0]$ and any $j\in\mathbb{Z}$, we have
			\begin{align} \label{est:dyadic-derivative-1}
				&\|\dot\Delta_ju(t)\|_{L^\infty}\lesssim 2^jM^{c_p},\quad \|\nabla\dot\Delta_ju(t)\|_{L^\infty}\lesssim 2^{2j}M^{c_p},\\ \label{est:dyadic-time-1}
				&\|\partial_t\dot\Delta_ju(t)\|_{L^\infty}\lesssim \bigl(2^{3j}+2^{(1+\frac{3}{p})j}(t-(t_0-T))^{-1+\frac{3}{2p}}\bigr)M^{c_p}.
			\end{align}
		\end{lemma}
		\begin{proof}
 Estimates in \eqref{est:dyadic-derivative-1} follow directly from  Lemma \ref{lem:Bernstein} and the assumption \eqref{assumption}. Since
			$$\partial_t\dot\Delta_ju=\Delta \dot\Delta_ju+\mathbb{P}\nabla\cdot\dot\Delta_j(u\otimes u),
			$$
in order to prove \eqref{est:dyadic-time-1},
it suffices to handle the term $\nabla\cdot\dot\Delta_j(u\otimes u)$.
 Yet in view of \eqref{decompv} and \eqref{def:Vi}, we write \begin{equation} \label{u=v+V1}
				u=v+V_1,
			\end{equation}
		so that
			\begin{equation} \label{decompposition_v+V1}
				\dot{\Delta}_j (u \otimes u) = \dot{\Delta}_j (v \otimes v)+ \dot{\Delta}_j (V_1 \otimes v) +  \dot{\Delta}_j (v \otimes V_1) + \dot{\Delta}_j (V_1 \otimes V_1).
			\end{equation}
	By using Bony's decomposition \eqref{form:decomposition}, we have
			\begin{align*} \Delta_j (w \otimes \tilde w)(t)
				=  \dot{\Delta}_j \Bigl( &\sum_{|j'-j| \leq 4} (\dot S_{j'-1}w \otimes \dot\Delta_{j'} \tilde w)(t)  \\
&+ \sum_{|j'-j| \leq 4} (\dot S_{j'-1} \tilde w \otimes \dot\Delta_{j'} w)(t) + \sum_{\ell \geq j-4} (\dot{\Delta}_{\ell} w \otimes \tilde \Delta_\ell \tilde w)(t)\Bigr),
			\end{align*}
			where $w,\tilde w \in \{v,V_1\}$. Then it follows from Lemma \ref{lem:Bernstein}, \eqref{S3eq5} and \eqref{S3eq7} that
			\begin{equation*}
				\begin{aligned}
					\sum_{|j'-j|\leq 4}\|\dot{\Delta}_j(\dot{S}_{j'-1}w\otimes\dot\Delta_{j'} \tilde w)(t)\|_{L^\infty}
&\lesssim \sum_{|j'-j|\leq 4} \Bigl(\sum_{\ell\leq j'-2}\|\dot\Delta_\ell w(t)\|_{L^\infty} \Bigr) \|\dot\Delta_{j'} \tilde w(t)\|_{L^\infty}\\
					&\lesssim 2^{2j}\| w(t) \|_{\dot B_{p,\infty}^{-1+\frac{3}{p}}}\| \tilde w(t)\|_{\dot B_{p,\infty}^{-1+\frac{3}{p}}}\lesssim 2^{2j}M^{c_p}.
				\end{aligned}
			\end{equation*}
For the remaining terms in Bony's decomposition of $v \otimes v$, $v \otimes V_1$ and $V_1 \otimes v$,
we get, by using \eqref{S3eq5} and \eqref{est:vbesov}, that for $w \in \{v,V_1\}$,
			\begin{align*}
\Bigl\|\sum_{l\geq j-4}\dot{\Delta}_j(\dot\Delta_l v \otimes\tilde\Delta_{l} w)(t) \Bigr\|_{L^\infty}
				&\lesssim 2^{3j}\sum_{\substack{\ell\geq j-4}}\|\dot\Delta_\ell v(t)\|_{L^{\frac{p}{p-1}}}\|\tilde\Delta_{\ell} w(t)\|_{L^p}  \\
				&\lesssim 2^{3j}\sum_{\substack{\ell\geq j-4}}2^{\left(\frac{3}{p}-2\right)\ell}\|v(t)\|_{\dot B_{1,\infty}^2}2^{\left(1-\frac{3}{p}\right)\ell}\| w(t) \|_{\dot B_{p,\infty}^{-1+\frac{3}{p}}}\lesssim 2^{2j}M^{c_p}.
			\end{align*}
For the remaining  terms in $V_1 \otimes V_1$, we get, by using \eqref{est:ulindecay} and \eqref{S3eq5}, that
			\begin{align*}
				\Bigl\|\sum_{\substack{\ell\geq j-4}}\dot{\Delta}_j(\dot\Delta_\ell V_1 \otimes\tilde\Delta_{\ell}V_1)(t) \Bigr\|_{L^\infty} &\lesssim 2^{\frac{6}{p}j} \sum_{\substack{\ell\geq j-4}} \|\dot\Delta_\ell V_1(t)\|_{L^{p}}\|\tilde\Delta_{\ell} V_1(t)\|_{L^p}\\
				&\lesssim 2^{\frac{6}{p}j}\sum_{k=1}^m \sum_{\substack{\ell\geq j-4}} 2^{\left(2-\frac{6}{p}\right)\ell}e^{-c_{k-1}(t-(t_0-T))2^{2\ell}}M^{c_p} \\
				&\lesssim 2^{\frac{3}{p}j}(t-(t_0-T))^{-1+\frac{3}{2p}}M^{c_p}.
			\end{align*}

By summarizing the above estimates and using \eqref{decompposition_v+V1}, we obtain \eqref{est:dyadic-time-1}.
	 This completes the proof of Lemma \ref{lem:ptw}.
		\end{proof}

\begin{remark}\label{rmk:ptw} Let $\omega:=\nabla\times u$ be the vorticity of the fluid, then $\omega$ satisfies
		\begin{equation} \label{eq:vorticity}
			\partial_t \omega = \Delta \omega - u \cdot \nabla\omega + \omega \cdot \nabla u.	
		\end{equation}
Furthermore, it is obvious to observe from Lemmas  \ref{lem:Bernstein} and \ref{lem:ptw} that, for any $t \in [t_0-T,t_0]$,
			\begin{align} \label{est:dyadic-vorticity-1}
				&\|\dot\Delta_j \omega(t)\|_{L^\infty}\lesssim 2^{2j}M^{c_p},\quad \|\nabla\dot\Delta_j \omega(t) \|_{L^\infty}\lesssim 2^{3j}M^{c_p},\\ \label{est:dyadic-vorticity-2}
				&\|\partial_t\dot\Delta_j \omega(t)\|_{L^\infty}\lesssim \bigl(2^{4j}+2^{(2+\frac{3}{p})j}(t-(t_0-T))^{-1+\frac{3}{2p}}\bigr)M^{c_p}.
			\end{align}
\end{remark}

Next we show that the total speed of $u$ is bounded. Thanks to the estimates on the components $u_{kL}$ of $u$ (see \eqref{S2eq12}), it is sufficient to deal with the component $v$. For this purpose, we decompose $v$ further into a linear part, a nonlinear part and a force term, and then establish estimates on the homogeneous dyadic blocks of these terms.

		\begin{proposition}[Bounded total speed]\label{S4lem2}  Under the assumptions of Lemma \ref{lem:ptw}, for any interval $I$ in $[t_0-T/2,t_0]$, we have
			\begin{equation} \label{est:L1Linfty}
				\| u \|_{L^1(I; L^\infty)} \lesssim M^{c_p}|I|^{\frac{1}{2}}.	
			\end{equation}
		\end{proposition}
		\begin{proof}
			Without loss of generality, we may assume that $I=[0,1] \subset [t_0-T/2,t_0]$, which implies that $[-1,1] \subset [t_0-T,t_0]$. Restarting the decomposition \eqref{decompv} at time $-1$, with
$u_{1L}(-1)=u(-1)$, $u_{kL}(-1)=0$ for $k\geq2$, and $v(-1)=0$,
we may apply the estimates of Section \ref{sec:basic} with $t+1$
in place of $t$.
			We first deduce from \eqref{condVi} that, for any multi-index $\alpha \in \N^3$ with $0 \leq |\alpha| \leq N_0$ (where $N_0$ is a large positive integer) and $r \geq p$,
			\begin{equation} \label{est:V1}
				\| \nabla^\alpha  V_1(t) \|_{L^r}
				\lesssim (t+1)^{-(\frac{|\alpha|}{2} + \frac{1}{2} - \frac{3}{2r})} M^{c_p}, \quad \forall t \in (-1,1],
			\end{equation}	
			while it follows from \eqref{est:LinftyL2v} that
			\begin{equation} \label{est:v-L2}
				\sup_{t \in (-1,1)}\bigl((t+1)^{-\frac{1}{4}}\| v(t) \|_{L^2}\bigr) +
				\sup_{t \in (-1,1)}\bigl((t+1)^{-\frac{1}{4}}\| \nabla v \|_{L^2(Q_t)}\bigr) \lesssim M^{c_p}.
			\end{equation}
			Applying Plancherel's equality yields
			\begin{equation} \label{est:sum-Dj-unlin}
				\sum_{j \in \Z} 2^{2j} \| \dot{\Delta}_j v \|_{L^2([-1/2,1] \times \R^3)}^2 \leq M^{c_p},
			\end{equation}
		whereas applying	 Sobolev embedding inequality gives
			\begin{equation} \label{est:unlin-L2L6} \| v \|_{L^2([-1/2,1];  L^6)} \lesssim M^{c_p}.
			\end{equation}
Thanks to \eqref{u=v+V1} and \eqref{est:V1}, it suffices to deal with the related estimate of $v$ in order to prove \eqref{est:L1Linfty}.
Indeed in view of \eqref{prob:NS-toy}, we write
\begin{equation} \label{S4eq1} \begin{split}
v(t) = &e^{(t+\frac{1}{2})\Delta} v(-1/2) + v_N(-1/2,t) + v_L(-1/2,t) + v_F(-1/2,t),
			\end{split} \end{equation}  where $ v_N(-1/2,t),$ $v_L(-1/2,t)$ and $v_F(-1/2,t)$ are given by \eqref{form:vF}.
		Below we handle term by term above.

			\noindent \textbf{The estimate of the initial term.} We deduce from \cite[Lemma 2.4]{BahCheDan} and \eqref{S3eq7} that
			\begin{equation}\label{est:vint}
				\begin{aligned}
					\| \dot{\Delta}_j e^{(t+\frac{1}{2})\Delta} v(-1/2) \|_{L^1([0,1]; L^\infty)}
					&\lesssim \Bigl\| e^{-c_02^{2j}(t+\frac{1}{2})} 2^j \| v(-1/2)\|_{\dot{B}_{p,\infty}^{-1+\frac{3}{p}}} \Bigr\|_{L_t^1([0,1])} \\
					&\lesssim 2^j e^{-\frac{c_0}{2}2^{2j}}M^{c_p}\lesssim M^{c_p}d_{j},
				\end{aligned}
			\end{equation}
			where $\{d_{j}\}_{j\in\mathbb{Z}}$ denotes a nonnegative element of $\ell^1(\Z)$ satisfying $\| \{d_{j}\}_{j\in\mathbb{Z}} \|_{\ell^1(\Z)} \leq 1$.

			\noindent \textbf{The estimate of the nonlinear term.} We first get, by applying \cite[Lemma 2.4]{BahCheDan}, that
			\begin{equation}\label{est:vnlin-a}
				\begin{aligned}
					\| \dot{\Delta}_j v_N(-1/2,t) \|_{L^1([0,1]; L^\infty)}
					&\lesssim \int_0^1 \Bigl( \int_{-1/2}^t e^{-c_02^{2j}(t-\tau)}  d\tau  \Bigr) 2^j \| \dot \Delta_j (v \otimes v)(\tau) \|_{L^\infty} \,dt \\
					&\lesssim 2^{-j} \| \dot \Delta_j (v \otimes v) \|_{L^1([-1/2,1]; L^\infty)}.
				\end{aligned}
			\end{equation}
	By using Bony's decomposition, we write
			\begin{align} \label{decompose-v}
				\dot{\Delta}_j(v \otimes v) =2\sum_{|j-j'|\leq 4}\dot{\Delta}_j (\dot{S}_{j'-1} v \otimes \dot{\Delta}_{j'} v) + \sum_{\ell\geq j-4} \dot{\Delta}_j(\dot{\Delta}_{\ell} v\otimes \tilde {\Delta}_\ell v),
			\end{align}
			where $\tilde \Delta_j = \dot{\Delta}_{j+1} + \dot{\Delta}_j + \dot{\Delta}_{j-1}$.

It follows from
			 Lemma \ref{lem:Bernstein} and \eqref{est:sum-Dj-unlin} that
			\begin{equation*}
				\begin{aligned}
					&\sum_{|j-j'|\leq 4} \| \dot{\Delta}_j(\dot{S}_{j'-1} v \otimes \dot{\Delta}_{j'} v) \|_{L^1([-1/2,1]; L^\infty)}\\
					&\lesssim \sum_{|j-j'|\leq 4}  \sum_{\ell\leq j'-2}\| \dot{\Delta}_\ell v \|_{L^2([-1/2,1]; L^\infty)}
\|\dot{\Delta}_{j'} v\|_{L^2([-1/2,1]; L^\infty)}\\
					&\lesssim c_{j,2}^22^j\|\nabla v\|^2_{L^2([-1/2,1] \times \R^3)} \lesssim d_{j} 2^j M^{c_p}.
				\end{aligned}
			\end{equation*}
Similarly we get, by using  Lemma \ref{lem:Bernstein} and \eqref{est:sum-Dj-unlin}, that
			\begin{align*}
				\Bigl\|\sum_{\ell\geq j-4} \dot{\Delta}_j(\dot{\Delta}_{\ell} v\otimes \tilde{\Delta}_\ell v) \Bigr\|_{L^1([-1/2,1]; L^\infty)}\lesssim& 2^{3j}\sum_{\ell\geq j-4} \|\dot{\Delta}_{\ell} v\otimes \tilde{\Delta}_\ell v\|_{L^1([-1/2,1]; L^1)}\\
				\lesssim& c_{j,2}^2 2^j\|\nabla v\|^2_{L^2([-1/2,1] \times \R^3)} \lesssim d_{j} 2^j M^{c_p}.
			\end{align*}
By substituting the above estimates into \eqref{est:vnlin-a}, we obtain
			\begin{equation} \label{est:vnlin}
				\| \dot{\Delta}_j v_N(-1/2,t) \|_{L^1([-1/2,1]; L^\infty)}\lesssim d_{j}M^{c_p}.
			\end{equation}
			
			\noindent \textbf{The estimate of the linear terms.} Notice that $\dot{\Delta}_j = \dot{\Delta}_j \tilde \Delta_j$, hence in view of
\eqref{form:vF},
 we get, by applying Lemma \ref{lem:Bernstein}, \cite[Lemma 2.4]{BahCheDan}, \eqref{est:V1}, and \eqref{est:unlin-L2L6}, that
			\begin{equation} \label{est:vlin}
				\begin{aligned}
					&\| \dot{\Delta}_j v_L(-1/2,t) \|_{L^1([-1/2,1]; L^\infty)}\\
					&\lesssim 2^{j}\min\{2^{-2j},1\} \bigl(\| \tilde \Delta_j (v \otimes V_1) \|_{L^1([-1/2,1]; L^\infty)}
 + \| \tilde \Delta_j (V_1 \otimes v) \|_{L^1([-1/2,1]; L^\infty)}\bigr) \\
					&\lesssim \min\bigl\{2^{-\frac{j}{2}}, 2^{\frac{3j}{2}}\bigr\}
\bigl(\| v \otimes V_1\|_{L^1([-1/2,1]; L^6)} + \| V_1 \otimes v\|_{L^1([-1/2,1]; L^6)} \bigr) \\
					&\lesssim \min\bigl\{2^{-\frac{j}{2}}, 2^{\frac{3j}{2}}\bigr\} \| v \|_{L^2([-1/2,1]; L^6)}
 \| V_1 \|_{L^2([-1/2,1]; L^\infty)} \lesssim    d_{j} M^{c_p}.
				\end{aligned}
			\end{equation}

			\noindent \textbf{The estimate of the force term.} Similar to \eqref{est:vlin}, we have
\begin{equation} \label{S4eq2}
			\begin{aligned}
				&\| \dot{\Delta}_j v_F(-1/2,t) \|_{L^1([0,1]; L^\infty)} \\
				&\lesssim 2^{j} \min\bigl\{2^{-2j},1\bigr\} \bigl(\| \tilde \Delta_j (V_3 \otimes V_1) \|_{L^1([-1/2,1]; L^\infty)}+ \| \tilde \Delta_j (V_2 \otimes V_3) \|_{L^1([-1/2,1]; L^\infty)}\bigr).
			\end{aligned}\end{equation}
		Yet	by using H\"older's inequality and the estimates
			\begin{align*}
				\| V_1(t) \|_{L^\infty} + \| V_2(t) \|_{L^\infty} + \| V_3(t) \|_{L^6} \lesssim M^{c_p}, \quad \forall t \in [0,1],
			\end{align*}
			we find
			\begin{align*}
				&\| \tilde \Delta_j (V_3 \otimes V_1) \|_{L^1([-1/2,1]; L^\infty)} + \| \tilde \Delta_j (V_2 \otimes V_3) \|_{L^1([-1/2,1]; L^\infty)} \\
				&\lesssim 2^{\frac{j}{2}} \bigl(\| V_3 \otimes V_1\|_{L^1([-1/2,1]; L^6)} + \| V_2 \otimes V_3\|_{L^1([-1/2,1]; L^6)}\bigr) \\
				&\lesssim 2^{\frac{j}{2}} \| V_3 \|_{L^2([-1/2,1]; L^6)} \bigl(\| V_1 \|_{L^2([-1/2,1]; L^\infty)} +\| V_2 \|_{L^2([-1/2,1]; L^\infty)}\bigr) \lesssim    2^{\frac{j}{2}}M^{c_p}.
			\end{align*}
			By inserting the above estimates into \eqref{S4eq2}, we achieve
			\begin{equation}\label{est:vforce}
				\| \dot{\Delta}_j v_F(-1/2,t) \|_{L^1([0,1]; L^\infty)} \lesssim \min\bigl\{2^{-\frac{j}{2}}, 2^{\frac{3j}{2}}\bigr\}  M^{c_p}\lesssim d_{j}M^{c_p}.
			\end{equation}

	In view of \eqref{S4eq1}, we get, by summarizing the estimates \eqref{est:vint}, \eqref{est:vnlin}, \eqref{est:vlin} and \eqref{est:vforce} that
			$$ \| v \|_{L^1([0,1]; L^\infty)}\lesssim M^{c_p},
			$$
			which, together with  \eqref{est:V1}, ensures
			$$
			\| u \|_{L^1([0,1];L^\infty)} \lesssim M^{c_p}.
			$$
			This completes the proof of Proposition \ref{S4lem2}.	
		\end{proof}
		
	Next we adapt the idea in \cite{Tao_20} to establish epochs of regularity for $u$. The proof relies on a careful analysis of the enstrophy-type quantity for $v$ and is more intricate than that in \cite{Tao_20} due to the involvement of several additional terms. By combining estimates on the linear components $u_{kL}$, $k=1,\ldots,m$, and on the Oseen kernel $\Gamma$ (see \eqref{S2eq13}), we are able to handle all the involved terms, thereby obtaining the desired bounds.
		
		\begin{proposition}[Epochs of regularity] \label{epochs}  Under the assumptions of Lemma \ref{lem:ptw},
			 for any interval $I \subset [t_0-T/2,t_0]$ there is a subinterval $I' \subset I$ with
			 $cM^{-c_p}|I|\leq |I'|\leq M^{-c_p}|I|$ such that for any $\al=(\al_1,\al_2,\al_3)$ with $|\al|=0,1$, we have
			\begin{equation} \label{est:grad-u} \| \nabla^\al u \|_{L^\infty (I' \times \R^3)} \lesssim M^{c_p}|I|^{-\frac{|\alpha|+1}{2}},
			\end{equation}
			and
			\begin{equation} \label{est:grad-omega} \| \nabla^\alpha \omega \|_{L^\infty (I' \times \R^3)} \lesssim M^{c_p}|I|^{-\frac{|\alpha|+2}{2}},
			\end{equation}
		where $\omega = \nabla \times u$ is the vorticity of the fluid (see Remark \ref{rmk:ptw}).
		\end{proposition}
		\begin{proof}
			By rescaling and time translation, we may assume that $I=[0,1]$ and $[-1,1] \subset [t_0-T,t_0]$. Below we divide the proof
into the following steps:
			
			\noindent \textbf{The estimate of $u$.} Let $v$ solve \eqref{prob:NS-toy} on $[t_0-T,t_0].$
			We define the enstrophy-type quantity
			$$
			\cE(t) := \frac{1}{2} \int_{\R^3} |\nabla v(t,x)|^2 dx.
			$$
			Then we get, by taking the gradient of \eqref{prob:NS-toy} and then taking $L^2$  inner product of the resulting equation with $\nabla v,$ that
			\begin{align*} \partial_t \cE(t) = &- \int_{\R^3} |\nabla^2 v(t)|^2 dx + \int_{\R^3} \Delta v(t) \nabla \cdot (v \otimes v)(t) dx \\
&+ \int_{\R^3} \Delta v(t)  (\nabla \cdot (v \otimes V_1 + V_1 \otimes v)(t) dx + 	\int_{\R^3} \Delta v(t)  \nabla \cdot (V_3 \otimes V_1 + V_2 \otimes V_3)(t) dx,
			\end{align*}
			where $V_i$, $i=1,2,3$, are given by \eqref{def:Vi}. By applying Young's inequality, we find that
\begin{equation}\label{S4eq3}
			\begin{split}
				\partial_t \cE(t) &\leq  - \frac{1}{2} \| \nabla^2 v(t) \|_{L^2}^2 + 16 \bigl(\| \nabla \cdot (v \otimes v)(t) \|_{L^2}^2 +  \|  (V_1 \cdot \nabla v + v \cdot \nabla V_1)(t) \|_{L^2}^2 \\
				&\quad  + \|  (V_1 \cdot \nabla V_3 + V_3 \cdot \nabla V_2)(t) \|_{L^2}^2\bigr).
			\end{split}\end{equation}
It is easy to observe that
	\begin{align*}
 \| \nabla \cdot (v \otimes v)(t) \|_{L^2} \lesssim& \| v(t) \|_{L^6} \| \nabla v(t) \|_{L^3}\\
 \lesssim &\| \nabla v(t) \|_{L^2}\| \nabla v(t)\|_{L^2}^{\frac{1}{2}} \| \nabla^2 v(t)\|_{L^2}^{\frac{1}{2}}
 \lesssim \cE(t)^{\frac{3}{4}}\| \nabla^2 v(t)\|_{L^2}^{\frac{1}{2}}.
 \end{align*}
	Recalling that $V_i$, $i=1,2,3$,  satisfy \eqref{condVi}.	By using Gagliardo-Nirenberg inequality, for $t \in [0,1]$, one has
			\begin{align*}
				\|  (V_1 \cdot \nabla v + v \cdot \nabla V_1)(t) \|_{L^2} &\leq \| V_1(t) \|_{L^p} \| \nabla v(t) \|_{L^{\frac{2p}{p-2}}} + \| v(t) \|_{L^{\frac{2p}{p-2}}} \| \nabla V_1(t) \|_{L^p}	\\
				& \lesssim \| V_1(t) \|_{L^p} \| \nabla v(t) \|_{L^2}^{\frac{p-3}{p}} \| \nabla^2 v(t) \|_{L^2}^{\frac{3}{p}} + \|v(t)\|_{L^2}^{\frac{p-3}{p}} \| \nabla v(t) \|_{L^2}^{\frac{3}{p}} \| \nabla V_1(t)\|_{L^p} \\
				&\lesssim  M^{c_p} \cE(t)^{\frac{p-3}{2p}}\| \nabla^2 v(t) \|_{L^2}^{\frac{3}{p}} + M^{c_p} \cE(t)^{\frac{3}{2p}},
			\end{align*}
			and
			\begin{align*}
				\|  (V_1 \cdot \nabla V_3 + V_3 \cdot \nabla V_2)(t) \|_{L^2} &\leq 	 \| V_1(t) \|_{L^p} \| \nabla V_3(t) \|_{L^{\frac{2p}{p-2}}} + \| V_3(t) \|_{L^{\frac{2p}{p-2}}} \| \nabla V_2(t) \|_{L^p} \lesssim M^{c_p}. 
				&
			\end{align*}
			By substituting the above estimates into \eqref{S4eq3} and using Young's inequality, we obtain
			\begin{equation} \label{est:E-1}
				\partial_t \cE(t) + \| \nabla^2 v\|_{L^2}^2 \lesssim \cE(t)^{3} + M^{c_p} \cE(t) + M^{c_p}.
			\end{equation}
		On the other hand, notice from \eqref{est:v-L2}  (we may assume $c_p$ large enough) that
			$ \int_0^1 \cE(t) dt \lesssim M^{c_p}.
			$
			So that it follows from the pigeonhole principle that there exists a time $t_1 \in [0,1/2]$ with $\cE(t_1) \lesssim M^{c_p}$.
Then we deduce from \eqref{est:E-1} that $\cE(t) \lesssim M^{c_p}$ for $t \in [t_1,t_1 + cM^{-3c_p}] = [\tau(0),\tau(1)]$ where $\tau(\sigma) = t_1 + cM^{-3c_p}\sigma$ and $c>0$ is a sufficiently small constant. By inserting this estimate into \eqref{est:E-1}, we obtain
			$$ \partial_t \cE(t) + \| \nabla^2 v(t)\|_{L^2}^2 \lesssim M^{3c_p}, \quad t \in [\tau(0),\tau(1)],
			$$
	which implies
			\begin{equation} \label{est:nabla-D2-u(m+1)} \| \nabla v\|_{L^\infty([\tau(0),\tau(1)]; L^2)} +  \| \nabla^2 v\|_{L^2([\tau(0),\tau(1)] \times \R^3)} \lesssim M^{c_p}.
			\end{equation}
		Applying the Gagliardo-Nirenberg inequality and  the Sobolev inequality  yields
			$$
			\|v\|_{L^4([\tau(0),\tau(1)]; L^\infty)}\lesssim M^{c_p},
			$$
		which, together with \eqref{est:nabla-D2-u(m+1)}, ensures that
			\begin{equation} \label{est:nablav-L2L6} \| v \|_{L^4([\tau(0),\tau(1)];  L^\infty)} \lesssim M^{c_p} \andf
 \| \nabla v \|_{L^2([\tau(0),\tau(1)];  L^6)} \lesssim M^{c_p}.
			\end{equation}
In view of \eqref{prob:NS-toy},  for $t \in [\tau(0.1), \tau(1)]$, we write
		\begin{equation}	\label{S4eq4}
				v(t) = e^{(t-\tau(0))\Delta}v(\tau(0))+v_N(\tau(0),t)+v_L(\tau(0),t)+v_F(\tau(0),t),
			\end{equation}
where $ v_N(\tau(0),t),$ $v_L(\tau(0),t)$ and $v_F(\tau(0),t)$ are given by \eqref{form:vF}.

			We first observe that
			\begin{align*}
				\| e^{(t-\tau(0))\Delta}v(\tau(0)) \|_{L^\infty([\tau(0.1),\tau(1)]\times \R^3)}
				&\lesssim (\tau(0.1) -\tau(0))^{-\frac{1}{4}}\|v(\tau(0))\|_{L^6}  \lesssim M^{c_p}.
			\end{align*}
By using Oseen kernel $\Gamma(t,x)$ (see \eqref{S2eq13} below), we find
\begin{align*}
\|v_N(\tau(0),t)\|_{L^\infty}\leq \int_{\tau(0)}^t \|\Gamma(t-s,\cdot)\|_{L^1}\|(v \otimes v)(s)\|_{L^\infty}\,ds\lesssim
\int_{\tau(0)}^t(t-s)^{-\frac12}\|v(s)\|_{L^\infty}^2\,ds,
\end{align*}
			which implies
			\begin{align*}
	\|v_N(\tau(0),t)\|_{L^8([\tau(0.1),\tau(1)];L^\infty)}		
				&\lesssim (\tau(1)-\tau(0))^{\frac{1}{8}}M^{c_p} \lesssim M^{c_p}.
			\end{align*}
The last two terms in \eqref{S4eq4} share the same estimate.
		
As a consequence, we obtain
			$$
			\| v \|_{L^8([\tau(0.1),\tau(1)];  L^\infty)} \lesssim M^{c_p}.	
			$$
			This together with \eqref{u=v+V1} and \eqref{est:V1} ensures that
			$$
			\| u \|_{L^8([\tau(0.1),\tau(1)];  L^\infty)}\lesssim M^{c_p}.	
			$$
		By	repeating the above process, we obtain for $t \in [\tau(0.2),\tau(1)]$,
			\begin{equation} \label{est:u,v-LinftyLinfty} \|(v,u)\|_{L^\infty([\tau(0.2),\tau(1)] \times \R^3)} \lesssim M^{c_p}.	
			\end{equation}
			
			\noindent \textbf{The estimate of $\nabla u$.} We shall prove in this step that
			\begin{equation} \| \nabla v \|_{L^\infty([\tau(0.4),\tau(1)] \times \R^3)} \lesssim M^{c_p}.	
			\end{equation}
In order to do so, in view of \eqref{prob:NS-toy},
			for $t \in [\tau(0.3),\tau(1)]$, we write
			\begin{align*}
			v(t) &= e^{(t-\tau(0.2))\Delta}v(\tau(0.2)) +v_N(\tau(0.2),t)+v_L(\tau(0.2),t)+v_F(\tau(0.2),t),
			\end{align*}
where $ v_N(\tau(0.2),t),$ $v_L(\tau(0.2),t)$ and $v_F(\tau(0.2),t)$ are given by \eqref{form:vF}.\\
We first observe from \eqref{est:u,v-LinftyLinfty} that
			\begin{align*}
				\| \nabla e^{(t-\tau(0.2))\Delta}v(\tau(0.2))  \|_{L^\infty} &\lesssim (t-\tau(0.2))^{-\frac{1}{2}} \| v(\tau(0.2))\|_{L^\infty} \\
				&\lesssim (\tau(0.3) - \tau(0.2))^{-\frac{1}{2}}M^{c_p} \lesssim M^{c_p}.
			\end{align*}
By using Oseen Kernel estimate \eqref{S2eq13},	 \eqref{est:nablav-L2L6} and \eqref{est:u,v-LinftyLinfty}, we find
			\begin{align*}
				 \|\nabla v_N(\tau(0.2),t)\|_{L^4([\tau(0.3),\tau(1)]; L^\infty)}
				&\lesssim \Bigl\| \int_{\tau(0.2)}^t (t-s)^{-\frac{3}{4}} \| \nabla \cdot (v \otimes v)(s) \|_{L^6}ds \Bigr \|_{L^4([\tau(0.3),\tau(1)] )}\\
				&\lesssim  \| \nabla \cdot (v \otimes v) \|_{L^2([\tau(0.2),\tau(1)];  L^6)} \\
				&\lesssim \| \nabla v \|_{L^2([\tau(0),\tau(1)];  L^6)} \| v \|_{L^\infty([\tau(0.2),\tau(1)] \times \R^3)} \lesssim M^{c_p}.
			\end{align*}
			Similarly, we get, by  using \eqref{est:unlin-L2L6}, \eqref{est:nablav-L2L6} and the estimate on $V_i$, $i=1,2,3$, in \eqref{condVi}, that
			\begin{align*}
			 \|&\nabla v_L(\tau(0.2),t)\|_{L^4([\tau(0.3),\tau(1)]; L^\infty)}\\
				&\lesssim \Bigl\| \int_{\tau(0.2)}^t (t-s)^{-\frac{3}{4}} \| (V_1 \cdot \nabla v + v \cdot \nabla V_1)(s) \|_{L^6}ds \Bigr \|_{L_t^4([\tau(0.3),\tau(1)] )}\\
				&\lesssim \| V_1 \|_{L^\infty([\tau(0.2),\tau(1)] \times \R^3)}\| \nabla v \|_{L^2([\tau(0),\tau(1)]; L^6)}\\
&\quad + \| v \|_{L^2([\tau(0),\tau(1)]; L^6)}\| \nabla V_1 \|_{L^\infty ([\tau(0.2),\tau(1)] \times \R^3)} \lesssim M^{c_p},
			\end{align*}
			and
			\begin{align*}
				\|&\nabla v_F(\tau(0.2),t)\|_{L^4([\tau(0.3),\tau(1)]; L^\infty)}\\
				&\lesssim \Bigl\| \int_{\tau(0.2)}^t (t-s)^{-\frac{3}{4}} \| (V_1 \cdot \nabla V_3 + V_3 \cdot \nabla V_2)(s) \|_{L^6}ds \Bigr \|_{L_t^4([\tau(0.3),\tau(1)] )}\\
				&\lesssim \| V_1 \|_{L^\infty([\tau(0.2),\tau(1)] \times \R^3)}\| \nabla V_3 \|_{L^2([\tau(0),\tau(1)];  L^6)}  \\
				&\quad + \| V_3 \|_{L^2([\tau(0),\tau(1)];  L^6)} \| \nabla V_2 \|_{L^\infty([\tau(0.2),\tau(1)] \times \R^3)} \lesssim M^{c_p}.
			\end{align*}
		As a result, it comes out
			$$
			\| \nabla v \|_{L^4([\tau(0.3),\tau(1)]; L^\infty)} \lesssim M^{c_p},
			$$
which, together with \eqref{u=v+V1} and \eqref{est:V1}, ensures that		
		\begin{equation} \label{S4eq7}
			\| \nabla u \|_{L^4([\tau(0.3),\tau(1)]; L^\infty)} \lesssim M^{c_p}.
			\end{equation}
For $t \in [\tau(0.4),\tau(1)]$, we have
			$$ \nabla u(t) = \nabla e^{(t-\tau(0.3))\Delta}u(\tau(0.3)) - \int_{\tau(0.3)}^t \nabla e^{(t-s)\Delta} \bP \nabla \cdot (u \otimes u)(s)ds.
			$$
It is obvious that
			\begin{align*}
				\| \nabla e^{(t-\tau(0.3))\Delta}u(\tau(0.3))  \|_{L^\infty} \lesssim M^{c_p}.
			\end{align*}
By proceeding along the same line to the estimate of $\|\nabla v \|_{L^4([\tau(0.3),\tau(1)]; L^\infty)},$ we infer		
			\begin{align*}
		\Bigl\|& \int_{\tau(0.3)}^t \nabla e^{(t-s)\Delta} \bP \nabla \cdot (u \otimes u)(s)ds \Bigr\|_{L^\infty}
				\leq 	\int_{\tau(0.3)}^t (t-s)^{-\frac{1}{2}} \| \nabla \cdot (u \otimes u)(s)\|_{L^\infty}ds \\
				&\lesssim \Bigl( \int_{\tau(0.3)}^t (t-s)^{-\frac{2}{3}}ds \Bigr)^{\frac{3}{4}} \| \nabla \cdot (u \otimes u) \|_{L^4([\tau(0.3),\tau(1)]; L^\infty)} \lesssim M^{c_p}.
			\end{align*}	
			Therefore, we achieve
			$$ \| \nabla u\|_{L^\infty ([\tau(0.4),\tau(1)] \times \R^3)} \lesssim M^{c_p}.
			$$
			
			\noindent \textbf{The estimates of $\omega$.} By a similar argument as in the proof of Proposition (iii) in \cite{Tao_20}, we
can prove \eqref{est:grad-omega}. The detail is omitted. We thus complete the proof of Proposition \ref{epochs}.
		\end{proof}
		
		Thanks to estimates on $u_{kL}$, $k=1,\ldots,m$, and $v$, we can also adapt the idea in \cite{Tao_20} to control the homogeneous dyadic blocks with high and intermediate frequencies for the linear part, nonlinear part and the force term appearing in the decomposition of $v$, and hence show the back propagation phenomenon. In the following lemma, we will use the constants $M_i$ which satisfy \eqref{Mi}.
		\begin{proposition}[Back propagation] \label{backprop}
			Let $u:[t_0-T,t_0] \times \R^3 \to \R^3$ be a classical solution of \eqref{prob:NS} and satisfy \eqref{assumption}.
Let $(t_1,x_1) \in [t_0-T/2,t_0] \times \R^3$ and $2^{j_1} \geq M_3 T^{-1/2}$ be such that
			\begin{equation} \label{est:backprop-1}
				|\dot{\Delta}_{j_1}u(t_1,x_1)| \geq M_1^{-1}2^{j_1}.	
			\end{equation}
			Then there exist $(t_2,x_2) \in [t_0-T,t_1] \times \R^3$ and $j_2 \in \Z$ with $2^{j_2} \in [M_2^{-1}2^{j_1},M_2 2^{j_1}]$ such that
		\begin{equation}\label{est:backprop-2}
\begin{split}
& M_3^{-1}2^{-2j_1} \leq t_1 - t_2 \leq M_3 2^{-2j_1},\  |x_2 - x_1| \leq M_4 2^{-j_1}, \
			 |\dot{\Delta}_{j_2}u(t_2,x_2)| \geq M_1^{-1}2^{j_2}.	
\end{split}
			\end{equation}
		\end{proposition}
		\begin{proof}
			We may normalize $j_1=0$ and $(t_1,x_1)=(0,0)$	so that $t_0 - T \leq -T/2 \leq -M_{3}^2/2$. In particular,
we have $[-2M_3,0] \subset [t_0 - T, t_0]$, and the assumption \eqref{est:backprop-1} reduces to
			\begin{equation} \label{assump:contradition} |\dot{\Delta}_0 u(0,0)| \geq M_1^{-1}.
			\end{equation}
Let us suppose by contradiction that \eqref{est:backprop-2} is not correct. Then we have
			\begin{equation} \label{supposition} \| \dot{\Delta}_{j} u \|_{L^\infty([-M_3,-M_3^{-1}] \times B(0,M_4))} \leq M_1^{-1}2^j,
			\end{equation}
			for all $j$ such that $2^j \in [M_2^{-1},M_2]$.\\
It follows from  Lemma \ref{lem:ptw}  that for any $t \in [-M_3^{-1},0]$, $x \in B(0,M_4)$ and $2^j\in[M_2^{-1},M_2]$
			\begin{align*}
				|\dot{\Delta}_j u(t,x)| &\leq |\dot{\Delta}_j u(-M_3^{-1},x)| + \int_{-M_3^{-1}}^t |\partial_t \dot{\Delta}_j u(\tau,x)|d\tau \lesssim M_1^{-1} 2^j,
			\end{align*}
that is
			\begin{equation*} \label{est:Deltaj-u-LinftyLinfty}
				\| \dot{\Delta}_{j} u \|_{L^\infty([-M_3,0] \times B(0,M_4))} \leq  M_1^{-1}2^j.
			\end{equation*}
We observe from \eqref{est:ulindecay} that for $2^j \in [M_2^{-1}, M_2]$ and $t \in [-M_3,0]$ that
			$$ \| \dot\Delta_j V_1(t) \|_{ L^\infty(\R^3)} \lesssim \sum_{k=1}^m 2^j e^{-c_{k-1}(t+2M_3)2^{2j}}M^{2^{k-1}} \lesssim  M_1^{-1}2^j.
			$$
		We	thus  obtain for $2^j \in [M_2^{-1},M_2]$ that
			\begin{equation} \label{est:Deltaj-v-LinftyLinfty}
				\| \dot{\Delta}_{j} v \|_{L^\infty ([-M_3,0] \times B(0,M_4))} \leq  M_1^{-1}2^j.
			\end{equation}
			The rest of the proof is split into several steps.
			
			\noindent {\sc Step 1.} The case for $2^j \geq M_2^{-1}$.

In what follows, we shall make a small modification of notations. We denote $u_{1L}:=e^{(t+2M_3)\Delta}u(-2M_3)),$
and similar modification for $u_{iL},$ $ 1<i\leq m$ (recall that $m=[p]+3$) defined by the equations \eqref{prob:lin} on $[-2M_3,\infty)$ with vanishing
initial data at  ${t=-2M_3}$. Let $v$ be determined by \eqref{prob:NS-toy} with initial data $v|_{t=-2M_3}=0$, and $V_1,V_2,V_3$ be the same as \eqref{def:Vi}. Then  \eqref{u=v+V1} still holds. Moreover, along the same line to the proof of \eqref{condVi} and \eqref{est:ulindecay}, we have
			\begin{equation}\label{est:Katonew}
				\begin{aligned}
					&\sum_{i=1}^3\sup_{\substack{t\in[-2M_3,0]\\ 0\leq |\al|\leq N_0, q\in[p,\infty]}}(t+2M_3)^{\frac{|\alpha|}{2}+\frac{1}{2}\left(1-\frac{3}{q}\right)}\|\nabla^\al V_i(t)\|_{L^q}\lesssim M^{c_p},\\
					&\|\dot\Delta_ju_{kL}(t)\|_{L^r}\lesssim 2^{(1-\frac{3}{r})j}e^{-c_{p}(t+2M_3)2^{2j}}M^{c_{p}},\quad 1\leq k\leq m,\quad r\geq p.
				\end{aligned}
			\end{equation}
			For $v$, the estimates \eqref{est:v-Besov} and \eqref{est:Deltaj-v-LinftyLinfty} still hold, and it follows from  the proof of \eqref{est:LinftyL2v} and Lemma \ref{refksv} that
			\begin{equation}\label{LinftyL2vnew}
				\begin{aligned}
				&\sup_{t\in(-2M_3,0]}(t+2M_3)^{-\frac{1}{2}}\|v(t)\|^2_{L^2}+\sup_{t\in(-2M_3,0]}(t+2M_3)^{-\frac{1}{2}}\|\nabla v\|^2_{L^2([-2M_3,t]\times \mathbb{R}^3)} \lesssim M^{c_p},\\
					&\|\dot{\Delta}_jv\|_{L^\infty([-2M_3,0];  L^{1}(B(0,M_4)))}\lesssim 2^{-2j} \|v\|_{L^\infty ([-2M_3,0]; \dot B_{1,\infty}^2)} \lesssim 2^{-2j}M^{c_p}.
				\end{aligned}
			\end{equation}
			To estimate $V_1$, we get, by using  \eqref{est:Katonew}, that
			\begin{align*}
				\| \dot{\Delta}_j V_1(t)\|_{L_x^{1}(B(0,M_4))} & = \sum_{k=1}^m \| \dot{\Delta}_j u_{kL}(t) \|_{L_x^{1}(B(0,M_4))}\lesssim \sum_{k=1}^m \| \1_{B(0,M_4)}\|_{L_x^{\frac{p}{p-1}}(\mathbb{R}^3)} \| \dot{\Delta}_j u_{iL} (t)\|_{ L_x^p( \mathbb{R}^3)} \\
				&\lesssim \sum_{k=1}^m M_4^{\frac{3(p-1)}{p}} 2^{j\frac{p-3}{p}} e^{-c_{k-1}  (t+2M_3) 2^{2j}} M^{c_p}\\
&\lesssim \sum_{k=1}^m\bigl(M_4^{\frac{3(p-1)}{p} }e^{-\frac{c_{k-1}}{2}M_3M_2^{-2}  }\bigr) \bigl(2^{(3-\frac{3}{p})j}e^{-\frac{c_{k-1}}{2}M_32^{2j}  }\bigr)2^{-2j}M^{c_p} \\
				&\lesssim 2^{-2j}M^{c_p} , \quad \forall t \in [-M_3,0].
			\end{align*}
where we used \eqref{Mi}, \eqref{est:ulindecay} and $2^{j}\geq M_2^{-1}$ in the last line. This, together
  with the second inequality of \eqref{LinftyL2vnew} and \eqref{u=v+V1}, ensures that, for $2^j \geq M_2^{-1}$,
			\begin{equation}\label{est:2^j>M^{-1}}
				\|\dot{\Delta}_ju\|_{L^\infty([-M_3,0];  L^{1}(B(0,M_4)))}\lesssim 2^{-2j}M^{c_p}.
			\end{equation}

			\noindent {\sc Step 2.} The case for $M_2^{-\frac{1}{2}} \leq 2^j \leq M_2^{\frac{1}{2}}$.

	Once again we shall use the decomposition \eqref{u=v+V1}.	By \eqref{est:Katonew}, we observe that
			\begin{equation} \label{est:LinftyL4-V1}
				\begin{aligned}
					\|\dot{\Delta}_j V_1 \|_{L^\infty([-M_3/2,0];  L^{4}(B(0,9M_1)))}
&\lesssim \sum_{k=1}^m \|\dot{\Delta}_j u_{kL}\|_{L^\infty([-M_3/2,0];  L^{4}(B(0,9M_1)))}  \\
					&\lesssim M_1^{\frac{3}{4}}\sum_{k=1}^m\|\dot\Delta_ju_{kL}\|_{L^\infty([-M_3/2,0] \times \mathbb{R}^3)}\\
 &\lesssim M_1^{\frac{3}{4}}M_2^{\frac12}e^{-c_pM_3M_2^{-\frac12}}M^{c_p}
\lesssim M_4^{-50c_p}.
				\end{aligned}
			\end{equation}

Similarly to \eqref{S4eq4}, we write
		\begin{equation}	\label{S4eq8}
				v(t) =	e^{(t+M_3)\Delta} v(-M_3)+v_N(-M_3,t)+v_L(-M_3,t)+v_F(-M_3,t),
			\end{equation}
where $ v_N(-M_3,t),$ $v_L(-M_3,t),$ $v_F(-M_3,t)$ are given by \eqref{form:vF}.

Below we handle term by term in \eqref{S4eq8}.

			\noindent \textbf{The  estimate of  $\dot{\Delta}_j e^{(t+M_3)\Delta} v(-M_3)$.}
We get, by applying
H\"older's inequality, \cite[Lemma 2.4]{BahCheDan} and \eqref{est:v-Besov}, that for $2^j \in [M_2^{-\frac{1}{2}},M_2^{\frac{1}{2}}]$ and for any $t \in [-M_3/2,0]$,
			\begin{align*}
				\| \dot{\Delta}_j  e^{(t+M_3)\Delta} v(-M_3)\|_{L^4(B(0,8M_1))} &\leq M_1^{\frac{3}{4}}\| e^{(t+M_3)\Delta} \dot{\Delta}_j  v(-M_3) \|_{L^\infty(\R^3)} \\
				&\lesssim M_1^{\frac{3}{4}} e^{-c_0(t+M_3)2^{2j}}2^{\frac{3}{p}j}\| \dot{\Delta}_j v(-M_3)\|_{L^p(\mathbb{R}^3)} \\
				&\lesssim M_1^{\frac{3}{4}}e^{-\frac{1}{2}c_0M_3M_2^{-1}}M_2^{\frac{1}{2}}M^{c_p} \lesssim M_4^{-50c_p},
			\end{align*}	
			which implies
			\begin{equation}\label{est:LinftyL4v0}
				\begin{aligned}
					\|\dot{\Delta}_j e^{(t+M_3)\Delta} v(-M_3)\|_{L^\infty([-M_3/2,0]; L^{4}(B(0,8M_1)))}\lesssim M_4^{-50c_p}.
				\end{aligned}
			\end{equation}	

			\noindent \textbf{The estimate of $\dot{\Delta}_j v_N(-M_3,t)$.} In view of \eqref{form:vF},
we deduce from
			Lemma \ref{localY}  that, for any $q\in[1,\infty]$
			\begin{equation}\label{est:Duhamel}
\begin{split}
				&\|\dot{\Delta}_j v_N(t)\|_{L_x^q(B(0,8M_1))} \\
				&\lesssim \int_{-M_3}^t 2^je^{-c_0(t-s)2^{2j}}\bigl(\|\dot\Delta_j(v\otimes v)(s)\|_{L_x^q(B(0,9M_1))}
+M_4^{-50c_p}\|v(s)\|^2_{L^2} \bigr)ds.
\end{split}
			\end{equation}
We use the decomposition \eqref{decompose-v} to deal with the estimate of $\dot\Delta_j(v\otimes v).$
		For the para-product part, we get, by using Lemma \ref{localY}, that
			\begin{align*}
				&\sum_{|j-j'|\leq 4}\|\dot{\Delta}_j(\dot S_{j'-1}v\otimes\dot{\Delta}_{j'}v)\|_{L^\infty([-M_3,0];  L^{4}(B(0,9M_1)))}\\
				&\lesssim \sum_{|j-j'|\leq 4}\bigl( \|\dot S_{j'-1}v\otimes\dot{\Delta}_{j'}v\|_{L^\infty([-M_3,0];  L^{4}(B(0,10M_1)))} + M_4^{-50c_p}\|\dot S_{j'-1}v\otimes \dot\Delta_{j'}v\|_{L^\infty([-M_3,0];  L^{1})}\bigr)\\
				&\lesssim \sum_{|j-j'|\leq 4}\Bigl(\sum_{\ell\leq j+3}\|\dot{\Delta}_\ell v\|_{L^\infty([-M_3,0]; L^{4}(B(0,10M_1)))}\Bigr)\|\dot{\Delta}_{j'}v\|_{L^{\infty}([-M_3,0] \times B(0,10M_1))}\\
				&\quad\quad+M_4^{-50c_p}\|v\|^2_{L^\infty([-M_3,0];  L^{2})}\\
				&\lesssim 2^{\frac{5}{4}j}M_1^{-1}M^{2c_p}+M_4^{-50c_p}M_3^{\frac12} M^{c_p}
\quad\mbox{[by using \eqref{est:Deltaj-v-LinftyLinfty} and \eqref{LinftyL2vnew}]}\\
				&\lesssim 2^{\frac{5}{4}j}M_1^{-1}M^{2c_p}+M_4^{-49c_p},
			\end{align*}
where we used \eqref{LinftyL2vnew} so that
$$\sum_{\ell\leq j+3}\|\dot{\Delta}_\ell v\|_{L^\infty([-M_3,0]; L^{4}(B(0,10M_1)))}\lesssim 2^{\frac{j}4} M^{c_p}.$$
Similarly for the remaining term in \eqref{decompose-v}, we obtain
			\begin{align*}
				&\bigl\|\dot{\Delta}_j\bigl(\sum_{\substack{\ell\geq j-4}}\dot{\Delta}_\ell v\otimes\tilde{\Delta}_{\ell}v\bigr)
\bigr\|_{L^\infty([-M_3,0];  L^{4}(B(0,9M_1)))}\\
				&\lesssim 2^{\frac{9j}{4}}\sum_{\substack{\ell\geq j-4\\ 2^\ell\leq M_2}}\|\dot{\Delta}_\ell v\|_{L^\infty([-M_3,0]; L^{1}(B(0,10M_1)))}\|\tilde{\Delta}_\ell v\|_{L^\infty ([-M_3,0] \times B(0,10M_1))}\ [\mbox{using (\ref{LinftyL2vnew}, \ref{est:Deltaj-v-LinftyLinfty})}]\\
				&+2^{\frac{9}{4}j}\sum_{\substack{\ell\geq j-4\\ 2^\ell\geq M_2}}\|\dot{\Delta}_\ell v\|_{L^\infty([-M_3,0]; L^{2}(B(0,10M_1)))}\|\tilde{\Delta}_\ell v\|_{L^\infty([-M_3,0]; L^{2}(B(0,10M_1)))}\ [\mbox{using (\ref{LinftyL2vnew})}]\\
				&+M_4^{-50c_p}2^{\frac{9}{4}j}\sum_{\substack{\ell\geq j-4}}\|\dot\Delta_\ell v\|_{L^\infty([-M_3,0];  L^{2})}
\|\tilde\Delta_{\ell}v\|_{L^\infty ([-M_3,0];L^{2})}\ [\mbox{using (\ref{LinftyL2vnew})}]\\
				&\lesssim 2^{\frac{9}{4}j}\Bigl(\sum_{\substack{\ell \geq j-4\\ 2^\ell\leq M_2}}2^{-2\ell} M^{c_p} 2^\ell M_1^{-1} +  \sum_{\substack{\ell\geq j-4\\ 2^\ell\geq M_2}} (2^{\frac{3}{2}\ell} 2^{-2\ell} M^{c_p})^2 + M_4^{-50c_p}\sum_{\substack{l\geq j-4}}(2^{\frac{3}{2}l} 2^{-2l} M^{c_p})^2\Bigr)\\
				&\lesssim 2^{\frac{5}{4}j}M_1^{-1}M^{2c_p}+ 2^{\frac{9}{4}j} M_2^{-1}M^{2c_p} +  M_4^{-50c_p}2^{-j}M^{2c_p}\lesssim 2^{\frac{5}{4}j}M_1^{-1}M^{2c_p}+M_4^{-49c_p}.
\end{align*}
		Then by virtue of \eqref{decompose-v}, we get, by summarizing the above estimates, that
			\begin{equation*}
				\|\dot{\Delta}_j(v\otimes v)\|_{L^\infty ([-M_3,0]; L^{4}(B(0,9M_1)))}\lesssim 2^{\frac{5}{4}j}M_1^{-1}M^{2c_p}.
			\end{equation*}
	By inserting the above estimate  and \eqref{LinftyL2vnew} into 	 \eqref{est:Duhamel}, we obtain
			\begin{equation}\label{est:L4vN}
				\|\dot{\Delta}_jv_{N}(-M_3,t)\|_{L^\infty([-M_3/2,0];  L^{4}(B(0,8M_1)))}\lesssim 2^{\frac{j}{4}}M_1^{-1}M^{2c_p}.
			\end{equation}
			
			\noindent \textbf{The estimate of $\dot{\Delta}_j v_L(-M_3,t)$.}	In view of  \eqref{def:Vi}, we get
			by applying Lemma \ref{localY}, that for any $ q\in[1,\infty]$
			\begin{equation}\label{est:DuhamelL}
				\begin{aligned}
					\|\dot{\Delta}_j v_L(-M_3,t)\|_{L_x^q(B(0,8M_1))}\lesssim
\int_{-M_3}^t2^je^{-c_0(t-s)2^{2j}}&\left(\|\dot\Delta_j(v\otimes V_1)(s)\|_{L^q(B(0,9M_1))}\right.\\
					&\left.\ +M_4^{-50c_p}\|v(s)\|_{L^2}\|V_1(s)\|_{L^p} \right)ds.
				\end{aligned}
			\end{equation}
		 We shall use Bony's  decomposition \eqref{form:decomposition} to handle $\dot\Delta_j(v\otimes V_1).$
In view of \eqref{est:Katonew},  the high frequency part of $u_{kL}$ has exponential decay with time, from which  and
 Lemma \ref{localY}, we deduce that
			\begin{align*}
				&\sum_{|j'-j|\leq 4}\| \tilde{\Delta}_{j}(\dot S_{j'-1}v\otimes\dot\Delta_{j'}V_1)\|_{L^\infty([-M_3,0];  L^{4}(B(0,9M_1)))} \\
				&\lesssim\sum_{i=1}^m \sum_{|j'-j|\leq 4}\|\tilde{\Delta}_j(\dot S_{j'-1}v\otimes\dot{\Delta}_{j'}u_{iL})\|_{L^\infty([-M_3,0];  L^{4}(B(0,9M_1)))}\lesssim 2^{\frac{5}{4}j}M_1^{-1}M^{2c_p},
\end{align*}
and
		\begin{align*}		&\|\tilde{\Delta}_j(\sum_{\substack{\ell\geq j-4}}\dot{\Delta}_\ell v \otimes\tilde{\Delta}_{\ell}V_1)\|_{L^\infty([-M_3,0];  L^{4}(B(0,9M_1)))} \\
				&\lesssim \sum_{i=1}^m \|\tilde{\Delta}_j(\sum_{\substack{\ell\geq j-4}}\dot{\Delta}_\ell v\otimes\tilde{\Delta}_{\ell}u_{iL})\|_{L_t^\infty L_x^{4}([-M_3,0] \times B(0,9M_1))}\lesssim 2^{\frac{5}{4}j}M_1^{-1}M^{2c_p}.
			\end{align*}
			In order to deal with the term $\sum_{|j'-j|\leq 4}\tilde{\Delta}_j(\dot S_{j'-1}u_{iL}\otimes\dot{\Delta}_{j'}v)$, we split the frequencies into two ranges $2^\ell \leq M_2^{-1}$ and $2^\ell \geq M_2^{-1}$ . Precisely, for any $t\in[-M_3,0]$, we have
			\begin{align*}
				&\sum_{|j'-j|\leq 4}\|\tilde{\Delta}_j(\dot S_{j'-1}V_1 \otimes\dot{\Delta}_{j'}v)(t)\|_{ L^{4}(B(0,9M_1))}\\
				&\lesssim \sum_{i=1}^m \sum_{|j'-j|\leq 4} \Bigl(\sum_{\substack{\ell\leq j+4\\ 2^\ell\leq M_2^{-1}}}\|\dot{\Delta}_\ell u_{iL}(t)\|_{ L^{\infty}( B(0,10M_1))}\\
&\qquad\qquad\qquad\quad+\sum_{\substack{\ell\leq j+4\\ 2^\ell\geq M_2^{-1}}}\|\dot{\Delta}_\ell u_{iL}(t)\|_{ L^{\infty}( B(0,10M_1))}\Bigr)\|\dot{\Delta}_{j'}v\|_{ L^{4}(B(0,10M_1))}\ [\mbox{using (\ref{est:Katonew},\ref{LinftyL2vnew})}]\\
				&+M_4^{-50c_p} \sum_{i=1}^m\sum_{|j'-j|\leq 4} \sum_{\ell \leq j'-2}\|\dot \Delta_{\ell} u_{iL}(t)\|_{L^p}\|\dot\Delta_{j'}v(t)\|_{L^p}\ [\mbox{using  \eqref{est:v-Besov}}]\\
				&\lesssim \sum_{i=1}^m \Bigl(\sum_{\substack{l\leq j+4\\ 2^\ell\leq M_2^{-1}}} 2^{\ell}M^{c_p} 2^{\frac{9}{4}\ell} 2^{-2\ell} M^{c_p} +  \sum_{\substack{\ell\leq j+4\\ 2^\ell\geq M_2^{-1}}} e^{-c_{i-1}(t+M_3)2^{2\ell}}M^{c_p} 2^{\frac{9}{4}\ell} 2^{-2\ell} M^{c_p} \\
				&\qquad \quad + M_4^{-50c_p}  \sum_{\ell \leq j'-2} 2^{\left(1-\frac{3}{p}\right)\ell}M^{c_p} 2^{\left(1-\frac{3}{p}\right)j}M^{c_p}\Bigr) \ [\mbox{using}\ \ M_2^{-\frac{1}{2}} \leq 2^j \leq M_2^{\frac{1}{2}}] \\
				&\lesssim M_2^{-1}2^{\frac{5}{4}j}M^{2c_p}+e^{-c_{m-1}2^{-1}M_3M_2^{-2}}2^{\frac{1}{4}j}M_2M^{2c_p}+M_4^{-45c_p},\quad \forall t\in[-M_3/2,0].
			\end{align*}

	Hence by virtue of \eqref{form:decomposition} and by substituting the above estimates into  \eqref{est:DuhamelL}, we obtain
			\begin{equation}\label{est:L4vL}
				\|\dot{\Delta}_jv_{L}(-M_3,t)\|_{L^\infty([-M_3/2,0];  L^{4}(B(0,8M_1)))}\lesssim 2^{\frac{j}{4}}M_1^{-1}M^{2c_p}.
			\end{equation}
			
			\noindent \textbf{The estimate of $\dot{\Delta}_j v_F(-M_3,t)$}.  Notice from \eqref{est:Katonew} that
			$\dot\Delta_ju_{iL}$ has exponential decay in time, we find
			\begin{equation}\label{est:L4vF}
				\|\dot{\Delta}_jv_F(-M_3,t)\|_{L^\infty([-M_3/2,0];  L^{4}(B(0,8M_1)))}\lesssim M_4^{-50c_p}.
			\end{equation}
		By summarizing \eqref{est:LinftyL4v0}, \eqref{est:L4vN}, \eqref{est:L4vL} and \eqref{est:L4vF}, we obtain, for $M_2^{-\frac{1}{2}} \leq 2^j \leq M_2^{\frac{1}{2}}$, that
			\begin{equation}\label{est:LinftyL4-v}
				\| \dot{\Delta}_j v \|_{L^\infty([-M_3/2,0];  L^{4}(B(0,8M_1)))} \lesssim 2^{\frac{j}{4}} M_1^{-1} M^{2c_p},
			\end{equation}
			which, together with \eqref{est:LinftyL4-V1}, ensures that for this range of $j$,
			\begin{equation}\label{est:small1}
				\| \dot{\Delta}_j u \|_{L^\infty([-M_3/2,0];  L^{4}(B(0,8M_1)))} \lesssim 2^{\frac{j}{4}} M_1^{-1} M^{2c_p}.
			\end{equation}
			
			\noindent {\sc Step 3.}  The case for  $M_2^{-\frac{1}{3}} \leq 2^j \leq M_2^{\frac{1}{3}}$.

			Once again, we shall use the decomposition $u=v+V_1.$  To deal with the estimate
of $v,$ similar to \eqref{S4eq8}, we write \begin{equation}	\label{S4eq9}
				v(t) =	e^{\left(t+{M_3}/2\right)\Delta} v\bigl(-{M_3}/2\bigr) +v_N\bigl(-{M_3}/2,t\bigr)+v_L\bigl(-{M_3}/2,t\bigr)+v_F\bigl(-{M_3}/2,t\bigr),
			\end{equation}
where $ v_N(-M_3/2,t),$ $v_L(-M_3/2,t),$ $v_F(-M_3/2,t)$ are given by \eqref{form:vF}.
 By using a similar argument as in \textbf{Step 2}, we can estimate the terms $\dot{\Delta}_j V_1$, $\dot{\Delta}_j v_0$ and $\dot{\Delta}_j v_F$ as	 follows
			\begin{equation} \label{est:L2v0}
				\begin{aligned}
					&\|\dot{\Delta}_j V_1\|_{L^\infty([-M_3/3,0];  L^{2}(B(0,3M_1)))}
					+ \|\dot{\Delta}_j v_F(-M_3/2,t)\|_{L^\infty([-M_3/3,0];  L^{2}(B(0,3M_1)))}\\
&+ \bigl\|\dot{\Delta}_j e^{\left(t+{M_3}/2\right)\Delta} v\bigl(-{M_3}/2\bigr) \bigr\|_{L^\infty([-M_3/3,0];  L^{2}(B(0,3M_1)))} \lesssim M_4^{-50c_p}.
				\end{aligned}
			\end{equation}

			To estimate	 $\dot{\Delta}_j v_N\bigl(-{M_3}/2,t\bigr)$, we use the decomposition \eqref{decompose-v} to deal with the estimate of $\dot\Delta_j(v\otimes v).$
		For the para-product part, we get, by using Lemma \ref{localY}, that
			\begin{align*}
				&\sum_{|j-j'|\leq 4}\|\dot{\Delta}_j(S_{j'-1}v\otimes\dot{\Delta}_{j'}v)\|_{L^\infty([-M_3/2,0];  L^{2}(B(0,5M_1)))}\\
				&\lesssim \sum_{|j-j'|\leq 4}\sum_{\ell\leq j+4}\|\dot{\Delta}_\ell v\|_{L^\infty([-M_3/2,0];  L^{4}(B(0,6M_1)))}
\|\dot{\Delta}_{j'} v\|_{L^\infty ([-M_3/2,0]; L^4(B(0,6M_1)))}\ [(\ref{LinftyL2vnew},\ref{est:LinftyL4-v})]\\
				&\quad+M_4^{-50c_p}\sum_{|j-j'|\leq 4}\sum_{\ell \leq j+4}\| \dot{\Delta}_\ell  v\|_{{L^\infty([-M_3/2,0] \times \mathbb{R}^3)}}\|\dot\Delta_{j'} v\|_{L^\infty ([-M_3/2,0]; L^p)}\ [\mbox{using}\ (\ref{est:Deltaj-v-LinftyLinfty},\ref{est:v-Besov})]\\
				&\lesssim \sum_{|j-j'|\leq 4}\sum_{\ell\leq j+4}2^{\frac{9}{4}\ell} 2^{-2\ell} M^{c_p} 2^{\frac{j'}{4}}M_1^{-1}M^{2c_p} + M_4^{-50c_p} \sum_{|j-j'|\leq 4}\sum_{\ell \leq j+4} 2^\ell M_1^{-1} 2^{\left(1-\frac{3}{p}\right)j'}M^{c_p} \\
				&\lesssim 2^{\frac{j}{2}}M_1^{-1}M^{3c_p}.
			\end{align*}
Along the same line, to handle the remaining term in \eqref{decompose-v},
		 we split further the frequencies into two ranges $2^\ell \leq M_2^{\frac{1}{2}}$ and $2^\ell \geq M_2^{\frac{1}{2}}$ as follows
			\begin{align*}
				&\|\dot{\Delta}_j(\sum_{\substack{\ell\geq j-4}}\dot{\Delta}_\ell v\otimes\tilde{\Delta}_{\ell}v)
\|_{L^\infty([-M_3/2,0]; L^{2}(B(0,5M_1)))}\\
				&\lesssim 2^{\frac{3}{4}j}\sum_{\substack{\ell\geq j-4\\ 2^\ell\leq M_2^{\frac{1}{2}}}}
\|\dot{\Delta}_\ell v\|_{L^\infty([-M_3/2,0]; L^{2}(B(0,6M_1)))}\|\tilde{\Delta}_{\ell}v\|_{L^\infty ([-M_3/2,0]; L^{4}(B(0,6M_1)))}\\
				&\quad +2^{\frac{3}{2}j}\sum_{\substack{\ell\geq j-4\\ 2^\ell\geq M_2^{\frac{1}{2}}}}
\|\dot{\Delta}_\ell v\|_{L^\infty([-M_3/2,0];  L^{2}(B(0,6M_1)))}\|\tilde{\Delta}_{\ell}v\|_{L^\infty([-M_3/2,0]; L^{2}(B(0,6M_1)))}\\
				&\quad+M_4^{-50c_p}\sum_{\substack{\ell\geq j-4}}\|\dot{\Delta}_\ell v\|_{{L^\infty([-M_3/2,0];  L^p)}} \| \tilde\Delta_{\ell}v\|_{{L^\infty ([-M_3/2,0];L^{\frac{p}{p-1}})}}\\
				&\lesssim 2^{\frac{3}{4}j}\sum_{\substack{\ell\geq j-4\\ 2^\ell\leq M_2^{\frac{1}{2}}}} 2^{\frac{3}{2}\ell} 2^{-2l} M^{c_p} 2^{\frac{l}{4}}M_1^{-1}M^{2c_p} + 2^{\frac{3}{2}j} \sum_{\substack{\ell\geq j-4\\ 2^\ell\geq M_2^{\frac{1}{2}}}}\bigl(2^{\frac{3}{2}\ell}2^{-2l}M^{c_p}\bigr)^2 \\
				&\quad + M_4^{-50c_p} \sum_{\ell \geq j-4} 2^{\left(1-\frac{3}{p}\right)\ell}M^{c_p} 2^{\frac{3}{p}\ell}2^{-2\ell}M^{c_p} \\
				&\lesssim  2^{\frac{j}{2}}M_1^{-1}M^{3c_p} + 2^{j}M_2^{-\frac{1}{4}}M^{2c_p} + M_4^{-50c_p}2^{-j}M^{2c_p} \lesssim 2^{\frac{j}{2}}M_1^{-1}M^{3c_p},
			\end{align*}
			for $2^j \in [M_2^{-\frac{1}{3}},M_2^{\frac{1}{3}}]$.\\
            We thus deduce from \eqref{decompose-v} that
			\begin{equation*}
				\|\dot{\Delta}_j(v\otimes v)\|_{L^\infty([-M_3/2,0]; L^{2}(B(0,6M_1)))}\lesssim 2^{\frac{j}{2}}M_1^{-1}M^{3c_p},
			\end{equation*}
			from which and similar version of \eqref{est:Duhamel}, we infer that
			\begin{equation}\label{est:L2vN}
				\|\dot{\Delta}_jv_{N}\bigl(-{M_3}/2,t\bigr)\|_{L^\infty([-M_3/3,0];  L^{2}(B(0,5M_1)))}\lesssim 2^{-\frac{j}{2}}M_1^{-1}M^{3c_p}.
			\end{equation}
			
			To estimate $v_L\bigl(-{M_3}/2,t\bigr)$, we shall first use Bony's decomposition \eqref{form:decomposition} for $v\otimes V_1.$ By virtue of \eqref{est:Katonew} and \eqref{LinftyL2vnew}, we can show that
			\begin{align*}
				&\sum_{|j-j'|\leq 4}\|\dot{\Delta}_j(\dot S_{j'-1}v\otimes\dot{\Delta}_{j'}V_1)\|_{L^\infty([-M_3/2,0];  L^{2}(B(0,6M_1)))}\lesssim 2^{\frac{j}{2}}M_1^{-1}M^{3c_p},\\
				&\|\dot{\Delta}_j(\sum_{\substack{\ell\geq j-4}}\dot{\Delta}_\ell v\otimes\tilde{\Delta}_{\ell}V_1)\|_{L^\infty([-M_3/2,0];  L^{2}(B(0,6M_1)))}\lesssim 2^{\frac{j}{2}}M_1^{-1}M^{3c_p}.
			\end{align*}
			To deal with the para-product term $T_{V_1}v$, we get, by using a similar argument as in \textbf{Step 2}, that
			\begin{align*}
				&\sum_{|j-j'|\leq 4}\|\dot{\Delta}_j(\dot S_{j'-1}V_1 \otimes\dot{\Delta}_{j'}v)\|_{ L^\infty([-M_3/2,0]; L^{2}(B(0,6M_1)))}\\
				&\lesssim \sum_{\substack{\ell\leq j+4\\ 2^\ell\leq M_2^{-\frac{1}{3}}}}\|\dot{\Delta}_\ell V_1
\|_{L^\infty([-M_3/2,0] \times B(0,7M_1))}\|\dot{\Delta}_j v\|_{ L^\infty([-M_3/2,0];  L^{2}(B(0,7M_1)))}\\
				&\quad +\sum_{\substack{\ell\leq j+4\\ 2^\ell\geq M_2^{-\frac{1}{3}}}}\|\dot{\Delta}_\ell V_1\|_{ L^{\infty}([-M_3/2,0] \times B(0,7M_1))}\|\dot{\Delta}_jv\|_{L^\infty([-M_3/2,0];  L^{2}(B(0,7M_1)))}\\
				&\quad +M_4^{-50c_p}\|V_1\|_{{L^\infty([-M_3/2,0];  L^p)}}\|\dot\Delta_jv\|_{{L^\infty([-M_3/2,0];  L^p)}}\\
				&\lesssim M_2^{-\frac{1}{3}}2^{-\frac{j}{2}}M^{2c_p}+e^{-c_{m-1}\frac{M_3}{6}M_2^{-\frac{1}{3}}}2^{\frac{j}{2}}M^{2c_p}+M_4^{-45c_p} \lesssim 2^{-\frac{j}{2}}M_2^{-\frac{1}{3}}M^{2c_p}.
			\end{align*}
By summing up the above estimates and using a similar version of \eqref{est:DuhamelL}, we deduce that
			\begin{equation}\label{est:L2vL}
				\|\dot{\Delta}_jv_{L}\bigl(-{M_3}/2,t\bigr)\|_{L^\infty ([-M_3/3,0];L^{2}(B(0,5M_1)))}\lesssim 2^{-\frac{j}{2}}M_1^{-1}M^{2c_p}.
			\end{equation}

	By summarizing the estimates \eqref{est:L2v0}, \eqref{est:L2vN} and \eqref{est:L2vL},
we conclude  that for this range of $j$,
			\begin{equation}\label{est:small2}
				\| \dot{\Delta}_j u \|_{L^\infty ([-M_3/3,0]; L^{2}(B(0,5M_1)))} \lesssim 2^{-\frac{j}{2}} M_1^{-1} M^{3c_p}.
			\end{equation}
			
			\noindent {\sc Step 4.} The case for $M_2^{-\frac{1}{4}} \leq 2^j \leq M_2^{\frac{1}{4}}$.

			Again we shall use the decomposition $u=v+V_1.$  To deal with the estimate
of $v,$ similar to \eqref{S4eq9}, we write \begin{equation}	\label{S4eq10}
				v(t) =	e^{\left(t+{M_3}/3\right)\Delta} v\bigl(-{M_3}/3\bigr)+v_N\bigl(-{M_3}/3,t\bigr)+v_L\bigl(-{M_3}/3,t\bigr)+v_F\bigl(-{M_3}/3,t\bigr),
			\end{equation}
where $ v_N(-M_3/3,t),$ $v_L(-M_3/3,t)$ and $v_F(-M_3/3,t)$ are given by \eqref{form:vF}.

 We first get, by using a similar argument as in \textbf{Step 2}, that
			\begin{equation} \label{est:L3/2v0}
				\begin{aligned}
					&\|\dot{\Delta}_j V_1\|_{L^\infty([-M_3/4,0];  L^{\frac{4}{3}}(B(0,M_1)))}+\|\dot{\Delta}_j v_F(-{M_3}/3,t)\|_{L^\infty ([-M_3/4,0];L^{\frac{4}{3}}(B(0,M_1)))}\\
&+\bigl\|\dot{\Delta}_je^{\left(t+{M_3}/3\right)\Delta} v\bigl(-{M_3}/3\bigr)\bigr\|_{L^\infty([-M_3/4,0];  L^{\frac{4}{3}}(B(0,M_1)))}   \lesssim M_4^{-50}.
				\end{aligned}
			\end{equation}

To estimate	 $\dot{\Delta}_j v_N\bigl(-{M_3}/3,t\bigr)$, we use the decomposition \eqref{decompose-v} to deal with the estimate of $\dot\Delta_j(v\otimes v).$
		For the para-product part, we get, by using Lemma \ref{localY}, that
			\begin{align*}
				&\sum_{|j-j'|\leq 4}\|\dot{\Delta}_j(\dot S_{j'-1}v\otimes\dot{\Delta}_{j'}v)\|_{L^\infty ([-M_3/3,0]; L^{\frac{4}{3}}(B(0,2M_1)))}\\
				&\lesssim \sum_{|j-j'|\leq 4}\sum_{\ell\leq j+4}\|\dot{\Delta}_\ell v\|_{L^\infty([-M_3/3,0];  L^{4}(B(0,3M_1)))}
\|\dot{\Delta}_{j'} v\|_{L^\infty([-M_3/3,0];  L^2(B(0,3M_1)))}\\
				&\quad+M_4^{-50c_p}\sum_{|j-j'|\leq 4} \sum_{\ell \leq j+4} \|\dot{\Delta}_\ell v\|_{{L^\infty([-M_3/3,0] \times \mathbb{R}^3)}}\|\dot\Delta_{j'} v\|_{{L^\infty([-M_3/3,0] \times \mathbb{R}^3)}} \\
				&\lesssim \sum_{|j-j'|\leq 4}\sum_{\ell\leq j+4} 2^{\frac{\ell}{4}}M^{c_p} 2^{-\frac{j'}{2}}M_1^{-1}M^{3c_p}
+ M_4^{-50c_p} \sum_{|j-j'|\leq 4} \sum_{\ell \leq j+4} 2^{\frac{3}{p}\ell}2^{\left(1-\frac{3}{p}\right)\ell}M^{c_p} M_1^{-1}2^{j'} \\
				&\lesssim 2^{-\frac{j}{4}}M_1^{-1}M^{4c_p}.
			\end{align*}
Similarly for the remaining term in	\eqref{decompose-v}, we have
			\begin{align*}
				&\bigl\|\dot{\Delta}_j\bigl(\sum_{\substack{\ell\geq j-4}}\dot{\Delta}_\ell v\otimes\tilde{\Delta}_{\ell}v\bigr)
\bigr\|_{L^\infty([-M_3/3,0];  L^{\frac{4}{3}}(B(0,2M_1)))}\\
				&\lesssim 2^{\frac{3}{4}j}\sum_{\ell\geq j-4}\|\dot{\Delta}_\ell v\|_{L^\infty ([-M_3/3,0]; L^{2}(B(0,3M_1)))}
\|\tilde{\Delta}_\ell v\|_{L^\infty ([-M_3/3,0]; L^{2}(B(0,3M_1)))}\\
				&\quad+M_4^{-50c_p}\sum_{\substack{\ell\geq j-4}}\|\dot\Delta_\ell v\otimes\tilde\Delta_{\ell}
v\|_{{L^\infty([-M_3/3,0];  L^1)}}\\
				&\lesssim 2^{\frac{3}{4}j} \sum_{\substack{\ell \geq j-4 \\ 2^\ell \leq M_2^{\frac{1}{3}}}} \bigl(2^{-\frac{\ell}{2}}M_1^{-1}M^{3c_p}\bigr)^2 +  2^j \sum_{\substack{\ell \geq j-4 \\ 2^\ell \geq M_2^{\frac{1}{3}}}} \bigl(2^{\frac{3}{2}\ell}2^{-2\ell}M^{c_p}\bigr)^2 + M_4^{-45c_p} \lesssim 2^{-\frac{j}{4}} M_1^{-2}M^{6c_p}.
			\end{align*}
	In view of \eqref{decompose-v}, we get, by summing up	 the preceding two estimates, that
			\begin{equation*}
				\|\dot{\Delta}_j(v\otimes v)\|_{L^\infty ([-M_3/3,0]; L^{\frac{4}{3}}(B(0,2M_1)))}\lesssim 2^{-\frac{j}{4}}M_1^{-1}M^{6c_p},
			\end{equation*}
		from 	which and a  similar version of \eqref{est:Duhamel}, we deduce that
			\begin{equation}\label{est:L3/2vN}
				\|\dot{\Delta}_jv_{N}(-{M_3}/3,t)\|_{L^\infty([-M_3/4,0];  L^{\frac{4}{3}}(B(0,M_1)))}\lesssim 2^{-\frac{5}{4}j}M_1^{-1}M^{6c_p}.
			\end{equation}

			To deal with the estimate $v_L\bigl(-{M_3}/3,t\bigr)$, we shall use Bony's decomposition \eqref{form:decomposition} for $v\otimes V_1.$
Indeed we get, by 	\eqref{est:Katonew} and \eqref{LinftyL2vnew}, that
			\begin{align*}
				&\sum_{|j-j'|\leq 4}\|\dot{\Delta}_j(\dot S_{j'-1}v\otimes\dot{\Delta}_{j'}V_1)\|_{L^\infty([-M_3/3,0];  L^{\frac{4}{3}}(B(0,2M_1)))}\lesssim 2^{-\frac{j}{4}}M_1^{-1}M^{4c_p},\\
				&\|\tilde{\Delta}_j(\sum_{\substack{\ell\geq j-4}}\dot{\Delta}_\ell v\otimes\tilde{\Delta}_{\ell}
V_1)\|_{L^\infty ([-M_3/3,0]; L^{\frac{4}{3}}(B(0,2M_1)))}\lesssim 2^{-\frac{j}{4}}M_1^{-1}M^{4c_p}.
			\end{align*}
			For $\dot{\Delta}_j(T_{V_1}v)$, we deduce that for $2^j \in [M_2^{-\frac{1}{4}},M_2^{\frac{1}{4}}]$ and for any $t\in[-{M_3}/{3},0]$,
			\begin{align*}
				&\sum_{|j-j'|\leq 4}\|\dot{\Delta}_j(\dot S_{j'-1}V_1\otimes\dot{\Delta}_{j'}v)(t)\|_{ L^{\frac{4}{3}}(B(0,2 M_1))}\\
				&\lesssim \sum_{\substack{\ell\leq j+4\\ 2^\ell\leq M_2^{-\frac{1}{3}}}}\|\dot{\Delta}_\ell V_1(t)\|_{ L^{\infty}( B(0,3M_1))}\|\dot{\Delta}_jv\|_{ L^{\frac{4}{3}}(B(0,3M_1))}\\
				&\quad +\sum_{\substack{\ell\leq j+4\\ 2^\ell\geq M_2^{-\frac{1}{3}}}}\|\dot{\Delta}_\ell V_1(t)\|_{ L^{\infty}( B(0,3M_1))}\|\dot{\Delta}_jv(t)\|_{ L^{\frac{4}{3}}( B(0,3M_1))}\\
				&\quad+M_4^{-50c_p}\|\dot\Delta_jv\|_{{L^\infty([-M_3/3,0] \times \mathbb{R}^3)}}
\|\dot S_jV_1\|_{{L^\infty([-M_3/3,0];  L^p)}}\\
				&\lesssim M_2^{-\frac{1}{3}}2^{-\frac{5}{4}j}M^{2c_p}+ 2^j e^{-c_{m-1}\left(t+\frac{M_3}{3}\right)M_2^{-\frac{2}{3}}}M^{c_p} 2^{-\frac{5}{4}j}M^{c_p} + M_4^{-45c_p} \lesssim 2^{-\frac{j}{4}}M_1^{-1}M^{2c_p}.
			\end{align*}
By summarizing the above estimates and using a similar version of \eqref{est:DuhamelL}, we obtain
			\begin{equation} \label{est:L3/2vL}
				\|\dot{\Delta}_jv_{L}(-{M_3}/3,t)\|_{L^\infty([-M_3/4,0];  L^{\frac{4}{3}}(B(0,M_1)))}\lesssim 2^{-\frac{5}{4}j}M_1^{-1}M^{6c_p}.
			\end{equation}

	In view of \eqref{S4eq10}, by summing up the estimates \eqref{est:L3/2v0}, \eqref{est:L3/2vN} and \eqref{est:L3/2vL}, we achieve
			\begin{equation*}\label{est:LinftyL3/2-v}
				\| \dot{\Delta}_j v \|_{L^\infty ([-M_3/4,0]; L^{\frac{4}{3}}(B(0,M_1)))} \lesssim 2^{-\frac{5}{4}j} M_1^{-1} M^{4c_p},
			\end{equation*}
which, together with the estimate of $V_1$ in \eqref{est:L3/2v0}, ensures that for $2^j \in [M_2^{-\frac{1}{4}},M_2^{\frac{1}{4}}]$
			\begin{equation}\label{est:small3}
				\| \dot{\Delta}_j u \|_{L^\infty([-M_3/4,0];  L^{\frac{4}{3}}(B(0,M_1)))} \lesssim 2^{-\frac{5}{4}j} M_1^{-1} M^{4c_p}.
			\end{equation}

			\noindent {\sc Step 5.} We are going to show the contradiction with \eqref{assump:contradition}. Indeed
 we first get, by using Lemma \ref{lem:Bernstein}, \eqref{assumption} and Lemma \ref{localY}, that
			\begin{align*}
				|\dot{\Delta}_0 u(0,0)| &\leq |e^{\frac{M_3}{5}\Delta} \dot{\Delta}_0 u(-M_3/5,0)|
 + \int_{-M_3/5}^0 |e^{s\Delta} \tilde{\Delta}_0 \bP \nabla \cdot   \dot\Delta_0 (u \otimes u)(s,0)|ds  \\
				&\lesssim e^{-c_0\frac{M_3}{5}}M^{c_p} + \int_{-M_3/5}^0 e^{c_0 s} \bigl(\| \tilde \Delta_0 (u \otimes u)(s)|\|_{L^1(B(0,M_1/2))} + M_1^{-50} \bigr)ds,
			\end{align*}
	from which, \eqref{assump:contradition} and the pigeonhole principle, we deduce that there exists $s\in [-M_3/5,0]$ such that
			\begin{equation} \label{est:contra}
				M_1^{-1} \lesssim \| \dot \Delta_0 (u \otimes u)(s)|\|_{L^1(B(0,M_1/2))}.
			\end{equation}
			By using the Bony decomposition \eqref{form:decomposition}, we have
			\begin{equation} \label{est:split-dyadic0} \begin{aligned}
					\tilde \Delta_0 (u \otimes u)(s) = 2 \sum_{|j|\leq 4} \dot \Delta_{0} (\dot{S}_{j-1} u(s) \otimes \dot{\Delta}_j u(s)) + \sum_{\substack{\ell\geq -4}} \dot \Delta_{0} (\dot{\Delta}_{\ell} u(s) \otimes \tilde{\Delta}_{\ell} u(s)).
				\end{aligned} \end{equation}
Notice that  if $2^{j'} \leq M_2^{-\frac{1}{4}},$  it follows from Lemma \ref{lem:Bernstein} and \eqref{assumption} that
 \begin{align*}
 \| \dot{\Delta}_{j'} u(s) \|_{L^4(B(0,M_1))}\leq &M_1^{\frac34}\| \dot{\Delta}_{j'} u(s) \|_{L^\infty(B(0,M_1))}\lesssim  M_1^{\frac34}2^{\frac3p j'}\| \dot{\Delta}_{j'} u(s) \|_{L^p(B(0,M_1))}\\
 \lesssim & M_1^{\frac34}2^{j'}\|u(s)\|_{\dot B^{-1+\frac3p}_{p,\infty}}\lesssim  M_1^{\frac34}2^{j'} M^{c_p},
 \end{align*}
from which and \eqref{est:small1}, we infer for $|j|<4,$
				\begin{align*}
					\| \dot S_{j-1} u(s)\|_{L^4(B(0,M_1))}
					&\leq \Bigl(\sum_{j' < 2,\, 2^{j'} \geq M_2^{-\frac{1}{4}}}  + \sum_{j' < 2,\, 2^{j'} \leq M_2^{-\frac{1}{4}}} \Bigr)\| \dot{\Delta}_{j'} u(s) \|_{L^4(B(0,M_1))}\\
					&\lesssim M_1^{-1}M^{2c_p}+ \sum_{j' < -4,\, 2^{j'} \leq M_2^{-\frac{1}{4}}} 2^{j'}M^{c_p} M_1^{\frac{3}{4}} \lesssim M_1^{-1}M^{3c_p}.
				\end{align*}
This,  together with \eqref{assumption} and \eqref{est:small3}, implies that, for all $s\in[-{M_3}/{5},0]$,
				\begin{equation} \label{est:dydic0-1} \begin{aligned}
						&\| \sum_{|j|\leq 4}\dot \Delta_{0} (\dot{S}_{j-1} u(s) \otimes \dot{\Delta}_{j} u(s))  \|_{L^1(B(0,M_1/2))} \\
						&\lesssim \sum_{|j|\leq 4} \| \dot{S}_{j-1} u(s))    \|_{L^4(B(0,M_1))} \| \dot{\Delta}_{j} u(s)  \|_{L^\frac{4}{3}(B(0,M_1))} 	 \\
						&\quad+M_4^{-50c_p}\sum_{|j|\leq 4}\|\dot S_{j-1}u\|_{{L^\infty([-M_3/5,0] \times \mathbb{R}^3)}}\|\dot\Delta_ju\|_{{L^\infty([-M_3/5,0] \times \mathbb{R}^3)}} \lesssim M_1^{-2}M^{6c_p}.	
					\end{aligned} \end{equation}
					
		For the remaining term in \eqref{est:split-dyadic0}, since $u=v+V_1$, we write
                        \begin{align*}
                            &\bigl\| \sum_{\substack{\ell\geq -4}} \dot \Delta_{0} (\dot{\Delta}_{\ell} u(s) \otimes \tilde{\Delta}_{\ell} u(s)) \bigr\|_{L^1(B(0,M_1/2))}\\
                            &\lesssim \sum_{\substack{\ell\geq -4}}\Bigl(\bigl\|  \dot \Delta_{0} (\dot{\Delta}_{\ell} v(s) \otimes \tilde{\Delta}_{\ell} v(s)) \bigr\|_{L^1(B(0,M_1/2))}+\bigl\| \dot \Delta_{0} (\dot{\Delta}_{\ell} V_1(s) \otimes \tilde{\Delta}_{\ell} V_1(s)) \bigr\|_{L^1(B(0,M_1/2))}\\
                            &\quad +\bigl\|  \dot \Delta_{0} (\dot{\Delta}_{\ell} V_1(s) \otimes \tilde{\Delta}_{\ell} v(s) )\bigr\|_{L^1(B(0,M_1/2))}+\bigl\|\dot \Delta_{0} (\dot{\Delta}_{\ell} v(s) \otimes \tilde{\Delta}_{\ell} V_1(s)) \bigr\|_{L^1(B(0,M_1/2))}\Bigr).
                        \end{align*}
It follows from \eqref{est:Katonew}, \eqref{est:2^j>M^{-1}}, \eqref{est:L2v0} and \eqref{est:small2} that for any $s\in[-\frac{M_3}{5},0]$,
                        \begin{align*}
                          &\sum_{\substack{\ell\geq -4}}\Bigl(\bigl\|  \dot \Delta_{0} (\dot{\Delta}_{\ell} v(s) \otimes \tilde{\Delta}_{\ell} v(s)) \bigr\|_{L^1(B(0,M_1/2))}+  \bigl\| \dot \Delta_{0} (\dot{\Delta}_{\ell} V_1(s) \otimes \tilde{\Delta}_{\ell} v(s)) \bigr\|_{L^1(B(0,M_1/2))}\Bigr)\\
                          &\lesssim \sum_{\substack{\ell\geq -4\\ 2^\ell \leq M_2^{\frac{1}{10}}}}\| \dot{\Delta}_{\ell} v(s)\|_{L^2(B(0,M_1))} \bigl( \| \tilde{\Delta}_{\ell} v(s)) \|_{L^2(B(0,M_1))}  +  \| \tilde{\Delta}_{\ell} V_1(s)) \|_{L^2(B(0,M_1))}\bigr)\\
&\quad + \sum_{\substack{\ell\geq -4\\ 2^\ell \geq M_2^{\frac{1}{10}}}}\| \dot{\Delta}_{\ell}
v(s)\|_{L^{1}(B(0,M_1))}\bigl( \| \tilde{\Delta}_{\ell} v(s)) \|_{L^\infty(B(0,M_1))} +  \| \tilde{\Delta}_{\ell} V_1(s)) \|_{L^\infty(B(0,M_1))}\bigr) \\
							&\quad+M_4^{-50c_p}\sum_{\substack{\ell\geq -4}}\left(\|\dot\Delta_\ell v\otimes\tilde\Delta_{\ell}v\|_{L^\infty ([-{M_3}/{5},0];L^1)}+\|\dot\Delta_\ell v\otimes\tilde\Delta_{\ell}V_1\|_{L^\infty ([-{M_3}/{5},0];L^1)}\right) \\
							&\lesssim M_1^{-2}M^{6c_p} + M_2^{-\frac{1}{2 }}M^{5c_p} + M_4^{-45c_p} \lesssim M_1^{-2}M^{6c_p}.
                        \end{align*}
                        Similar estimate can be obtained for $ \sum_{\substack{\ell\geq -4}}\bigl\| \dot \Delta_{0} (\dot{\Delta}_{\ell} V_1(s) \otimes \tilde{\Delta}_{\ell} v(s)) \bigr\|_{L^1(B(0,M_1/2))}$. Furthermore, by a similar argument as in \textbf{Step 2}, we deduce that
\begin{align*}
  \sum_{\substack{\ell\geq -4}}   \bigl\| \dot \Delta_{0} (\dot{\Delta}_{\ell} V_1(s) \otimes \tilde{\Delta}_{\ell} V_1(s)) \bigr\|_{L^1(B(0,M_1/2))}\lesssim M_4^{-45c_p}.
\end{align*}

                By summarizing the above estimates, we arrive at
                \begin{equation} \label{est:dydic0-3} \begin{aligned}
							\bigl\| \sum_{\substack{\ell\geq -4}} \dot \Delta_{0} (\dot{\Delta}_{\ell} u(s) \otimes \tilde{\Delta}_{\ell} u(s)) \bigr\|_{L^1(B(0,M_1/2))}\lesssim M_1^{-2}M^{6c_p}.
						\end{aligned} \end{equation}

					By	inserting \eqref{est:split-dyadic0}, \eqref{est:dydic0-1} and \eqref{est:dydic0-3} into \eqref{est:contra},
 we obtain $M_1^{-1} \lesssim M_1^{-2}M^{7c_p}$, which  contradicts the definition of  $M_1$.
 We thus conclude the proof of Proposition \ref{backprop}.
					\end{proof}

					\begin{proposition}[Iterated back propagation] \label{itbkpro}
						Let $u:[t_0-T,t_0] \times \R^3 \to \R^3$ be a classical solution of \eqref{prob:NS}
and satisfy \eqref{assumption}. Let $x_0\in\mathbb{R}^3$ and $j_0\in\mathbb{Z}$ be such that	
						\begin{equation*}
							|\dot{\Delta}_{j_0}u(t_0,x_0)|\geq M_1^{-1}2^{j_0},
						\end{equation*}
						then for any $T_1 \in [M_4 2^{-2j_0}, M_4^{-1}T]$, there exist
						\begin{equation*}
							(t_1,x_1)\in[t_0-T_1,t_0-M_3^{-1}T_1]\times B(x_0,CM_4^{c_p}T_1^{\frac{1}{2}}) \quad \text{and} \quad
							2^{j_1}\asymp M_3^{c_p}T_1^{-\frac{1}{2}}
						\end{equation*}
						such that
						\begin{equation*}
							|\dot{\Delta}_{j_1}u(t_1,x_1)|\geq M_1^{-1}2^{j_1}.
						\end{equation*}
					\end{proposition}
					\begin{proof}
						The proof is similar to that of \cite[Proposition 3.1 (v)]{Tao_20}, with a slight change regarding the pointwise estimate of $\partial_t\dot\Delta_ju$ in Lemma \ref{lem:ptw},  hence we omit the details here.
					\end{proof}
					
					Finally, we establish the annuli of regularity for $u$. 
					
					\begin{proposition}[Annuli of regularity] \label{anreg}
						Let $u:[t_0-T,t_0] \times \R^3 \to \R^3$ be a classical solution of \eqref{prob:NS} and
 satisfy \eqref{assumption}. Let $0<T'<\frac{T}{2}$, $x_0\in \mathbb{R}^3$, and $R_0>(T')^{\frac{1}{2}}$. Then there exists a scale
						\begin{equation*}
							R_0\leq R\leq \exp(M_6^{c_p})R_0
						\end{equation*}
such that any multi-index $\alpha \in \N^3$ with $|\al|=0,1$, we have
						\begin{equation}
\begin{split}
&\|\nabla^\al u\|_{L^\infty(\Omega)}\lesssim M_6^{-2}(T')^{-\frac{|\alpha|+1}{2}}\andf
						\|\nabla^\al \omega \|_{L^\infty(\Omega)}\lesssim M_6^{-2}(T')^{-\frac{|\alpha|+2}{2}} \\
&\with \Omega :=\{(t,x)\in[t_0-T',t_0]\times \mathbb{R}^3:R\leq |x-x_0|\leq M_6R \}.
\end{split}
						\end{equation}
					\end{proposition}
					\begin{proof}
						By a scaling argument, we may take $t_0=T'=1$, $T=2$ and $x_0=0$. We shall use the decomposition $u=v+V_1$ with $v$ solving \eqref{prob:NS-toy} and $V_1$ being defined by \eqref{def:Vi}.
						 Let us denote $$W:=\nabla\times v.$$
We observe that
						\begin{equation}\label{eq:cross}
							\nabla\times(f\cdot\nabla g)=f\cdot\nabla(\nabla\times g)+R,\quad \nabla\times (f\cdot\nabla f)=f\cdot\nabla(\nabla\times f)-(\nabla\times f)\cdot\nabla f,
						\end{equation}
						where $R$ consists of terms involving $\pm\nabla f\otimes\nabla g$. Then in view of \eqref{prob:NS-toy}, $W$ solves
						\begin{equation}\label{eq:W}
							\partial_tW-\Delta W=-v\cdot\nabla W+W\cdot\nabla v-V_1\cdot\nabla W-v\cdot\nabla(\nabla\times V_1)+F+G,
						\end{equation}
where $V_1,V_2,V_3$ are defined in  \eqref{def:Vi}, $F$ is the force term consisting of terms of the following form:
						\begin{equation*}
							\pm\partial_{i}V_3\otimes\partial_j V_1, \quad \pm V_1\otimes\partial_{ij}V_3,\quad \pm V_3\otimes \partial_{ij}V_2,\quad\pm\partial_jV_3\otimes \partial_{i}V_2, \quad 1 \leq i,j \leq 3,
						\end{equation*}
						and $G$ is the lower-order term consisting of terms of the following form:
						\begin{align*}
							\pm\partial_{i}v\otimes\partial_jV_1, \quad 1\leq i,j\leq 3.
						\end{align*}

						By virtue of \eqref{est:v-L2}, there exists $t_1\in(-{1}/{2},-\frac{1}{4})$ so that
						\begin{equation*}
							\int_{\mathbb{R}^3}|\nabla v(t_1,x)|^2 dx\lesssim M^{c_p}.
						\end{equation*}
Furthermore, since both $t_1$ and $-3/4$ are separated from the initial time, \eqref{S2eq12} gives
						\begin{align*}
							\|\nabla v(t_1)\|_{L^2}
							+\sum_{i=1}^m\sum_{|\al|\leq k_0+5}
							\Bigl(\|\nabla^\al u_{iL}(t_1)\|_{L^p}
							+\|\nabla^\al u_{iL}(-3/4)\|_{L^p}\Bigr)\lesssim M^{c_p},
						\end{align*}
						for some $k_0 \in\mathbb{N}$ large enough.  Then by the pigeonhole argument, we can find a scale $R$ in the range
						\begin{equation} \label{est:scale-R}
							M_6^{10000c_p}R_0\leq R\leq \exp(M_6^{c_p})R_0,
						\end{equation}
						with $R_0\geq 1$, such that
						\begin{equation}\label{sm1}
							\int_{\CC(M_6^{-10}R, M_6^{10}R)}\Bigl(|\nabla v(t_1,x)|^2
							+\sum_{i=1}^m\sum_{|\al|\leq k_0+5}
							\bigl(|\nabla^\al u_{iL}(t_1,x)|^p+|\nabla^\al u_{iL}(-3/4,x)|^p\bigr)\Bigr)dx
							\leq M_6^{-10000c_p}.
						\end{equation}
					Let us fix this $R$. Then we get, by using Sobolev inequalities, that
						\begin{equation*}
							\sup_{x\in \CC(M_6^{-9}R, M_6^9R)}\Bigl(\sum_{i=1}^m\sum_{|\al|\leq k_0}|\nabla^\al u_{iL}(t_1,x)|\Bigr)\lesssim M_6^{-1000c_p}.
						\end{equation*}
					We now claim that
						\begin{equation}\label{pwsm1}
							\sup_{t \in [t_1,1]}\sup_{x\in \CC(M_6^{-9}R, M_6^9R)}\Bigl(\sum_{i=1}^m\sum_{|\al|\leq k_0}|\nabla^\al u_{iL}(t,x)|\Bigr)\lesssim M_6^{-1000c_p}.
						\end{equation}
The proof will be postponed in Appendix \ref{Appb}.

				Next		we shall apply the energy method to equation \eqref{eq:W}. Let
us define the time-dependent radii, which is introduced in \cite{Tao13}.
\begin{equation}\label{Rrange}
						\begin{split}
							&R_-(t):=R_-+C_0\int_{t_1}^t(M_6+\|u(\tau)\|_{L^\infty})d\tau,\\
							&R_+(t):=R_+-C_0\int_{t_1}^t(M_6+\|u(\tau)\|_{L^\infty})d\tau, \\
&\text{with } R_-\in[M_6^{-8}R,2M_6^{-8}R] \andf  R_+\in [{M_6^{8}R}/{2},M_6^{8}R],
						\end{split}
					\end{equation}
						be scales to be determined later. Then
						for any $t\in[t_1,1]$, it follows from  the bounded total speed \eqref{est:L1Linfty}, \eqref{Rrange} and \eqref{est:scale-R}, that
						\begin{equation*}
							R_-(t)\in[M_6^{-8}R,3M_6^{-8}R],\quad R_+(t)\in [{M_6^{8}R}/{3},M_6^{8}R].
						\end{equation*}
						For $t\in[t_1,1]$, we define the local enstrophy
						\begin{equation*}
							E(t):=\frac{1}{2}\int_{\R^3}| W(t,x)|^2\eta(t,x)dx,
						\end{equation*}
						where $\eta$ is the time-varying cutoff function as in  \cite{Tao_20}, given by
						\begin{equation} \label{def:eta}
							\eta(t,x):=\max\left\{\min\bigl\{M_6,|x|-R_-(t), R_+(t)-|x| \bigr\},0 \right\}.
						\end{equation}
 It is easy to observe that  $\supp(\eta)\subset \mathcal{C}(R_-,R_+)=\{y:R_-\leq |x|\leq R_+ \}$, and $\eta=M_6$ in $\mathcal{C}(R_-+M_6,R_+-M_6)$.
 In particular, it follows from \eqref{sm1} that
						\begin{equation} \label{est:E(t_1)}
							E(t_1)\lesssim M_6^{-999c_p}.
						\end{equation}
We get, by taking $L^2$ inner product of  \eqref{eq:W} with $W$, that
						\begin{equation} \label{est:partial_tE-Yi}
							\partial_t E(t)+Y_1(t)+Y_2(t)\leq \sum_{i=3}^{8}Y_i(t),	
						\end{equation}
						where
						\begin{equation*}\label{def:Yi}
							\begin{aligned}
								&Y_1(t):=\int_{\mathbb{R}^3}|\nabla W(t,x)|^2\eta(t,x) dx,\\
								&Y_2(t):=-\frac{1}{2}\int_{\mathbb{R}^3}| W(t,x)|^2\partial_t \eta(t,x)dx,\\
								&Y_3(t):=\frac{1}{2}\int_{\mathbb{R}^3}| W(t,x)|^2\Delta\eta(t,x) dx,\\
								&Y_4(t):=\frac{1}{2}\Bigl|\int_{\mathbb{R}^3}| W(t,x)|^2u(t,x)\cdot\nabla \eta(t,x)dx\Bigr|,\\
								&Y_5(t):=\sum_{1\leq i,j\leq 3}\Bigl|\int_{\mathbb{R}^3}(v(t,x) W(t,x)\partial_{ij}V_1(t,x))\eta(t,x) dx\Bigr|,\\
								&Y_6(t):=\Bigl|\int_{\mathbb{R}^3}W(t,x)(W\cdot\nabla )v(t,x))\eta(t,x) dx\Bigr|,\\
								&Y_7(t):=\sum_{1\leq i,j\leq 3}\Bigl|\int_{\mathbb{R}^3}(\partial_i v(t,x) W(t,x)\partial_jV_1(t,x))\eta(t,x) dx\Bigr|,\\
								&Y_8(t):=\Bigl|\int_{\mathbb{R}^3}( W(t,x)\cdot F(t,x))\eta(t,x) dx\Bigr|.\\
							\end{aligned}
						\end{equation*}
					Below	we shall control $Y_i(t)$ for $3 \leq i \leq 8$ in terms of $Y_1(t),$ $ Y_2(t)$ and $E(t)$.  We first  deduce from \eqref{def:eta} that
						\begin{equation*}
							-\partial_t\eta=C_0(M_6+\|u\|_{L^\infty})|\nabla\eta|,
						\end{equation*}
						hence $Y_2\geq 0$ and we have
						\begin{equation*}
							Y_4(t)\leq C_0^{-1}Y_2(t).
						\end{equation*}
It follows from \eqref{est:v-L2}, \eqref{pwsm1} and H\"older's inequality that
						\begin{equation*}
							\begin{aligned}
								Y_7(t)&\lesssim E(t)+\sum_{1\leq i,j\leq 3}\int_{\mathbb{R}^3}\left|(\partial_iv\partial_j V_1)(t,x)\right|^2\eta(t,x)dx\\
								&\lesssim E(t)+M_6^{-100c_p}Y_9(t) \with Y_9(t):= \int_{\mathbb{R}^3}|\nabla v(t,x)|^2dx.
							\end{aligned}
						\end{equation*}
						Notice from   \eqref{est:v-L2} that
						\begin{equation}\label{est:Y9}
							\int_{t_1}^1Y_9(t)dt\leq M^{c_p}.
						\end{equation}
						Similarly, for $Y_5,Y_8$, we have
						\begin{equation*}
							Y_5(t)+Y_8(t)\leq \frac{1}{10}E(t)+M_6^{-10c_p}.
						\end{equation*}
						Next we deal with the heat flux term $Y_3(t)$ and the most ``non-linear" term $Y_6$. Since these two terms are of the same nature as the corresponding terms $Y_3(t)$ and $Y_6(t)$ in \cite[page 169]{Tao_20}, we can show that there exist $R_-$ and $R_+$ in the range \eqref{Rrange} such that
						\begin{equation} \label{est:Y3-1}
							\int_{t_1}^1Y_3(t)dt\lesssim M_6^{-10c_p},
						\end{equation}
						and
						\begin{equation} \label{est:Y6-1}
							Y_6(t)\leq \frac{1}{2}Y_1(t)+C\Bigl(\frac{1}{2C_0}Y_2(t)+M_6^{-2}+E(t)+Y_1(t)(E(t)^{\frac{1}{2}}+E(t)^2)\Bigr).
						\end{equation}
						Note that both \eqref{est:Y3-1} and \eqref{est:Y6-1} can be obtained by using the same  proof as \cite[estimate (3.44) and estimate at the end of page 173]{Tao_20}, hence we omit their proof here. We emphasize that the proof in \cite{Tao_20} only requires the $L_T^\infty ([t_0-\frac{T}{2},t_0];L^2)\cap L_T^2([t_0-\frac{T}{2},t_0];\dot H^1)$ bound of the non-linear term, which is also valid in our case due to \eqref{est:LinftyL2v}.

			By substituting the  above estimates for $Y_i$, $i=3,\ldots,9$, into \eqref{est:partial_tE-Yi} and taking $C_0\geq C$, we obtain
						\begin{equation*}
							\partial_tE(t)+\frac{1}{4}Y_1(t)\lesssim E(t)+Y_3(t)+M_6^{-2}+(E(t)^2+E(t)^{\frac{1}{2}})Y_1(t)+M_6^{-100c_p}Y_9(t),
						\end{equation*}
				from which, \eqref{est:E(t_1)}, \eqref{est:Y9} and \eqref{est:Y3-1},  we get, by a standard continuity argument, that
						\begin{equation} \label{est:E+Y1}
							E(t)+\int_{t_1}^1Y_1(\tau)d\tau\lesssim M_6^{-2}, \quad \forall t\in[t_1,1].
						\end{equation}
Since $W=\nabla\times v$, we get, by the standard elliptic estimate, that
						\begin{equation}\label{ellest}
							\|\nabla v(t)\|_{L^2(B)}^2\lesssim \|W(t)\|_{L^2(3B)}^2+r_B^{-2}\|v(t)\|_{L^2(3B)}^2,
						\end{equation}
						for any ball $B$ with radius $r_B$.

  By covering the annulus $\mathcal{C}(M_6^{-7}R,M_6^{7}R)$
 with  balls of radius $M_6^{2c_p}$, and using  \eqref{est:v-L2}, \eqref{est:E+Y1} and \eqref{ellest}, we find
						\begin{equation} \label{est:nablav-annulus}
							\int_{\mathcal{C}(M_6^{-7}R,M_6^7R)}|\nabla v(t,x)|^2dx\lesssim M_6^{-2}, \quad \forall t\in[t_1,1].
						\end{equation}
					Along the same line, we have
						\begin{equation*}
							\|\nabla^2 v(t)\|_{L^2(B)}^2\lesssim \|\nabla W(t)\|_{L^2(3B)}^2+r_B^{-2}\|\nabla v(t)\|_{L^2(3B)}^2,
						\end{equation*}
						and by the same argument as above, using \eqref{est:E+Y1} and \eqref{est:nablav-annulus},	we obtain
						\begin{equation*}
							\int_{t_1}^1 \int_{\mathcal{C}(M_6^{-6}R,M_6^6R)}|\nabla^2 v(t,x)|^2dxdt\lesssim M_6^{-2}.
						\end{equation*}
					Applying  the Gagliardo-Nirenberg inequality, we obtain
						\begin{equation*}
							\|v\|_{L^4((t_1,1);  L^\infty(\mathcal{C}(M_6^{-5}R,M_6^5R)))}\lesssim M_6^{-2},
						\end{equation*}
						which, together  with the fact: $u=v+V_1,$ and \eqref{pwsm1}, implies
						\begin{equation*}
							\|u\|_{L^4((t_1,1);  L^\infty(\mathcal{C}(M_6^{-5}R,M_6^5R)))}\lesssim M_6^{-2}.
						\end{equation*}
To proceed, we write
						\begin{equation*}
							u(t,x)=e^{(t-t_1)\Delta} u(t_1,x)-\int_{t_1}^t \int_{\R^3} \Gamma(t-s,x-y) (u\otimes u)(s,y)dyds, \quad x\in \mathcal{C}(M_6^{-4}R,M_6^4R),
						\end{equation*}
						where $\Gamma$ denotes the Oseen kernel which satisfies \eqref{S2eq13}.
						By splitting the integral domain into $|y|\leq M_6^{-5}R$, $M_6^{-5}R\leq |y|\leq M_6^5R$ and $|y|\geq M_6^{5}R,$ and using the decomposition
						$$u\otimes u=v\otimes v + v\otimes V_1 + V_1\otimes v + V_1\otimes V_1 \in L^{\frac{p}{2}}+L^1,
						$$
						(due to \eqref{est:v-L2} and \eqref{est:V1}), we find
						\begin{equation*}
							\|u\|_{L^8((t_1,1); L^\infty(\mathcal{C}(M_6^{-4}R,M_6^4R)))}\lesssim M_6^{-2}.
						\end{equation*}
						Iteratively, we can show that
						\begin{equation*}
							\|u\|_{L^\infty((t_1,1)\times \mathcal{C}(M_6^{-3}R,M_6^3R))}\lesssim M_6^{-2}.
						\end{equation*}
						By using a similar argument as in the proof of Proposition \ref{epochs}, we obtain
						\begin{equation*}
							\|\nabla u\|_{L^4((t_1,1);  L^\infty( \mathcal{C}(M_6^{-2}R,M_6^2R)))}\lesssim M_6^{-2c_p}
						\andf
							\|\nabla u\|_{L^\infty((t_1,1)\times \mathcal{C}(M_6^{-1}R,2M_6R))}\lesssim M_6^{-2}.
						\end{equation*}
						Similarly we obtain
						\begin{equation*}
							\|\nabla \omega \|_{L^\infty((t_1,1)\times \mathcal{C}(R,M_6R))}\lesssim M_6^{-2}.
						\end{equation*}
We thus complete the proof of Proposition \ref{anreg}.
					\end{proof}

					\section{Proofs of the main results}\label{sec:mainthm}

					We shall first use the estimates derived in Section \ref{Sect4} to obtain a key lemma of the paper.

					\begin{lemma}[Key lemma] \label{mt}
	Let $u$ be a smooth solution of \eqref{prob:NS} on $[t_0-T,t_0]$ and satisfy \eqref{assumption} and \eqref{assumption-Du}.
 Assume that there exist $x_0\in \mathbb{R}^3$ and $j_0 \in \Z$ such that
						\begin{equation}\label{hfsmall}
							|\dot \Delta_{j_0}u(t_0,x_0)|\geq M_1^{-1}2^{j_0}.
						\end{equation}
						Then one has
						\begin{equation}\label{inftybound}
							T\,2^{2j_0}\leq \exp\Bigl(A\exp\bigl(\exp(M_6^{c_p})\bigr)\Bigr).
						\end{equation}
					\end{lemma}
					\begin{proof}
						By time and space translation, we may assume without loss of generality that $(t_0,x_0)=(0,0)$. Let $T_1$ be an arbitrary scale in the range
						\begin{equation}\label{rangeT1}
							M_42^{-2j_0}\leq T_1\leq M_4^{-1}T.
						\end{equation}
					Then it follows from Proposition \ref{itbkpro} that there exist $
						(t_1,x_1)\in [-T_1,-M_3^{-1}T_1]\times B(0,M_4^{c_p}T_1^{\frac{1}{2}})$	
						and $2^{j_1}\asymp M_3^{c_p}T_1^{-\frac{1}{2}}$
						such that
						\begin{equation*}
							|\dot \Delta_{j_1}u(t_1,x_1)|\gtrsim M_1^{-1}2^{j_1},
						\end{equation*}
						from which, the Biot-Savart law, Lemma \ref{localY} and the estimate
						$\|\dot\Delta_j \omega (t_1)\|_{L^\infty(\mathbb{R}^3)}\lesssim M 2^{2j}$ (see \eqref{est:dyadic-vorticity-1}), we infer that
						\begin{equation*}
							|\dot{\Delta }_{j_1}\omega(t_1,x'_1)|\gtrsim M_1^{-1}2^{2j_1},
						\end{equation*}
						for some $x'_1\in B(x_1,O(M_1 2^{-j_1}))\subset B(0,M_4^{c_p}T_1^{\frac{1}{2}})$.
						Then we deduce from  Lemma \ref{lem:ptw} and Remark \ref{rmk:ptw} that for all $(t,x)\in [t_1,t_1+M_1^{-c_p}2^{-2j_1}]\times B(x'_1,M_1^{-c_p}2^{-j_1})$,
						\begin{equation*}
						\|\partial_t\dot{\Delta}_j \omega (t)\|_{L^\infty}\lesssim M^{c_p}(2^{4j}+2^{(2+\frac{3}{p})j}(t+T)^{-1+\frac{3}{2p}}),\quad \|\nabla\dot{\Delta}_j\omega(t)\|_{L^\infty}\lesssim M^{c_p}2^{3j}.
						\end{equation*}
						So for any $x\in B(x'_1,M_1^{-c_p}2^{-j_1}),$ we deduce that
						\begin{equation*}
						|\dot{\Delta}_{j_1}\omega(t_1,x)|\gtrsim M_1^{-1}2^{2j_1}-M_1^{-c_p}2^{-2j_1}2^{4j_1}\gtrsim M_1^{-1}2^{2j_1}.
						\end{equation*}
						Furthermore, for all $(t,x)\in [t_1,t_1+M_1^{-c_p}2^{-2j_1}]\times B(x'_1,M_1^{-c_p}2^{-j_1})$,
						\begin{align*}
						|\dot{\Delta}_{j_1}\omega(t,x)|\gtrsim M_1^{-1}2^{2j_1}-(M_1^{-c_p}2^{-2j_1})^{\frac{3}{2p}} 2^{(2+\frac{3}{p})j_1}-2^{-2j_1}M_1^{-c_p}2^{4j_1}\gtrsim M_1^{-1}2^{2j_1},
						\end{align*}
						which leads to\begin{equation}\label{ptww}
							|\dot{\Delta}_{j_1}\omega(t,x)|\gtrsim M_1^{-1}2^{2j_1}.
						\end{equation}
                        We observe  from
						Proposition \ref{epochs} that the interval
						\[
						J:=[-T_1,-M_3^{-c_p}T_1]\cap
						[t_1,t_1+M_1^{-c_p}2^{-2j_1}]
						\]
						has length comparable to $M_3^{-c_p}T_1$ (after increasing the generic
						constant $c_p$), and it contains a regularity epoch $I'$. Shrinking this epoch if
						necessary, we may arrange that $|I'|=cM_3^{-c_p}T_1$. For any multi-index
						$\alpha \in \N^3$ with $|\alpha|=0,1$,
						\begin{equation}\label{est:I'w}
							\begin{aligned}
								&|\nabla^\alpha u(t,x)|\lesssim M_3^{c_p}T_1^{-\frac{1+|\alpha|}{2}}\quad  \forall (t,x) \in I' \times \R^3, \\
								&|\nabla^\alpha \omega(t,x)|\lesssim M_3^{c_p}T_1^{-\frac{2+|\alpha|}{2}}\quad  \forall (t,x) \in I' \times \R^3.
							\end{aligned}
						\end{equation}	
           For any positive constant $C_0$, by taking a subset $I''$ of $I'$
           with  $|I''|=C_0^{-1}M_3^{-\frac{2c_p}{2}}T_1$ (which we still denote by $I'$), we have, for any multi-index $\alpha \in \N^3$ with $|\alpha|=0,1$,
						\begin{equation}\label{ptwu}
							|\nabla^\alpha u(t,x)|\lesssim (C_0|I'|)^{-\frac{1+|\alpha|}{2}}\quad \forall (t,x) \in I' \times \R^3.
						\end{equation}
						Then we deduce from  Lemma \ref{localY} and \eqref{ptwu} that
						\begin{equation}\label{wsmall}
							\int_{B(0,M_4^{c_p}T_1^{\frac{1}{2}})}|\omega(t,x)|^2dx\gtrsim M_3^{-c_p}T_1^{-\frac{1}{2}},\quad \forall t\in I'.
						\end{equation}
						
						To proceed, we shall use the quantitative propagation property of $\omega$ stated in
						Lemma \ref{tech} below. Its proof uses Proposition \ref{anreg} and the two Carleman
						estimates, and follows from \cite[estimates (5.10) and (5.17)]{Tao_20}, so we omit the details here.

						\begin{lemma}\label{tech}
							Let $u:[-T,0]\times\mathbb{R}^3\rightarrow \mathbb{R}^3$ be a classical solution of \eqref{prob:NS}
 and satisfy \eqref{assumption}. We assume that for any time scale $T_1$ in the range \eqref{rangeT1}, there exists a time interval $I'$ such that \eqref{est:I'w}, \eqref{ptwu}, \eqref{wsmall} hold. Let  $j_0$ be the index in \eqref{hfsmall}  and $T_2$ be a scale such that
							\begin{equation*}\label{rangeT2}
								M_4^{2}2^{-2j_0}\leq T_2\leq M_4^{-1}T.
							\end{equation*}
						Then there exists a scale $R$ in the range
							\begin{equation}\label{Rrange2}
								M_6T_2^{\frac{1}{2}}\leq R\leq \exp(M_6^{c_p})T_2^{\frac{1}{2}}
							\end{equation}
							such that in the cylindrical annulus $\Omega=\bigl\{(t,x) \in[-T_2,0]\times \mathcal{C}(R,M_6R) \bigr\}$, we have, for any multi-index $\alpha \in \N^3$ with $|\alpha|=0,1$,
							\begin{equation} \label{est:uomega-2}
								|\nabla^\alpha u(t,x)|\lesssim M_6^{-2}T_2^{-\frac{|\alpha|+1}{2}},\quad |\nabla^\alpha \omega(t,x)|\lesssim M_6^{-2}T_2^{-\frac{|\alpha|+2}{2}},
							\end{equation}
							and	\begin{equation} \label{est:omega0}
								\int_{\mathcal{C}(2R,\frac{M_6}{2}R)}|\omega(0,x)|^2dx\gtrsim \exp(-\exp(M_6^{c_p}))T_2^{-\frac{1}{2}}.
							\end{equation}
						\end{lemma}
						\begin{remark}
						We remark that in the proof of Lemma \ref{mt} only assumption \eqref{assumption} has been employed so far. To proceed further, we are going to make use of assumption \eqref{assumption-Du}.	
						\end{remark}
						\noindent \textbf{Continuation of the proof of Lemma \ref{mt}.} By the choice of $R$ obeying \eqref{Rrange2}, the volume of the annulus $\mathcal{C}(2R,\frac{M_6}{2}R)$ is controlled by $\exp (M_6^{c_p})T_2^{\frac{3}{2}}$, then it follows from  the pigeonhole  principle and \eqref{est:omega0} that
 there exists a $x_* \in \mathcal{C}(2R,\frac{M_6}{2}R)$ such that
						\begin{equation}\label{ptww0}
							|\omega(0,x_*)|\gtrsim \exp\bigl(-\exp(M_6^{c_p})\bigr)T_2^{-1}.
						\end{equation}
						We now use the signed assumption \eqref{assumption-Du}. Set
						\begin{equation}\label{def:r}
							r:=\varepsilon_0T_2^{\frac12},
							\qquad \varepsilon_0:=\exp\bigl(-\exp(M_6^{c_p})\bigr),
						\end{equation}
						and introduce
						\begin{equation}\label{def:sigma-signed}
							\sigma:=1-\frac3p\in(0,1),\qquad p':=\frac{p}{p-1}\in\Bigl(1,\frac32\Bigr).
						\end{equation}
						The remainder of the proof is divided into four steps.

\medskip
\noindent \textbf{Step 1: a signed concentration estimate.}

\begin{lemma}\label{lem:signed}
	Fix $\theta\in C_c^\infty(B(0,1))$ with $\theta\geq0$ and $\int_{\R^3}\theta\,dy=1$. Under
	the hypotheses of Lemma \ref{tech}, let $T_2$ be a scale in the admissible range, let $R$
	and $x_*\in\mathcal{C}(2R,\frac{M_6}{2}R)$ be as in \eqref{Rrange2} and \eqref{ptww0}, and
	let $r$ be given by \eqref{def:r}. Then there exists a unit vector $e\in\Sph^2$ such that
	the vector field
	\begin{equation}\label{def:Psi-signed}
		\Psi(z):=\nabla\theta\Bigl(\frac{x_*-z}{r}\Bigr)\times e
	\end{equation}
	satisfies
	\begin{equation}\label{est:signed}
		\Bigl|\int_{\R^3}u(0,z)\cdot\Psi(z)\,dz\Bigr|\ \geq\ \tfrac{c_0}{2}\,\varepsilon_0^{2}\,T_2^{-\frac12}r^{3},
	\end{equation}
	together with
	\begin{equation}\label{prop:Psi-signed}
		\supp\Psi\subset B(x_*,r),\qquad \|\Psi\|_{L^\infty}\leq\|\nabla\theta\|_{L^\infty},
		\qquad \int_{\R^3}\Psi(z)\,dz=0 .
	\end{equation}
\end{lemma}

\begin{proof}
	Let $e:=\omega(0,x_*)/|\omega(0,x_*)|$, so that $|\omega(0,x_*)|\geq c_0\varepsilon_0T_2^{-1}$
	by \eqref{ptww0}, where $\varepsilon_0$ denotes the \emph{same} quantity
	$\exp(-\exp(M_6^{c_p}))$ as in \eqref{def:r}, so that $r=\varepsilon_0T_2^{\frac12}$. Since
	$x_*\in\mathcal{C}(2R,\frac{M_6}{2}R)$ and $r\leq\varepsilon_0M_6^{-1}R\ll R$, the ball
	$B(x_*,r)$ is contained in $\mathcal{C}(R,M_6R)$, where \eqref{est:uomega-2} gives
	$\|\nabla\omega(0,\cdot)\|_{L^\infty}\leq CM_6^{-2}T_2^{-\frac32}$. Hence, for
	$|y|\leq1$, using that $CM_6^{-2}\leq\frac{c_0}{2}$ because $M_6$ is large,
	\begin{equation*}
		|\omega(0,x_*-ry)-\omega(0,x_*)|\leq CrM_6^{-2}T_2^{-\frac32}
		=C M_6^{-2}\varepsilon_0T_2^{-1}\leq\tfrac12|\omega(0,x_*)| ,
	\end{equation*}
	so that $\omega(0,x_*-ry)\cdot e\geq\frac12|\omega(0,x_*)|\geq\frac{c_0}{2}\varepsilon_0T_2^{-1}$
	for all $|y| \leq 1$. Therefore, with $\psi:=\theta\,e$ and since $\theta\geq0$,
	$\int\theta=1$,
	\begin{equation}\label{est:pairing-omega}
		\int_{\R^3}\omega(0,x_*-ry)\cdot\psi(y)\,dy
		=\int_{\R^3}\bigl(\omega(0,x_*-ry)\cdot e\bigr)\theta(y)\,dy\geq\frac{c_0}{2}\varepsilon_0T_2^{-1}.
	\end{equation}
	Writing $F(y):=u(0,x_*-ry)$, we have $\nabla_y\times F=-r\,\omega(0,x_*-r\,\cdot)$, whence,
	by integration by parts and $\nabla\times\psi=\nabla\theta\times e$,
	\begin{equation*}
		\int_{\R^3}\omega(0,x_*-ry)\cdot\psi(y)\,dy=-\frac1r\int_{\R^3} (\nabla\times F)\cdot\psi\,dy
		=-\frac1r\int_{\R^3} F(y)\cdot\bigl(\nabla\theta(y)\times e\bigr)\,dy .
	\end{equation*}
	Combining this identity with \eqref{est:pairing-omega} and substituting $z=x_*-ry$ yields
	\eqref{est:signed}. The properties in \eqref{prop:Psi-signed} are immediate, the last one
	because $\int_{\R^3}\nabla\theta\,dy=0$.
\end{proof}

\medskip
\noindent \textbf{Step 2: the dual test functions and their tails.} We set
\begin{equation}\label{def:g-phi}
	g:=|D|^{-\sigma}u(0,\cdot)\andf \phi:=r^{-2}|D|^{\sigma}\Psi ,
\end{equation}
so that $\|g\|_{L^p}\leq A^{\frac1p}$ by \eqref{assumption-Du}.

\begin{lemma}\label{lem:phi}
	With the notation of Lemma \ref{lem:signed}, there holds
	\begin{align}
		\label{est:phi-norm}&\|\phi\|_{L^{p'}(\R^3)}\leq C_p,\\
		\label{est:phi-point}&|\phi(z)|\leq C_p\,r^{2}|z-x_*|^{-4-\sigma}
		\qquad\text{for }|z-x_*|\geq 2r,\\
		\label{est:phi-tail}&\|\phi\|_{L^{p'}(\{|z-x_*|>\Lambda r\})}\leq C_p\Lambda^{-2}
		\qquad\text{for all }\Lambda\geq2,\\
		\label{est:phi-pairing}&|\langle g,\phi\rangle|\geq\tfrac{c_0}{2}\varepsilon_0^{3}.
	\end{align}
\end{lemma}

\begin{proof}
	We write $\Psi(z)=\Theta\bigl(\frac{z-x_*}{r}\bigr)$ with $\Theta(y):=\nabla\theta(-y)\times e$,
	and we note that $\Theta\in C_c^\infty(B(0,1))$, $\int_{\R^3}\Theta\,dy=0$, and that all the
	bounds below are uniform with respect to $e\in\Sph^2$. By homogeneity,
	\begin{equation}\label{scaling-phi}
		|D|^\sigma\Psi=r^{-\sigma}\bigl(|D|^\sigma\Theta\bigr)\Bigl(\frac{\cdot-x_*}{r}\Bigr),
		\qquad\text{so that}\qquad
		\bigl\||D|^\sigma\Psi\bigr\|_{L^{p'}}=r^{-\sigma+\frac3{p'}}\bigl\||D|^\sigma\Theta\bigr\|_{L^{p'}} .
	\end{equation}
	Since $\sigma=1-\frac3p$ and $\frac3{p'}=3-\frac3p$, we have the crucial identity
	\begin{equation}\label{key-exponent}
		-\sigma+\frac{3}{p'}=2 ,
	\end{equation}
	which gives $\|\phi\|_{L^{p'}}=\||D|^\sigma\Theta\|_{L^{p'}}$. For $\sigma\in(0,1)$ we have
	\begin{equation}\label{frac-lap}
		|D|^\sigma f(z)=c_\sigma\,{\rm p.v.}\!\int_{\R^3}\frac{f(z)-f(w)}{|z-w|^{3+\sigma}}\,dw ,
	\end{equation}
	so that $|D|^\sigma\Theta$ is smooth and, for $|y|\geq2$, using $\Theta(y)=0$ and
	$\int_{\R^3}\Theta dz=0$,
	\begin{equation}\label{est:Theta-decay}
		\bigl||D|^\sigma\Theta(y)\bigr|
		=c_\sigma\Bigl|\int_{B(0,1)}\Theta(w)\Bigl[\frac{1}{|y-w|^{3+\sigma}}-\frac{1}{|y|^{3+\sigma}}\Bigr]dw\Bigr|
		\leq C\|\Theta\|_{L^1}|y|^{-4-\sigma} .
	\end{equation}
	The gain of one power in \eqref{est:Theta-decay} comes precisely from the cancellation
	$\int_{\R^3}\Theta dz=0$ and the mean value theorem applied to the kernel. Since
	$|D|^\sigma\Theta$ is also bounded on bounded sets, we have
	\[
		\bigl||D|^\sigma\Theta(y)\bigr|\lesssim\min\{1,|y|^{-4-\sigma}\}.
	\]
	Combining this estimate with the scaling identity \eqref{scaling-phi} gives
	\begin{equation}\label{eq5.21}
		\bigl||D|^{1-\frac3p}\Psi(z)\bigr|
		\lesssim r^{-1+\frac3p}
		\min\left\{1,\left(\frac{|z-x_*|}{r}\right)^{-5+\frac3p}\right\}.
	\end{equation}
	In particular $|D|^\sigma\Theta\in L^q(\R^3)$ for every $q>\frac{3}{4+\sigma}$, hence for
	$q=p'>1$, which proves \eqref{est:phi-norm}; and \eqref{est:phi-point} follows from
	\eqref{scaling-phi}, \eqref{est:Theta-decay} and \eqref{key-exponent}. Consequently, since
	$(4+\sigma)p'>4>3$,
	\begin{equation*}
		\|\phi\|_{L^{p'}(\{|z-x_*|>\Lambda r\})}
		\leq C_p\,r^{2}\Bigl(\int_{\Lambda r}^{\infty}s^{-(4+\sigma)p'+2}\,ds\Bigr)^{\frac1{p'}}
		= C_p\,r^{2}(\Lambda r)^{\frac{3}{p'}-4-\sigma}=C_p\Lambda^{-2},
	\end{equation*}
	because $\frac{3}{p'}-4-\sigma=-2$ by \eqref{key-exponent}; this is \eqref{est:phi-tail}.
	Finally, $u(0,\cdot)=|D|^{\sigma}g$ and $|D|^\sigma$ is self-adjoint, so that
	\begin{equation}\label{est:duality-id}
		\langle u(0,\cdot),\Psi\rangle=\langle g,|D|^\sigma\Psi\rangle=r^2\langle g,\phi\rangle
	\end{equation}
	(see Remark \ref{rmk:duality} below), and \eqref{est:signed} together with
	$r=\varepsilon_0T_2^{\frac12}$ gives
	$|\langle g,\phi\rangle|\geq\frac{c_0}{2}\varepsilon_0^2T_2^{-\frac12}r=\frac{c_0}{2}\varepsilon_0^3$.
\end{proof}

\begin{remark}\label{rmk:duality}
	The identity \eqref{est:duality-id} is justified by
	\begin{equation*}
		\langle u(0,\cdot),\Psi\rangle=\sum_{j\in\Z}\langle\dot\Delta_ju(0,\cdot),\Psi\rangle
		=\sum_{j\in\Z}\langle\dot\Delta_jg,|D|^\sigma\Psi\rangle,
	\end{equation*}
	each term being finite by Plancherel's theorem. The series converges absolutely: by the
	quasi-orthogonality of the dyadic blocks,
	$\langle\dot\Delta_jg,|D|^\sigma\Psi\rangle
	=\langle\dot\Delta_jg,\tilde\Delta_j|D|^\sigma\Psi\rangle$ with
	$\tilde\Delta_j=\sum_{|j'-j|\leq1}\dot\Delta_{j'}$, so that
	\begin{equation*}
		\sum_{j\in\Z}\bigl|\langle\dot\Delta_jg,|D|^\sigma\Psi\rangle\bigr|
		\leq\sup_{j\in\Z}\|\dot\Delta_jg\|_{L^p}\sum_{j\in\Z}\|\tilde\Delta_j|D|^\sigma\Psi\|_{L^{p'}}
		\lesssim M\,\bigl\||D|^\sigma\Psi\bigr\|_{\dot B^0_{p',1}} ;
	\end{equation*}
	indeed $\|\dot\Delta_jg\|_{L^p}=\bigl\||D|^{-\sigma}\dot\Delta_ju\bigr\|_{L^p}
	\lesssim2^{-j\sigma}\|\dot\Delta_ju\|_{L^p}\leq CM$ by \eqref{assumption}, while
	$|D|^\sigma\Psi\in\dot B^0_{p',1}$, since
	$\|\dot\Delta_j|D|^\sigma\Psi\|_{L^{p'}}\lesssim2^{j\sigma}\|\dot\Delta_j\Psi\|_{L^{p'}}$
	and $\|\dot\Delta_j\Psi\|_{L^{p'}}\lesssim r^{\frac{3}{p'}}\min\{2^jr,(2^jr)^{-N}\}$, the
	low-frequency bound being a consequence of $\int\Psi=0$.
\end{remark}

\medskip
\noindent \textbf{Step 3: the multiscale family and its almost orthogonality.} Let
$\rho:=\exp(M_6^{C_p})$, with $C_p$ large enough that \eqref{cond:rho} below holds. We may
assume that $T\,2^{2j_0}\geq M_4^{6}\rho^{2}$, since otherwise \eqref{inftybound} is trivially
satisfied, and we define the sequence of scales
\begin{equation}\label{def:scales-signed}
	T_2^1:=M_4^22^{-2j_0},\quad T_2^{i+1}:=\rho^2T_2^i,\quad T_2^{K}\leq M_4^{-1}T,
	\qquad K=1+\Bigl[\frac{\ln\bigl(M_4^{-3}T\,2^{2j_0}\bigr)}{2\ln\rho}\Bigr]
	\asymp\frac{\ln(2^{2j_0}T)}{\ln\rho} ,
\end{equation}
so that in particular all the scales $T_2^i$, $1\leq i\leq K$, belong to the admissible range
of Lemma \ref{tech}.
For each $i\in\{1,\dots,K\}$ we apply Lemma \ref{tech} and Lemma \ref{lem:signed} with
$T_2=T_2^i$; this provides $R_i$, $x_*^i$, $r_i:=\varepsilon_0(T_2^i)^{\frac12}$, $\Psi_i$ and,
through \eqref{def:g-phi}, $\phi_i$, such that
\begin{equation}\label{est:geometry-signed}
	R_i\in\bigl[M_6(T_2^i)^{\frac12},\exp({M_6^{c_p}})(T_2^i)^{\frac12}\bigr],\qquad
	x_*^i\in\mathcal{C}\Bigl(2R_i,\frac{M_6}{2}R_i\Bigr),
\end{equation}
and such that \eqref{est:phi-norm}--\eqref{est:phi-pairing} hold for $(\phi_i,x_*^i,r_i)$. We
set
\begin{equation}\label{def:ai}
	a_i:=2M_6(T_2^i)^{\frac12},\qquad \kappa:=\tfrac14\exp({M_6^{c_p}}),
\end{equation}
so that, by \eqref{est:geometry-signed}, $|x_*^i|\in[a_i,\kappa a_i]$ and $a_{i+1}=\rho\,a_i$,
and we observe that
\begin{equation}\label{def:Lambda0}
	\frac{a_i}{2}=M_6(T_2^i)^{\frac12}=\Lambda_0\,r_i,\qquad
	\Lambda_0:=M_6\varepsilon_0^{-1}\gg1 .
\end{equation}
The constant $C_p$ in the definition of $\rho$ is chosen so large that
\begin{equation}\label{cond:rho}
	\rho\geq4\kappa=\exp({M_6^{c_p}} ).
\end{equation}
Finally, we partition $\R^3$ into the $K$ pairwise disjoint regions
\begin{equation}\label{def:Ai}
	\mathcal{A}_i:=\bigl\{z\in\R^3:\ \tau_i\leq|z|<\tau_{i+1}\bigr\},\qquad
	\tau_1:=0,\quad\tau_i:=\frac{a_i}{2}\ \ (2\leq i\leq K),\quad\tau_{K+1}:=+\infty .
\end{equation}

\begin{lemma}\label{lem:ortho}
	Assume \eqref{cond:rho}. Then $B(x_*^i,\Lambda_0r_i)\subset\mathcal{A}_i$ for every
	$1\leq i\leq K$, and
	\begin{equation}\label{est:crossterms}
		\sum_{j\neq i}\|\phi_j\|_{L^{p'}(\mathcal{A}_i)}\leq C_p\Lambda_0^{-2},
		\qquad 1\leq i\leq K .
	\end{equation}
	Consequently, for any choice of signs $\theta_1,\dots,\theta_K\in\{-1,+1\}$, the vector
	field $\Phi:=\sum_{i=1}^K\theta_i\phi_i$ satisfies
	\begin{equation}\label{est:Phi}
		\|\Phi\|_{L^{p'}(\R^3)}\leq C_pK^{\frac1{p'}} .
	\end{equation}
\end{lemma}

\begin{proof}
	\emph{(i) Localisation.} By \eqref{def:ai} and \eqref{cond:rho} we have
	$|x_*^i|\leq\kappa a_i\leq\frac{\rho a_i}{4}=\frac{a_{i+1}}{4}<\tau_{i+1}$ and
	$|x_*^i|\geq a_i>\tau_i$, hence $x_*^i\in\mathcal{A}_i$ and
	\begin{equation*}
		\dist(x_*^i,\mathcal{A}_i^c)\geq\min\Bigl\{|x_*^i|-\frac{a_i}{2},\
		\frac{a_{i+1}}{2}-|x_*^i|\Bigr\}\geq\frac{a_i}{2}=\Lambda_0r_i ,
	\end{equation*}
	where for $i=1$ (resp. $i=K$) only the second (resp. the first) term occurs. This is the
	first assertion.

	\emph{(ii) The scales $j<i$.} Since $|x_*^j|\leq\kappa a_j=\kappa\rho^{\,j-i}a_i\leq
	\frac{a_i}{4}$ by \eqref{cond:rho}, we get
	$\dist(x_*^j,\mathcal{A}_i)\geq\frac{a_i}{2}-\frac{a_i}{4}=\frac{a_i}{4}
	=\frac{\Lambda_0}{2}\rho^{\,i-j}r_j$, so that \eqref{est:phi-tail} yields
	\begin{equation*}
		\|\phi_j\|_{L^{p'}(\mathcal{A}_i)}\leq C_p\Lambda_0^{-2}\rho^{-2(i-j)} .
	\end{equation*}

	\emph{(iii) The scales $j>i$.} Here $\mathcal{A}_i\subset B\bigl(0,\frac{a_{i+1}}{2}\bigr)$
	and $\dist(x_*^j,\mathcal{A}_i)\geq a_j-\frac{a_{i+1}}{2}\geq\frac{a_j}{2}
	=\frac12\rho^{\,j-i}a_i$ for $j\geq i+1$. For $j=i+1$ we simply use \eqref{est:phi-tail}
	with $\Lambda=\Lambda_0$. For $j\geq i+2$ we use instead the pointwise bound
	\eqref{est:phi-point} together with $|\mathcal{A}_i|^{\frac1{p'}}\leq C(\rho a_i)^{\frac3{p'}}$
	and $r_j=\frac{\varepsilon_0}{2M_6}\rho^{\,j-i}a_i$, which gives
	\begin{align*}
		\|\phi_j\|_{L^{p'}(\mathcal{A}_i)}
		&\leq C_p\,r_j^2\bigl(\tfrac12\rho^{\,j-i}a_i\bigr)^{-4-\sigma}(\rho a_i)^{\frac3{p'}}\\
		&\leq C_p\Lambda_0^{-2}\,\rho^{\frac{3}{p'}}\,\rho^{-(j-i)(2+\sigma)}\,
		a_i^{\,2-4-\sigma+\frac3{p'}}
		\leq C_p\Lambda_0^{-2}\,\rho^{\,3-2(j-i)},
	\end{align*}
	where we used \eqref{key-exponent} again, so that the power of $a_i$ vanishes, together with
	$\frac3{p'}<3$ and $2+\sigma>2$. The right-hand side is summable in $j\geq i+2$, with sum
	bounded by $C_p\Lambda_0^{-2}\rho^{-1}$.

	Adding the three cases we obtain \eqref{est:crossterms}. Finally, since the regions
	$\mathcal{A}_i$ are pairwise disjoint and cover $\R^3$, and since
	$|\Phi|\leq\sum_{j}|\phi_j|$ pointwise, we deduce from \eqref{est:phi-norm} and
	\eqref{est:crossterms} that
	\begin{equation*}
		\|\Phi\|^{p'}_{L^{p'}(\R^3)}=\sum_{i=1}^K\|\Phi\|^{p'}_{L^{p'}(\mathcal{A}_i)}
		\leq\sum_{i=1}^K\Bigl(\|\phi_i\|_{L^{p'}}+\sum_{j\neq i}\|\phi_j\|_{L^{p'}(\mathcal{A}_i)}\Bigr)^{p'}
		\leq K\bigl(C_p+C_p\Lambda_0^{-2}\bigr)^{p'},
	\end{equation*}
	which is \eqref{est:Phi}.
\end{proof}

\medskip
\noindent \textbf{Step 4: conclusion.} Let $g=|D|^{-\sigma}u(0,\cdot)$ and let
$\phi_1,\dots,\phi_K$ be as in Step 3. By virtue of \eqref{est:phi-pairing} we may define
$\theta_i:={\rm sgn}\,\langle g,\phi_i\rangle\in\{-1,1\}$ and
$\Phi:=\sum_{i=1}^K\theta_i\phi_i$. Then, by \eqref{est:phi-pairing}, H\"older's inequality,
\eqref{assumption-Du} and \eqref{est:Phi},
\begin{equation} \label{est:K-epsi}
	\tfrac{c_0}{2}\varepsilon_0^{3}K\ \leq\ \sum_{i=1}^K|\langle g,\phi_i\rangle|
	=\langle g,\Phi\rangle\ \leq\ \|g\|_{L^p}\|\Phi\|_{L^{p'}}\ \leq\ C_pA^{\frac1p}K^{\frac1{p'}},
\end{equation}
whence $K^{\frac1p}\leq2c_0^{-1}C_pA^{\frac1p}\varepsilon_0^{-3}$ and therefore
\begin{equation}\label{est:K-signed}
	K\ \leq\ \bigl(2c_0^{-1}C_p\bigr)^{p}\,\varepsilon_0^{-3p}A\ \leq\ A\exp\bigl(\exp(M_6^{c_p})\bigr).
\end{equation}
Recall that the generic constant $c_p$ may vary from on appearance to other (we enlarge the constant $c_p$ to absorb the factor $3p$ and the constant coefficient).
Since $K\asymp\ln(2^{2j_0}T)/\ln\rho$ and $\ln\rho=M_6^{c_p}$ by \eqref{def:scales-signed}, we
conclude that
\begin{equation*}
	\ln\bigl(2^{2j_0}T\bigr)\leq M_6^{c_p}A\exp\bigl(\exp(M_6^{c_p})\bigr)
	\leq A\exp\bigl(\exp(M_6^{c_p})\bigr),
\end{equation*}
that is, $T\,2^{2j_0}\leq\exp\bigl(A\exp(\exp(M_6^{c_p}))\bigr)$, which is exactly
\eqref{inftybound}. Note that in the above estimate the generic constant $c_p$ may vary from on appearance to other (again we increase the constant $c_p$ in the right-hand side to absorb $M_6^{c_p}$). This completes the proof of Lemma \ref{mt}.
\end{proof}

\begin{remark}\label{rmk:signed-comments}
	The above argument uses neither the positivity of the kernel of $|D|^{-1+\frac3p}$ nor any
	pointwise information on $g$; the almost orthogonality involves only the explicit functions
	$\phi_i$, so that no a priori control of $\|g\|_{L^p}$ is needed at that stage. Three
	ingredients make it work: the identity \eqref{key-exponent}, which says that the
	normalisation $\|\phi_i\|_{L^{p'}}\asymp1$ is exactly the critical one; the cancellation
	$\int\Psi_i=0$, which is automatic since $\Psi_i$ arises as a curl and which produces the
	gain $\Lambda^{-2}$ in \eqref{est:phi-tail}; and the geometric separation
	\eqref{est:geometry-signed} of the concentration scales.
\end{remark}

Now we are in a position to present the proof of the main results in the paper.
					
					\begin{proof}[{\textbf{Proof of Theorem \ref{thm:main1}}}]
						By rescaling, we may assume that  $t=1$ so that $T \geq 1$. Then we get, by applying Lemma \ref{mt} in the contrapositive,
that for any $j$ such that $j\geq [j_*]$ with
						\begin{equation*}
							2^{2j_*}:=\exp\Bigl(A\exp\bigl(\exp(M_6^{c_p})\bigr)\Bigr).
						\end{equation*}
						there holds
						\begin{equation*}
							\|\dot{\Delta }_ju\|_{L^\infty([{1}/{2},1]\times\mathbb{R}^3)}\leq M_1^{-1}2^j.
						\end{equation*}
						Next we still use the decomposition \eqref{u=v+V1} that  $u=V_1+v$.
 Correspondingly, we decompose   the vorticity $\omega=\omega^{\lin}+W$ with $\omega^\lin:=\nabla \times V_1$ and $W= \nabla \times v$. Note that 	 for any $k \geq 0$ and $q\geq p$, we have
						\begin{equation*}
							\|\nabla^k V_1\|_{L^\infty([\frac{1}{2},1];L^q)}\lesssim M^{c_p}.
						\end{equation*}
						Let
						\begin{equation*}
							\sE(t):=\frac{1}{2}\int_{\mathbb{R}^3}|W(t,x)|^2\,dx,
						\end{equation*}
then we get, by taking $L^2$ inner product of  \eqref{eq:W} with $W,$ that
						\begin{equation}\label{S5eq1}
							\partial_t \sE(t)+\sY_1(t)\leq \sY_2(t)+\sY_3(t)+\sY_4(t)+\sY_5(t),
						\end{equation}
						where
						\begin{align*}
							&\sY_1(t):=\frac{1}{2}\int_{\mathbb{R}^3}|\nabla W(t,x)|^2dx,\\
							&\sY_2(t):=\Bigl|\int_{\mathbb{R}^3}W(t,x)\cdot(W\cdot\nabla v)(t,x)dx\Bigr|,\\
							&\sY_3(t):=\Bigl|\int_{\mathbb{R}^3}W(t,x)\cdot(v\cdot\nabla\omega^{\mathrm{lin}})(t,x)dx\Bigr|,\\
							&\sY_4(t):=\sum_{1\leq i,j\leq 3}\Bigl|\int_{\mathbb{R}^3}W(t,x)\partial_iV_1(t,x)\partial_jv(t,x)dx\Bigr|,\\
							&\sY_5(t):=\Bigl|\int_{\mathbb{R}^3}W(t,x)\cdot H(t,x)dx\Bigr|,
						\end{align*}
						with $H$ consisting of terms of the form:
						\begin{equation*}
							\pm\partial_{i}V_3\otimes\partial_j V_1, \quad \pm V_1\otimes\partial_{ij}V_3,\quad \pm V_3\otimes \partial_{ij}V_2,\quad\pm\partial_jV_3\otimes \partial_{i}V_2, \quad 1 \leq i,j \leq 3.
						\end{equation*}
					The terms $\sY_3,\sY_4,\sY_5$ can be estimated as follows:
						\begin{align*}
							&|\sY_3(t)|\lesssim \int_{\mathbb{R}^3}(| W(t,x)|^2+|v(t,x)|^2)|\nabla^2V_1|dx\lesssim M^{c_p} \sE(t)+M^{c_p},\\
							&|\sY_4(t)|\leq \int_{\mathbb{R}^3}|\nabla W(t,x)\|\nabla V_1(t,x)\|v(t,x)|dx+\int_{\mathbb{R}^3}| W(t,x)\|\nabla^2V_1(t,x)\|v(t,x)|dx\\
							&\quad\quad\quad \leq \frac{1}{2}\sY_1(t)+M^{c_p}\sE(t)+M^{c_p},\\
							&|\sY_5(t)|\lesssim \sE(t)+M^{c_p}.
						\end{align*}
						
						Next we handle the term $\sY_2$. Indeed by applying  Biot-Savart's law and Bony's decomposition, we obtain
						\begin{equation*}
							\sY_2(t)\lesssim\sum_{j_1\sim j_2\gtrsim j_3}\|\dot{\Delta}_{j_1}W(t)\|_{L^2}\|\dot{\Delta}_{j_2}W(t)\|_{L^2}\|\dot{\Delta}_{j_3}W(t)\|_{L^\infty}.
						\end{equation*}
			It follows from \eqref{assumption} that, for every $j\in\mathbb{Z}$,
			$\|\dot{\Delta}_j W(t)\|_{L^\infty}\lesssim M^{c_p}2^{2j}$. For the high
			frequencies we first use Lemma \ref{mt} to obtain
			\[
			\|\dot\Delta_j\omega(t)\|_{L^\infty}
			\lesssim 2^j\|\dot\Delta_ju(t)\|_{L^\infty}
			\lesssim M_1^{-1}2^{2j},\qquad j\geq[j_*].
			\]
			Moreover, \eqref{est:ulindecay} and $t\in[1/2,1]$ imply
			\[
			\|\dot\Delta_j\omega^{\rm lin}(t)\|_{L^\infty}
			\lesssim \sum_{k=1}^m2^{2j}e^{-c_{k-1}t2^{2j}}M^{c_p}
			\leq M_1^{-1}2^{2j},\qquad j\geq[j_*],
			\]
			after increasing the constant in the definition of $j_*$. Since
			$W=\omega-\omega^{\rm lin}$, we conclude, after harmlessly enlarging $M_1$, that
			$\|\dot\Delta_jW(t)\|_{L^\infty}\lesssim M_1^{-1}2^{2j}$ for $j\geq[j_*]$.
			Therefore, we have
						\begin{equation*}
							\sum_{j_3\lesssim j_2}\|\dot{\Delta}_{j_3}W(t)\|_{L^\infty}\lesssim M_1^{-1}2^{2j_2}+M^{c_p}2^{2j_*},
						\end{equation*}
						which implies
						\begin{equation*}
							\sY_2(t)\lesssim \sum_{j_1}\|\dot{\Delta}_{j_1}W(t)\|_{L^2}^2\bigl(M_1^{-1}2^{2j_1}+M^{c_p}2^{2j_*}\bigr).
						\end{equation*}
	By Plancherel's theorem, we get
						\begin{equation*}
							\sY_1(t)\asymp  \sum_{j_1}\|\dot{\Delta}_{j_1}W(t)\|_{L^2}^22^{2j_1} \quad \text{and} \quad
							\sE(t)\asymp \sum_{j_1}\|\dot{\Delta}_{j_1}W(t)\|_{L^2}^2.
						\end{equation*}
	We thus obtain
						\begin{equation*}
							\sY_2(t)\lesssim M_1^{-1}\sY_1(t)+M^{c_p}2^{2j_*}\sE(t).
						\end{equation*}

By substituting the above estimates into \eqref{S5eq1}, we derive
						\begin{equation}\label{S5eq2}
							\partial_t \sE(t)+\frac{1}{2} \sY_1(t)\lesssim M^{c_p}2^{2j_*}\sE(t) + M^{c_p} .
						\end{equation}
						Notice from \eqref{est:v-L2} that
						\begin{equation*}
							\int_{\frac{1}{2}}^{1}\sE(t)dt\lesssim M^{c_p}.
						\end{equation*}
						Therefore, for any time sub-interval of $[1/2,1]$ with length less than $M^{-c_p}2^{-2j_*}$, by the pigeonhole principle, there exists a point $t_*$ such that $\sE(t_*)\lesssim M^{2c_p}2^{2j_*}$. By the above estimates and Gronwall's inequality, we derive
						\begin{equation*}
							\sE(t_2)\lesssim \sE(t_1) +M^{c_p},
						\end{equation*}
						whenever $\frac{1}{2}\leq t_1\leq t_2\leq 1$ is such that $|t_1-t_2|\leq cM^{-c_p}2^{-2j_*}$ holds. Indeed, Gronwall's inequality applied to \eqref{S5eq2} gives
						\begin{equation*}
							\sE(t_2)\lesssim\Bigl(\sE(t_1)+M^{c_p}(t_2-t_1)\Bigr)
							\exp\Bigl(CM^{c_p}2^{2j_*}(t_2-t_1)\Bigr),
						\end{equation*}
						and the claimed estimate follows by choosing the absolute constant $c>0$ sufficiently small.
  Therefore, $\sE(t)\lesssim M^{2c_p}2^{2j_*}$ holds for any $\frac{3}{4}\leq t\leq 1$. Then we get, by
 integrating \eqref{S5eq2} over $[3/4,1],$ that
						\begin{equation*}
							\int_{\frac{3}{4}}^{1}\sY_1(t)dt\lesssim M^{3c_p}2^{4j_*}.
						\end{equation*}
						We now make the dependence of the final regularity bootstrap explicit. Put
						\[
						Q:=2^{2j_*}=\exp\Bigl(A\exp\bigl(\exp(M_6^{c_p})\bigr)\Bigr).
						\]
						In the proof of Proposition \ref{epochs}, replace the bounds for $\cE$ and its
						dissipation by the two estimates just obtained for $\sE$ and $\sY_1$. The same
						Gagliardo--Nirenberg and Duhamel steps, on time intervals of length
						$cM^{-c_p}Q^{-1}$, give
						\begin{equation}\label{est:final-W1inf}
						\|u\|_{L^\infty([7/8,1]\times\R^3)}\lesssim M^{c_p}Q,
						\qquad
						\|\nabla u\|_{L^\infty([7/8,1]\times\R^3)}\lesssim M^{c_p}Q^2.
						\end{equation}
						Indeed, the $L^4_tL^\infty_x$ step uses
						$\sup\sE\lesssim M^{2c_p}Q$ and
						$\int\sY_1\lesssim M^{3c_p}Q^2$; the Oseen-kernel estimate then upgrades
						this to $L^\infty_{t,x}$, and one spatial derivative is obtained by repeating the
						Duhamel estimate on the next half-interval. Increasing the generic exponent $c_p$ in
						$M_6^{c_p}$ absorbs the polynomial factors $M^{c_p}$ into $Q$. Thus
						\eqref{est:final-W1inf} proves \eqref{est:main-quadruple} for
						$|\alpha|=0,1$.

						For completeness, the higher derivatives follow without losing the stated exponent.
						After the preceding absorption, write the right-hand sides of
						\eqref{est:final-W1inf} as $Q$ and $Q^2$. Set $\lambda=Q^{-1}$ and, on the
						backward interval $[1-\lambda^2,1]$, define
						\[
						\widetilde u(s,y):=\lambda u(1-\lambda^2+\lambda^2s,\lambda y),
						\qquad 0\leq s\leq1.
						\]
						Then $\|\widetilde u\|_{L^\infty W^{1,\infty}}\lesssim1$. Lemma
						\ref{LeA2} gives $\|\nabla^\alpha\widetilde u(1)\|_{L^\infty}
						\lesssim_{|\alpha|}1$, and scaling back yields
						\[
						\|\nabla^\alpha u(1)\|_{L^\infty} \lesssim_{|\alpha|}Q^{1+|\alpha|}
						\leq Q^{2^{|\alpha|}},\qquad |\alpha|\geq2.
						\]
						This is precisely \eqref{est:main-quadruple} at the normalized time. Scaling back
						gives the factor $t^{(1+|\alpha|)/2}$.
						This completes the proof of Theorem \ref{thm:main1}.
					\end{proof}
					
					\begin{proof}[{\textbf{Proof of Theorem \ref{thm:main2}}}]
						By a time rescaling, we may assume that $T_*=1$. Set
						\[
							X(t):=\|u(t)\|_{\dot B_{p,\infty}^{-1+\frac3p}},
							\qquad
							Y(t):=\bigl\||D|^{-1+\frac3p}u(t)\bigr\|_{L^p},
							\qquad
							L(t):=|\ln(1-t)|.
						\]
						Fix a positive absolute constant $b$. We argue by contradiction and assume that
						\eqref{est:quantitative-Besovinfty} fails. Since the quantity in
						\eqref{est:quantitative-Besovinfty} is nonnegative, there exist
						$t_0\in(\frac12,1)$ and $C_0>0$ such that $L$ is positive and increasing on
						$[t_0,1)$ and
						\begin{equation}\label{est:expD}
							\exp\bigl(\exp(X(t)^{c_p})\bigr)Y(t)\leq C_0L(t)^b,
							\qquad t\in[t_0,1).
						\end{equation}

						By the standard small-data theory, there exists $\varepsilon_0>0$ such that
						a solution is global whenever, at some positive time,
						$\|u\|_{\dot B_{p,\infty}^{-1+\frac3p}}\leq\varepsilon_0$. Moreover,
						\[
							\|u\|_{\dot B_{p,\infty}^{-1+\frac3p}}
							\lesssim \bigl\||D|^{-1+\frac3p}u\bigr\|_{L^p},
						\]
						so smallness of the second norm also implies global existence. Since $u$
						blows up at time $1$, after decreasing $\varepsilon_0$ if necessary, we have
						\begin{equation}\label{est:norms-lower}
							X(t)\geq\varepsilon_0,\qquad Y(t)\geq\varepsilon_0,
							\qquad 0<t<1.
						\end{equation}

						For $t\in[t_0,1)$, define
						\[
							M_t:=\max\left\{2,\sup_{\tau\in[t_0,t]}X(\tau)^{c_p}\right\},
							\qquad
							A_t:=\max\left\{2,\sup_{\tau\in[t_0,t]}Y(\tau)^p\right\}.
						\]
						Using the monotonicity of $L$, \eqref{est:expD}, and
						\eqref{est:norms-lower}, we obtain, for $t$ sufficiently close to $1$,
						\begin{equation}\label{est:MtAt}
							M_t\lesssim \ln(b\ln L(t)),
							\qquad
							A_t\lesssim L(t)^{bp}.
						\end{equation}
						Indeed, the first estimate follows by taking two logarithms in
						$\exp(\exp(X(\tau)))\leq C_0\varepsilon_0^{-1}L(t)^b$, while the
						second follows directly from \eqref{est:expD}.

						From
						\eqref{est:MtAt}, we infer that
						\[
							\ln\left(A_t\exp\bigl(\exp(M_t)\bigr)\right)
							=\ln A_t+\exp(M_t)=bO(\ln L(t)).
						\]
						Consequently,
						\begin{equation}\label{est:small-exponent}
							A_t\exp\bigl(\exp(M_t^{c_p})\bigr)=L(t)^{bO(1)}
							\leq \frac{L(t)}{20}
						\end{equation}
						for all $t$ sufficiently close to $1$, if we take $b$ small enough.

						Apply Theorem \ref{thm:main1}, after translating the initial time to $t_0$,
						on the interval $[t_0,t]$. Since $t-t_0$ is bounded from below when $t$ is
						close to $1$, \eqref{est:small-exponent} gives
						\[
							\|u(t)\|_{L^\infty}
							\lesssim
							\exp\left(A_t\exp\bigl(\exp(M_t^{c_p})\bigr)\right)
							\leq \exp\left(\frac{L(t)}{20}\right)
							=(1-t)^{-\frac1{20}}.
						\]
						After enlarging the implicit constant, we may write
						\[
							\|u(t)\|_{L^\infty}\lesssim(1-t)^{-\frac1{10}},
							\qquad t\in\left(\frac12,1\right).
						\]
						Hence $u\in L^2((\frac12,1);L^\infty)$, contradicting the
						Prodi--Serrin--Ladyzhenskaya regularity criterion
						\cite{Pro_59,Ser_62,Lad_67}. This proves Theorem \ref{thm:main2}.
					\end{proof}

					\appendix
\section{The estimates of Kato norms}\label{Appa}
We start the proof of \eqref{S2eq12} by the following Lemma:

\begin{lemma}\label{estker}
Let $F\in C^\infty(\R^3\setminus \{0\})$ satisfy, for any multi-index $\beta \in \N^3$,
		\begin{equation}\label{conddec}
			|\nabla^\beta F(\xi)|\lesssim |\xi|^{\sigma-|\beta|}e^{-|\xi|^2},\quad \forall \xi \in \R^3 \setminus \{0\},
		\end{equation}
		for some $\sigma>-3$, and let $G(x)=\mathcal{F}^{-1}(F)(x)$. Then we have
		\begin{equation}\label{conddecay}
			|G(x)|\lesssim \frac{1}{(1+|x|)^{3+\sigma}}, \quad \forall x \in \R^3.
		\end{equation}
	\end{lemma}
	\begin{proof}
By assumption \eqref{conddec}, we see that $F\in L_{\xi}^1$, and by the definition of $G$, we have $\|G\|_{L^\infty}\leq \|F\|_{L^1}$. Then it suffices   to prove that $|G(x)|\lesssim |x|^{-3-\sigma}$ for any $x \in \R^3 \setminus \{0\}$.  In order to do so, for any $\delta>0,$
 we write
		\begin{equation*} G(x)=(2\pi)^{-\frac{3}{2}}\int_{\mathbb{R}^3}F(\xi)\chi_\delta(\xi)e^{ix\cdot\xi}d\xi+(2\pi)^{-\frac{3}{2}}\int_{\mathbb{R}^3}F(\xi)(1-\chi_\delta(\xi))e^{ix\cdot\xi}d\xi=:\text{I}+\text{II},
		\end{equation*}
		where $\chi_\delta$ is a smooth cut-off function so that $\chi_\delta(\xi)=1$ if $|\xi|\leq \delta$, $\chi_\delta(\xi)=0$ if $|\xi|\geq 2\delta$, and $\|\nabla^k\chi_\delta\|_{L^\infty}\lesssim \delta^{-k}$. Then we deduce from \eqref{conddec} that
		\begin{equation*}
			|\text{I}|\lesssim \int_{|\xi|\leq 2\delta}|\xi|^{\sigma}d\xi\lesssim \delta^{3+\sigma}.
		\end{equation*}
 To handle $\text{II}$, we define the operator $\hat{\nabla}:=\frac{x\cdot\nabla_\xi}{|x|^2}$. Note that $\hat{\nabla }e^{ix\cdot\xi}=i e^{ix\cdot\xi}$, hence we get, by using integration by parts for $l=[3+\sigma]+1$ times, that
		\begin{equation*}
			\begin{aligned}
				|\text{II}|&\lesssim \left|\int_{\mathbb{R}^3}F(\xi)(1-\chi_\delta)(\xi)\hat{\nabla}^{l}e^{ix\cdot\xi}d\xi\right|\lesssim
\sum_{l_1+l_2=l}\int_{\mathbb{R}^3}\bigl|\hat{\nabla}^{l_1}F(\xi)\hat{\nabla}^{l_2}(1-\chi_\delta(\xi))\bigr|\,d\xi\\
				&\lesssim |x|^{-l}\Bigl( \sum_{k=0}^{l-1}\delta^{k-l}\int_{\delta\leq|\xi|\leq 2\delta}|\xi|^{\sigma-k}d\xi+\int_{|\xi|\geq \delta}|\xi|^{\sigma-l}d\xi\Bigr)\lesssim |x|^{-l}\delta^{\sigma+3-l}.
			\end{aligned}
		\end{equation*}
As a result, it comes out
$$
|G(x)|\lesssim \delta^{3+\sigma}+|x|^{-l}\delta^{\sigma+3-l}. $$
Taking $\delta=|x|^{-1}$ leads to \eqref{conddecay}, and we complete the proof  Lemma \ref{estker}.
	\end{proof}

Let us now present the proof of \eqref{S2eq12}.

	\begin{proof}[Proof of \eqref{S2eq12}] We follow the ideas in \cite{CHN1}. In view of \eqref{prob:lin}, we write
				\begin{equation}\label{linfor}
\begin{split}
				u_{(k+1)L}(t,x)=&-\int_0^te^{(t-\tau)\Delta}\nabla\cdot\bP\Bigl(u_{kL}\otimes u_{kL}+\sum_{i=1}^{k-1}(u_{iL}\otimes u_{kL}+u_{kL}\otimes u_{iL}) \Bigr)(\tau,x)d\tau\\
=&-\int_0^t\Gamma(t-\tau)\ast\Bigl(u_{kL}\otimes u_{kL}+\sum_{i=1}^{k-1}\bigl(u_{iL}\otimes u_{kL}+u_{kL}\otimes u_{iL}\bigr) \Bigr)(\tau,x)d\tau,
\end{split}
				\end{equation}
where $\Gamma(t,x)$ denotes Oseen kernel. For any multi-index $\alpha \in \N^3$, there holds
\begin{equation}\label{S2eq13}
|\nabla^\al\Gamma(t,x)|\leq \frac{C_\al}{(\sqrt{t}+|x|)^{4+|\al|}},\quad \forall t>0, \, x \in \R^3.
\end{equation}
We notice that estimate \eqref{S2eq13} with $t=1$ is a direct corollary of Lemma  \ref{estker}  (see also Lemma 2.1 in \cite{CHN1}), then by rescaling, we obtain \eqref{S2eq13} for any $t>0$.  As a consequence, for any $p \in [1,\infty]$, we get
    \begin{equation} \label{est:Lp-Gamma}
        \|\nabla^{\alpha}\Gamma(t)\|_{L^p}\leq C_\alpha t^{-\frac{|\alpha|}{2}-2+\frac{3}{2p}}, \quad \forall t>0.
    \end{equation}

				We shall prove \eqref{S2eq12} by induction argument.  It is easy to observe from \eqref{Besov-characterization} and
\eqref{prob:u1L} that
 \eqref{S2eq12}  holds for $k=1$.
 Inductively, assume that \eqref{S2eq12} holds for every $1\leq k\leq\ell$,
 where $1\leq\ell\leq m-1$. We shall prove it for $k=\ell+1$.
 By the inductive assumption, for every $k\leq\ell$ and $|\al|\leq N_0$, we have
	\begin{equation}\label{indasp}
				\sup_{\substack{t\in (0,T]\\ q\in \bigl[\max\{1,\frac{p}{k}\},\infty\bigr]}}\Bigl(t^{\frac{1}{2}\left(|\al|+1-\frac{3}{q}\right)}\|\nabla^\al u_{kL}(t)\|_{L^q}\Bigr)\lesssim  \bigl(1+\|u_0\|_{\dot B_{p,\infty}^{-1+\frac{3}{p}}}\bigr)^{2^{k-1}}.
				\end{equation}
To construct $u_{(\ell+1)L}$, we first consider the case $|\alpha|=0$. By virtue of
\eqref{linfor} and \eqref{est:Lp-Gamma}, for
$r=\max\bigl\{1,\frac{p}{\ell+1}\bigr\}$, we infer
				\begin{align*}
				\|u_{(\ell+1)L}(t)\|_{L^{r}}
				\lesssim&\int_0^t\|\Gamma(t-\tau)\|_{L^1}\Bigl(\|u_{\ell L}\|_{L^{\frac{pr}{p-r}}}\| u_{\ell L}\|_{L^{p}}
+ \sum_{i=1}^{\ell-1} \|u_{iL}\|_{L^p}\| u_{\ell L}\|_{L^{\frac{pr}{p-r}}}\Bigr) (\tau)\,d\tau.
\end{align*}
Observing that $\frac{pr}{p-r}\geq\frac{p}{\ell}$ due to $r\geq\frac{p}{\ell+1},$ we get, by applying \eqref{indasp}, that
\begin{align*}
				\|u_{(\ell+1)L}(t)\|_{L^{r}}
				\lesssim&\int_0^t(t-\tau)^{-\frac{1}{2}}\tau^{-1+\frac{3}{2r}}d\tau(1+\|u_0\|_{\dot B_{p,\infty}^{-1+\frac{3}{p}}})^{2^{\ell}}
				\lesssim t^{-\frac{1}{2}(1-\frac{3}{r})}\bigl(1+\|u_0\|_{\dot B_{p,\infty}^{-1+\frac{3}{p}}}\bigr)^{2^{\ell}}.
				\end{align*}
Along the same line, for $\epsilon$ being sufficiently small, one has
\begin{align*}
				\|u_{(\ell+1)L}(t)\|_{L^{\infty}}
				&\lesssim\int_0^t\|\Gamma(t-\tau)\|_{L^{\frac{1}{1-\epsilon}}}\Bigl(\|u_{\ell L}\|_{L^\infty}\| u_{\ell L}\|_{L^{\frac{1}{\epsilon}}}+ \sum_{i=1}^{\ell-1} \|u_{iL}\|_{L^\infty}\|u_{\ell L}\|_{L^{\frac{1}{\epsilon}}}\Bigr) (\tau)d\tau\\
				&\lesssim\int_0^t(t-\tau)^{-\frac{1}{2}-\frac{3}{2}\epsilon}\tau^{-1+\frac{3}{2}\epsilon}d\tau \bigl(1+\|u_0\|_{\dot B_{p,\infty}^{-1+\frac{3}{p}}}\bigr)^{2^{\ell}}\\
				&\lesssim  t^{-\frac{1}{2}}\bigl(1+\|u_0\|_{\dot B_{p,\infty}^{-1+\frac{3}{p}}}\bigr)^{2^{\ell}}.
				\end{align*}
By interpolating the above two inequalities, we arrive at
				\begin{equation}\label{S2eq17}
				\sup_{t\in (0,T]}\Bigl(t^{\frac{1}{2}\left(1-\frac{3}{q}\right)}\|u_{(\ell+1)L}(t)\|_{L^q}\Bigr)\lesssim  \bigl(1+\|u_0\|_{\dot B_{p,\infty}^{-1+\frac{3}{p}}}\bigr)^{2^{\ell}},\quad \forall q\geq \max\bigl\{1,\frac{p}{\ell+1}\bigr\}.
				\end{equation}
Next we are going to establish estimates on the derivatives of $u_{(\ell+1)L}$. In view of \eqref{linfor},  for $r=\max\bigl\{1,\frac{p}{\ell+1}\bigr\}$ and
$|\al|\leq N_0,$ we infer
				\begin{align*}
				&\|\nabla^\al u_{(\ell+1)L}(t)\|_{L^{r}} \\
				&\lesssim \int_0^{\frac{t}{2}}\|\nabla^\al \Gamma(t-\tau)\|_{L^1}\Bigl(\|u_{\ell L}\|_{L^{\frac{pr}{p-r}}}\| u_{\ell L}\|_{L^{p}} + \sum_{i=1}^{\ell-1} \|u_{iL}\|_{L^p}\| u_{\ell L}\|_{L^{\frac{pr}{p-r}}}\Bigr) (\tau)\,d\tau\\
				& +\sum_{\beta+\gamma=\al}\int_{\frac{t}{2}}^t\|\Gamma(t-\tau)\|_{L^1}\Bigl(\|\nabla^{\beta}u_{\ell L}\|_{L^{\frac{pr}{p-r}}}\| \nabla^{\gamma}u_{\ell L}\|_{L^{p}} + \sum_{i=1}^{\ell-1} \|\nabla^{\beta}u_{iL}\|_{L^p} \|\nabla^{\gamma}u_{\ell L}\|_{L^{\frac{pr}{p-r}}}\Bigr)(\tau)d\tau.
				\end{align*}
In the second term on the right-hand side of the above estimate, $\beta, \gamma \in \N^3$ are multi-indices such that $\beta + \gamma=\alpha$.  On one hand, since $\frac{pr}{p-r}\geq \frac{p}{\ell}$, we get, by applying \eqref{est:Lp-Gamma} and \eqref{indasp}, that
				\begin{align*}
				&\|\nabla^\alpha u_{(\ell+1)L}(t)\|_{L^{r}}\\
				&\lesssim\Bigl(\int_0^{\frac{t}{2}}(t-\tau)^{-\frac{1+|\al|
}{2}}\tau^{-1+\frac{3}{2r}}\,d\tau +\int_{\frac{t}{2}}^{t}(t-\tau)^{-\frac{1}{2}}\tau^{-1+\frac{3}{2r}-\frac{|\al|}{2}}d\tau\Bigr)\bigl(1+\|u_0\|_{\dot B_{p,\infty}^{-1+\frac{3}{p}}}\bigr)^{2^{\ell}}\\
				&\lesssim  t^{-\frac{|\al|}{2}-\frac{1}{2}(1-\frac{3}{r})}\bigl(1+\|u_0\|_{\dot B_{p,\infty}^{-1+\frac{3}{p}}}\bigr)^{2^{\ell}}.
				\end{align*}
On the other hand, for some $\epsilon$ small enough, we obtain
				\begin{align*}
				&\|\nabla^\al u_{(\ell+1)L}(t)\|_{L^{\infty}}\\
				&\lesssim\int_0^{\frac{t}{2}}\|\nabla^\al \Gamma(t-\tau)\|_{L^\frac{1}{1-\epsilon}}\Bigl(\|u_{\ell L}\|_{L^{\frac{1}{\epsilon}}}\| u_{\ell L}\|_{L^{\infty}}+ \sum_{i=1}^{\ell-1} \|u_{iL}\|_{L^\infty}\| u_{\ell L}\|_{L^{\frac{1}{\epsilon}}}\Bigr)(\tau)\,d\tau\\
				& +\sum_{\beta+\gamma=\alpha}\int_{\frac{t}{2}}^t\|\Gamma(t-\tau)\|_{L^\frac{1}{1-\epsilon}}\Bigl(\|\nabla^{\beta}u_{\ell L}\|_{L^{\frac{1}{\epsilon}}} \| \nabla^{\gamma}u_{\ell L}\|_{L^{\infty}} + \sum_{i=1}^{\ell-1} \|\nabla^{\beta}u_{iL}\|_{L^\infty} \|\nabla^{\gamma}u_{\ell L}\|_{L^{\frac{1}{\epsilon}}}\Bigr)(\tau)\,d\tau\\ &\lesssim\Bigl(\int_0^{\frac{t}{2}}(t-\tau)^{-\frac{1+|\al|}{2}-\frac{3}{2}\epsilon}\tau^{-1+\frac{3}{2}\epsilon}d\tau
+\int_{\frac{t}{2}}^{t}(t-\tau)^{-\frac{1}{2}-\frac{3}{2}\epsilon}\tau^{-1+\frac{3}{2}\epsilon-\frac{|\al|}{2}}\,d\tau\Bigr)\bigl(1+\|u_0\|_{\dot B_{p,\infty}^{-1+\frac{3}{p}}}\bigr)^{2^{\ell}}\\
				&\lesssim t^{-\frac{1+|\alpha|}{2}}\bigl(1+\|u_0\|_{\dot B_{p,\infty}^{-1+\frac{3}{p}}}\bigr)^{2^{\ell}}.
				\end{align*}
	By interpolating the above two inequalities, for any $ q\geq \max\bigl\{1,\frac{p}{\ell+1}\bigr\},$ we achieve
				\begin{equation}\label{S2eq20}
				\sup_{\substack{t\in(0,T]\\ |\al|\leq N_0}}\Bigl(t^{\frac{|\al|}{2}+\frac{1}{2}\left(1-\frac{3}{q}\right)}\|\nabla^{\al}u_{(\ell+1)L}(t)\|_{L^q}\Bigr)\lesssim  \bigl(1+\|u_0\|_{\dot B_{p,\infty}^{-1+\frac{3}{p}}}\bigr)^{2^{\ell}}.
				\end{equation}
By combining\eqref{S2eq17} and \eqref{S2eq20}, we derive \eqref{indasp} for $k=\ell+1$. By induction, we obtain \eqref{indasp} for any $1 \leq k \leq m$, thus we complete  the proof of \eqref{S2eq12}.
			\end{proof}
Using a similar argument as in the proof of \eqref{S2eq12}, we also obtain
\begin{lemma}\label{LeA2} Let $u\in L^\infty([0,1];W^{1,\infty}(\mathbb{R}^3))$ be a solution of \eqref{prob:NS} with $T=1$. Then for any $\alpha \in \N^3$ such that $|\alpha| \geq 2$, there holds
	\begin{equation}
	\sup_{t \in [\frac{1}{2},1]}	\|\nabla^{\alpha} u(t)\|_{L^\infty}\lesssim_{|\alpha|} 1 + \|u\|_{L^\infty([0,1]; W^{1,\infty})}^{2^{|\alpha|}}.
	\end{equation}
\end{lemma}
\section{Proof of the propagation estimate}\label{Appb}

\begin{proof}[Proof of \eqref{pwsm1}]
	 The strengthened choice of $R$ in \eqref{sm1} includes the linear components at
	 time $-3/4$. Applying the Sobolev inequality on the slightly smaller annulus therefore gives
						\begin{equation}\label{4.6initial}
						\sup_{x\in \CC(M_6^{-10 }R, M_6^{10}R)}\Bigl(\sum_{i=1}^m\sum_{|\al|\leq k_0}|\nabla^\al u_{iL}(-\frac{3}{4},x)|\Bigr)\lesssim M_6^{-1000c_p}.
							\end{equation}
By induction, we are going to  prove that
						\begin{equation}\label{induc4.6}
						\sup_{t \in [t^i,1]}\sup_{x\in \CC(M_6^{-\gamma_i}R, M_6^{\gamma_i}R)}
\Bigl(\sum_{|\al|\leq k_0}|\nabla^\al u_{iL}(t,x)|\Bigr)\lesssim M_6^{-1000c_p},
						\end{equation}
						for some $-\frac{3}{4}<t^1<t^2<...<t^m<t_1$ and $10>\gamma_1>\gamma_2>...>\gamma_m>9$, which implies \eqref{pwsm1}. Here $m$ is the same as \eqref{decompv} and only depends on $p$.

First we are going to prove \eqref{induc4.6} for $i=1$. Let $K$ be the heat kernel in $\R^3$, namely
$K(t,x)=(4\pi t)^{-\frac{3}{2}}e^{-\frac{|x|^2}{4t}}$ for $t>0$ and $x \in \R^3$. We write
	\begin{equation} \label{est:u1L} \begin{aligned}
	u_{1L}(t,x)=&\int_{\CC(M_6^{-10}R, M_6^{10}R)}K(t+3/4,x-y)u_{1L}(-3/4,y)\,dy \\
	& +\int_{\R^3 \setminus \CC(M_6^{-10}R, M_6^{10}R)}K(t+3/4,x-y)u_{1L}(-3/4,y)\,dy.
	\end{aligned} \end{equation}
For some $t^1 \in (-\frac{3}{4},-\frac{1}{2})$ and for any $t \in [t^1,1]$, $x \in \R^3$, it follows from  \eqref{4.6initial} that
				\begin{equation*}
						\bigl|\int_{\CC(M_6^{-10}R, M_6^{10}R)}K(t+3/4,x-y)u_{1L}(-3/4,y)\,dy\bigr|\lesssim M_6^{-1000c_p}.
				\end{equation*}
We observe that for some $9<\gamma_1<10$ and for any $t^1<t<1$, if $ x \in \CC(M_6^{-\gamma_1}R, M_6^{\gamma_1}R)$
 and $y \in \R^3 \setminus \CC(M_6^{-10}R, M_6^{10}R)$,
\begin{equation*}
						K(t+3/4,x-y) \lesssim \bigl(t+3/4\bigr)^{-\frac{3}{2}}e^{-c{M_6^{-10}R^2}}e^{-c\frac{|x-y|^2}{t+\frac{3}{4}}}.
\end{equation*}
Combining the above estimate and the estimate $\|u_{1L}(-3/4)\|_{L^\infty}\lesssim M^{c_p}$, we deduce that
			\begin{equation*}
						\bigl|\int_{\R^3 \setminus \CC(M_6^{-10}R,M_6^{10}R)}K(t+3/4,x-y)u_{1L}(-3/4,y)\,dy\bigr|\lesssim M^{c_p}e^{-cM_6^{-10}R^2}\lesssim M_6^{-1000c_p},
			\end{equation*}
where in the second estimate we used the fact that the scale $R$ belongs to the range \eqref{est:scale-R}.
Substituting the above estimates into \eqref{est:u1L} leads to estimate \eqref{induc4.6} for $i=1$ and $|\alpha|=0$.\\
By the derivative estimate of heat kernel $K$, we can also derive the estimates for $\nabla^\al u_{1L}$ with  $1\leq |\al|\leq k_0$. Thus we have proved \eqref{induc4.6} for $i=1$.
						
Next let us	assume that \eqref{induc4.6} holds for all $i\leq j-1$ for some $2 \leq j \leq m$,  we are going  to prove \eqref{induc4.6} for $i=j$. Note that $u_{jL}$ satisfies the system \eqref{prob:lin}, which we rewrite as follows
\begin{equation}\label{Apal}
\left\{ 	\begin{aligned}
			&\partial_t u_{jL} -\Delta u_{jL}+\nabla P_{jL}=\nabla\cdot F_j, \\
			&\nabla \cdot u_{jL}=0, \quad  u_{jL}\vert_{t=-\frac{3}{4}}=u_{jL}\bigl(-\frac{3}{4}\bigr),
		\end{aligned} \right.
						\end{equation}
						where
						\begin{equation*}
						F_j=-u_{(j-1)L}\otimes u_{(j-1)L}-\sum_{i=1}^{j-2}u_{(j-1)L}\otimes u_{iL}-\sum_{i=1}^{j-2}u_{iL}\otimes u_{(j-1)L}.
						\end{equation*}
We notice that
\begin{equation} \label{est:Fj}\|F_j\|_{L^\infty([t^{j-1},1];  C^{k}(\mathbb{R}^3))}\lesssim M^{c_p}
\andf \|F_j\|_{L^\infty([t^{j-1},1];  C^{k}(\mathcal{C}(M_6^{-\gamma_{j-1}}R,M_6^{\gamma_{j-1}}R)))}\lesssim M_6^{-1000c_p}.
\end{equation}

We are going to establish the desired estimate of $u_{jL}$ for $|\alpha|=0$. In view of \eqref{Apal}, we write
						\begin{equation} \label{est:ujL} \begin{aligned}
						 u_{jL}(t,x)&=\int_{\mathbb{R}^3}K(t-t^{j-1},x-y)u_{jL}(t^{j-1},y)dy \\
						 &-\int_{t^{j-1}}^t\int_{\mathbb{R}^3}\Gamma(t-s,x-y)F_j(s,y)dyds
						 =:u_{jL}^1+u_{jL}^2,
						\end{aligned} \end{equation}
						where $K$ is the heat kernel and $\Gamma$ is the Oseen kernel which satisfies \eqref{S2eq13}. We now split
the integral domain in the formula of $u_{jL}^\ell, \ell=1,2$, into two parts: $\CC(M_6^{-\gamma_{j-1}}R, M_6^{\gamma_{j-1}}R)$ and
$\R^3 \setminus \CC(M_6^{-\gamma_{j-1}}R, M_6^{\gamma_{j-1}}R)$, and denote the corresponding integral as $u_{jL}^{\ell k}$ for $k=1,2.$

To estimate the term $u_{jL}^{12}$, we note that for some $t^{j} \in (t^{j-1},t]$ and for any $x \in \CC(M_6^{-\gamma_{j}}R, M_6^{\gamma_{j}}R)$ and $y \in \R^3 \setminus \CC(M_6^{-\gamma_{j-1}}R, M_6^{\gamma_{j-1}}R)$, we have $$K(t-t^{j-1},x-y)\lesssim \bigl(t^{j}-t^{j-1}\bigr)^{-\frac{3}{2}}e^{-c\frac{M_6^{-\gamma_{j-1}}R^2}{1-t^{j-1}}}e^{-c_j|x-y|^2}.$$
 Yet it follows from  \eqref{est:scale-R} that
  $$e^{-c\frac{M_6^{-\gamma_{j-1}}R^2}{1-t^{j-1}}}\lesssim M_6^{\gamma_{j-1}}R^{-2}\lesssim M_6^{-1000c_p}.$$
  As a result, for any $t \in [t^{j},1]$ and $x \in \CC(M_6^{-\gamma_j}R, M_6^{\gamma_j}R)$, we obtain, by using the inductive assumption, that
  \begin{align*}
  |u_{jL}^{12}(t,x)| \leq \int_{\R^3 \setminus \CC(M_6^{-\gamma_{j-1}}R, M_6^{\gamma_{j-1}}R) }K(t-t^{j-1},x-y)|u_{jL}(t^{j-1},y)|dy \lesssim M_6^{-1000c_p}.
  \end{align*}
						
For $u_{jL}^{22}$, by using estimates \eqref{est:Fj} for $F_j$ and estimate \eqref{S2eq13} for $\Gamma$, we deduce, for $x \in \CC(M_6^{-\gamma_j}R, M_6^{\gamma_j}R)$ and $1>t>t^{j}>t^{j-1}$, that
			\begin{align*}
				|u_{jL}^{22}(t,x)| &= 		\bigl|\int_{t^{j-1}}^t\int_{\R^3 \setminus \CC(M_6^{-\gamma_{j-1}}R, M_6^{\gamma_{j-1}}R)}\Gamma(t-s,x-y)F_j(s,y)dyds\bigr|\\
				&\lesssim M^{c_p}\int_{t^{j-1}}^t\int_{\R^3 \setminus \CC(M_6^{-\gamma_{j-1}}R, M_6^{\gamma_{j-1}}R)}(t-s)^{-2}\Bigl(\frac{1}{1+\frac{|x-y|}{\sqrt{t-s}}}\Bigr)^4dyds\\
				&\lesssim M^{c_p}\int_{t^{j-1}}^t\int_{\R^3 \setminus \CC(M_6^{-\gamma_{j-1}}R, M_6^{\gamma_{j-1}}R)}(t-s)^{-\frac{7}{4}}\Bigl(\frac{1}{1+\frac{|y|}{\sqrt{t-s}}}\Bigr)^{\frac{7}{2}}|x-y|^{-\frac{1}{2}}dyds\\
				&\lesssim M_6^{11}R^{-\frac{1}{2}}(t-t^{j})^{\frac{3}{4}}\lesssim M_6^{-1000c_p},
						\end{align*}
						where we used  \eqref{est:scale-R} in the last step.
						
			By using Young's inequality and the inductive assumption, we obtain, for any $t \in [t^{j},1]$ and  $x \in \CC(M_6^{-\gamma_{j}}R, M_6^{\gamma_{j}}R)$, that
				$$ |u_{jL}^{11}(t,x)| + |u_{jL}^{21}(t,x)| \lesssim M_6^{-1000c_p}.
				$$
Inserting the above estimates into \eqref{est:ujL} leads to the estimate \eqref{induc4.6} with $i=j$ and $|\alpha|=0$.
						
						The desired  estimates  of $\nabla^{\alpha}u_{jL}$, for $1 \leq |\alpha| \leq k_0$, can be
derived along the same line. Thus we complete the proof of \eqref{induc4.6}.
						\end{proof}

\section*{Acknowledgement}

Quoc-Hung Nguyen is supported by the CAS Project for Young Scientists in Basic Research, Grant No. YSBR-031, and the NSFC under Grant Nos. 1251101538 and 12595282. P. Zhang is partially supported by the National Key R$\&$D Program of China under Grant No. 2021YFA1000800 and by the National Natural Science Foundation of China under Grant Nos. 12421001, 12494542, and 12288201.

\section*{Declarations}

\subsection*{Conflict of interest} The authors declare that there are no conflicts of interest.

\subsection*{Data availability}
This article has no associated data.

				\end{document}